\documentclass{article}
\usepackage{amsmath,amsfonts,amssymb,amsthm,amsopn,amscd}
\usepackage{bbm}
\usepackage{hyperref}
\usepackage{color}
\usepackage{subfigure}
\usepackage[dvips]{graphicx}

\newcommand{\Si}{S}
\newcommand{\R}{\mathbbm{R}}

\newcommand{\SOS}[1]{{}}

\begin{document}

\title{Codimension Two Bifurcations and Rythms in Neural Mass Models}
\author{Jonathan Touboul \and Olivier Faugeras}

\maketitle
\hrule
\tableofcontents
\vspace{.7cm}
\hrule
\vspace{.7cm}

\section*{Abstract}
{\bf
Temporal lobe epilepsy is one of the most common chronic neurological disorder characterized by the occurrence of spontaneous recurrent seizures which can be observed at the level of populations through electroencephalogram (EEG) recordings. This paper summarizes some preliminary works aimed to understand from a theoretical viewpoint the occurrence of this type of seizures and the origin of the oscillatory activity in some classical cortical column models. We relate these rhythmic activities to the structure of the set of periodic orbits in the models, and therefore to their  bifurcations. We will be mainly interested Jansen and Rit model, and study the codimension one, two and a codimension three bifurcations of  equilibria and cycles of this model. We can therefore understand the effect of the different biological parameters of the system of the apparition of epileptiform activity and observe the emergence of alpha, delta and theta sleep waves in a certain range of parameter. We then present a very quick study of Wendling and Chauvel's model which takes into account $\rm{GABA}_{\rm{A}}$ inhibitory postsynaptic currents.
}

\section{Introduction}
Epilepsy is a common chronic neurological disorder affecting 1\% of the world population. It is characterized by the recurrence of seizures that strongly alter the patient quality of life. Epileptic seizures are transient manifestations of abnormal brain activity characterized by excessive and highly synchronous firing in networks of neurons distributing over focal or more extended cerebral regions. These networks are often referred to as ``epileptogenic networks'' (REF).

 It can manifest as an alteration in mental state, tonic or clonic movements, convulsions, and various other psychic symptoms. Different types of seizures are distinguished according to whether the source of the seizure within the brain is localized (partial or focal onset seizures) or distributed (generalized seizures) further divided according to the possible loss of conciousness and to the possible effect on the body.

Indeed, although the exact mechanisms leading to the various forms of epilepsy are still largely unknown, it is now commonly accepted that seizures result from a change in the excitability of single neurons or populations of neurons (see \cite{kandel-schwartz-etal:00}). Indeed, it is known, from epilepsy research, that the nature of the interactions between neurons and the properties of neurons themselves are altered in epileptogenic networks, although these alterations are not completely described and understood. Most of the studies suggest that both hyperexcitability and hypersynchronization processes occur in epileptogenic networks. However, although these two concepts are now widely accepted, most neuroscientists agree on the fact that they must be better defined and that associated mechanisms still have to be quantitatively characterized. Indeed, epileptogenic networks are complex dynamical systems in which increased excitability and synchronization may be caused by many different nonlinear factors, ranging from subcellular to neuronal population level (REF REF).

In this context, the objective of this paper is show how theoretical studies performed in physiologically-plausible computational models of neuronal assemblies (``neural mass models'') can allow us to establish some relationships between excitability-related parameters in models and some characteristic electrophysiological patterns typically observed in local field potentials (LFPs) or in the EEG recorded under normal or epileptic conditions.

More particularly, we report results from a mathematical study (stability, bifurcation) of two well-established biologically-inspired macroscopic models of EEG generation. The first model, refereed to as the Jansen and RitÕs model \cite{jansen-zouridakis-etal:93,jansen-rit:95}, is a minimal model of a neuronal population comprising two sub-populations of cells: pyramidal neurons and interneurons. Despite its relative simplicity, this model was shown to simulate signals with realistic temporal dynamics as encountered is real EEG (background activity, alpha activity, sporadic or rhythmic epileptic spikes) \SOS{(REF : Zetterberg, wendling 2000)}.  The second model, refereed to as the Wendling and ChauvelÕs model \cite{wendling-bellanger-etal:00,wendling-bartolomei-etal:02, wendling-hernandez-etal:05}, is an extended version of the former. It was intended to reproduce some architectonic features of the hippocampus, a sub-cortical structure very often involved in temporal lobe epilepsy (TLE). This model which accounts for three sub-populations of cells (pyramidal neurons and two-types of interneurons) was shown to simulate fast oscillations (beta, low gamma frequency band) as encountered at the onset of seizures in TLE. This fast onset activity was not represented by the Jansen and RitÕs model.

\SOS{FW:
Il faudrait ensuite dŽvelopper un peu les aspects ÇÊmathematical studyÊÈ, notamment annoncer les points originaux. A ce stade, je ne peux pas encore le faire car je nÕai pas creusŽ en dŽtails la partie math. 
Par contre, il faut mentionner quÕ une telle Žtude nÕa jamais ŽtŽ conduite dans ce type de modle.
}

Seizures are characterized by single electrical transients called spikes at a slow time resolution (hundreds of milliseconds to seconds), macroscopic events which have to be distinguished from spikes of single nerve cells, that last only 1 or 2ms.The EEG represents a set of potential recorded by multiple electrodes on the surface of the scalp. The electrical signal recorded is a measure of the extracellular current flow from the summated activity of many neurons, distorted by the filtering, attenuated by layers of tissues and bone, representing the activity of neural population of the brain. EEG patterns are characterized by the frequency and the amplitude of electrical activity in the range of 1?30Hz (sometimes larger) with amplitudes of 20 ? 100?V, and the frequencies observed are divided into four groups: delta (0.5-4 Hz), theta (4-7Hz), alpha (8-13 Hz) and beta (13-30 Hz). As neuronal aggregates become synchronized, the amplitude of the summated current becomes larger. 

\SOS{Bizarre de finir sur ce paragraphe. Je le mettrai avant la description des modles: on termine sur les aspects math puis sur la description du papier.}

\section{Neural mass models}
Jansen's neural mass model was first introduced by Lopes Da Silva and colleagues in 1974, and studied further by Van Rotterdam \emph{et al.} in 1982 \cite{lopes-da-silva-hoeks-etal:74, lopes-da-silva-rotterdam-etal:76, rotterdam-lopes-da-silva-etal:82}. These authors developed a biologically inspired mathematical framework to simulate spontaneous electrical activities of neuronal assemblies measured for instance by EEG, with a particular interest for alpha activity. In their model, three neuronal populations interact with both excitatory and inhibitory connections. Jansen \emph{et al.} \cite{jansen-zouridakis-etal:93,jansen-rit:95} discovered that besides alpha activity, this model was also able to simulate evoked potentials, i.e. EEG-like activity observed after a sensory stimulation. More recently, Wendling and colleagues used this model to synthesize activities very similar to those observed in epileptic patients \cite{wendling-bellanger-etal:00}, and David and Friston studied connectivity between cortical areas with a similar framework \cite{david-friston:03,david-cosmelli-etal:04}. Nevertheless, one of the main issue of Jansen's model is that it is not able to produce all the rhythms in play in epileptic activity or observed in EEG recordings. This limitation of the model lead Wendling, Chauvel and their colleagues \cite{wendling-hernandez-etal:05, wendling-chauvel:08} to develop an extended Jansen model to better reproduce hippocampus activity. This model is based on neurological studies \cite{jefferys-whittington:96} revealing that gamma-frequency oscillations were linked with the inhibitory interneurons in the hippocampal networks, and that two types of $\rm{GABA}_{\rm{A}}$ inhibitory postsynaptic currents may play a crucial role in the formation of the nested theta/gamma rhythms in the hippocampal pyramidal cells. For these reasons, we will also study the 5-populations extended Wendling and Chauvel's model.

\subsection{Jansen and Rit's model}\label{ssect:JR}

\paragraph{Description of the model}

\begin{figure}[!htb]
\begin{center}
\subfigure[Populations involved in Jansen's model]{\includegraphics[width=.5\textwidth]{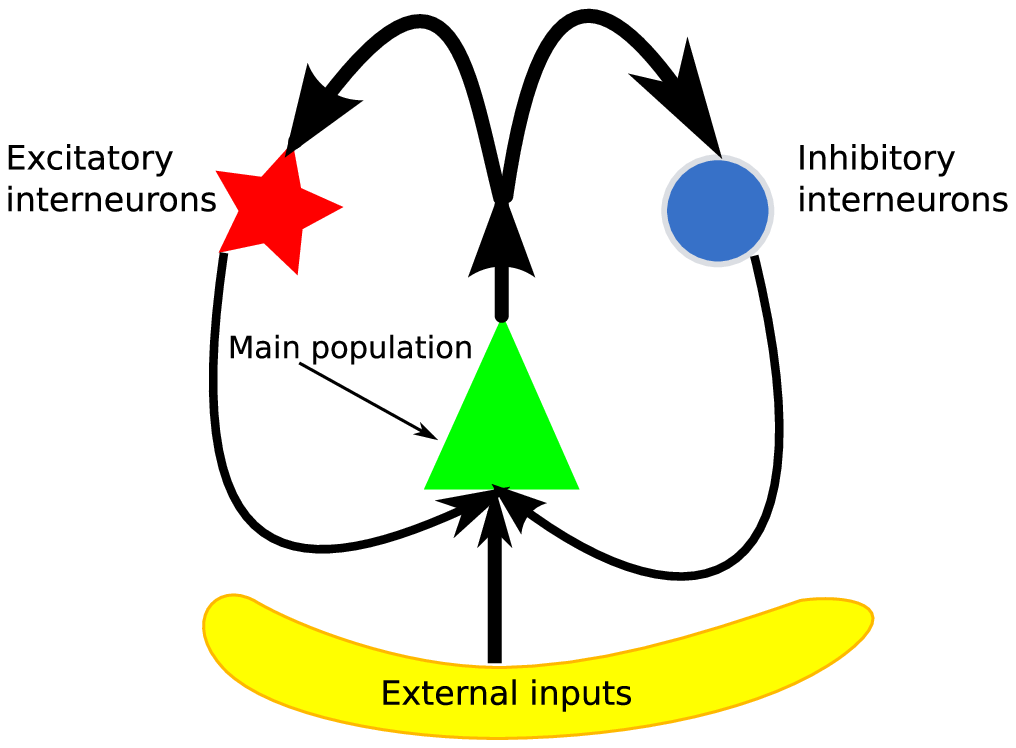}\label{sfig:JansenPops}}
\subfigure[Block diagram]{\includegraphics[width=.45\textwidth]{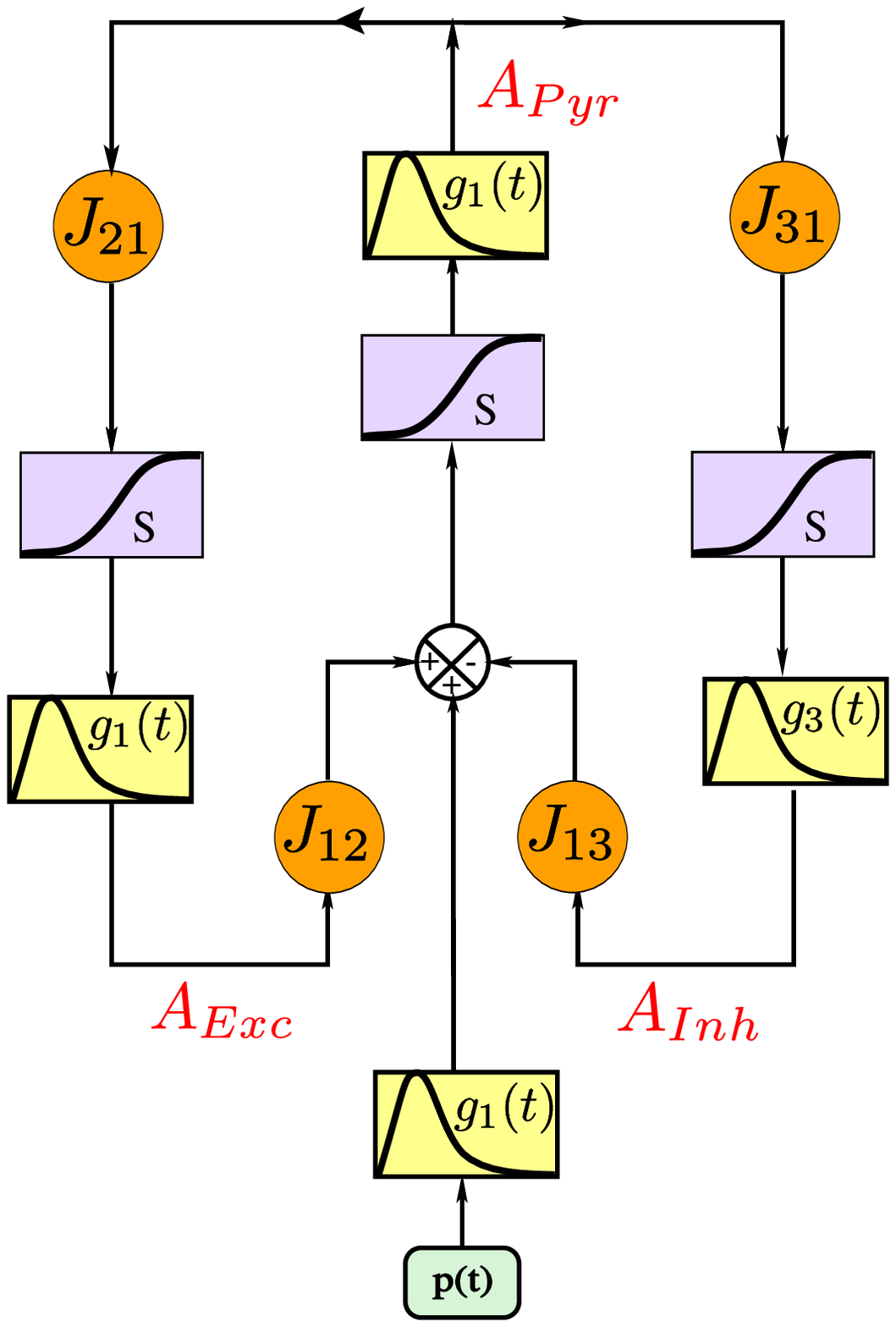}\label{sfig:Jansenblocks}}
\end{center}
\caption[Jansen and Rit's model]{(\ref{sfig:JansenPops})Neural mass model of a cortical unit: a pyramidal cells population interacts with an excitatory and an inhibitory population of interneurons. (\ref{sfig:Jansenblocks}) Block representation of a unit. $h$ boxes are synaptic transformations, $S$ boxes simulate cell bodies of neurons and is the sigmoidal transform of the membrane potential into an output firing rate. The constants $J_{i}$ account for the strength of the synaptic connections between populations. }
\label{jansen}
\end{figure}

\noindent The initial Jansen and Rit's model features a population of pyramidal neurons  (central part of
figure \ref{sfig:JansenPops}) that receive excitatory and inhibitory feedback from local inter-neurons 
and an excitatory input from neighboring cortical units and sub-cortical structures like the thalamus. Actually the excitatory feedback must be considered as coming from both local pyramidal neurons and genuine excitatory interneurons like spiny stellate cells. 

Figure \ref{sfig:Jansenblocks} is a block diagram representation of figure \ref{sfig:JansenPops} representing the mathematical operations performed inside such a cortical unit. The excitatory input is represented by an arbitrary average firing rate $p(t)$ which can be random (accounting for a non specific background activity) or deterministic, accounting for some specific activity in other cortical units. 

\paragraph{The postsynaptic system}
The functions $h_e(t)$ and $h_i(t)$ of figure \ref{sfig:Jansenblocks} are the average EPSP and IPSP. They convert the average input firing rate into an average excitatory or inhibitory post-synaptic potential. Following the works of van Rotterdam et al \cite{rotterdam-lopes-da-silva-etal:82} these transfer function are of type:
\[
h(t)=\left\{
\begin{array}{ll}
\alpha \beta te^{-\beta t} & t\geq0\\
0 & t<0
\end{array}
\right. .
\]
In other words, if $x(t)$ is the input to the system, its output
$y(t)$ is the convolution product $h\star x(t)$. The parameters $\alpha$ determines the maximal amplitude of the post-synaptic potentials and $\beta$ corresponds to the characteristic time of integration, which is mainly linked with the kinetics of synaptic transmission and with the averaged distributed delays in the dendritic tree.

These transfer functions are solutions of the second-order differential equation:
\begin{equation*}
\ddot{y}(t)=\alpha \beta x(t)-2\beta \dot{y}(t)-\beta^2y(t),
\end{equation*}
which can be conveniently rewritten as a system of two first-order equations

\begin{equation*}
\begin{cases}
  \dot{y}(t)&=z(t)\\
  \dot{z}(t)&=\alpha \beta x(t)-2\alpha z(t)-\alpha^2y(t)
\end{cases}
\end{equation*}

These two constants are different for the excitatory and inhibitory populations to fit the experimentally recorded EPSP and IPSP functions. The parameters $\alpha$ and $\beta$ have been adjusted by van Rotterdam \cite{rotterdam-lopes-da-silva-etal:82} to reproduce some basic properties of real post-synaptic potentials. Following their works, parameters can be set as follows:

\[
\begin{cases}
 \text{Excitatory population} & \alpha=:A=3.25\,mV \quad \beta=:a=100\,s^{-1}\\
 \text{Inhibitory population} & \alpha=:B=22\,mV \qquad \beta=:b=50\,s^{-1}
\end{cases}
\]

\paragraph{Firing rates}
The activity of the population is considered to be the related mean firing rate. This mean firing rate is modelled as a sigmoidal transformation of the average membrane potential (see, e.g. \cite{gerstner-kistler:02b}). The function $\mbox{Sigm}$ we chose in this model approximates the functions proposed by Freeman \cite{freeman:75}, and has the form:
\[
 \mbox{Sigm}(v)=\frac{\nu_{max}}{2}(1+\tanh{\frac{r}{2}(v-v_0)})=\frac{\nu_{max}}{1+e^{r(v_0-v)}},
\]
where $\nu_{max}$ is the maximum firing rate of the families of neurons,
$v_0$ is the value of the potential for which a $50\%$ firing rate is
achieved and $r$ is the slope of the sigmoid at $v_0$.

The excitability of cortical neurons can vary as a function of the action of several substances and $v_0$ could potentially take different values. In our model, we nevertheless consider a fixed value $v_0=6 mV$ as suggested by Jansen on the basis of experimental studies due to Freeman \cite{freeman:87}. The works of the latter also suggest that $\nu_{max}=5\,s^{-1}$ and $r=0.56\,mV^{-1}$, the values used by Jansen and Rit.

\paragraph{Interconnections} 
The three neural populations defined interact through excitatory and inhibitory synapses. The number of synapses established between two neuronal populations is denoted by  $J_{i}$ for $i=1\ldots 4$ as in diagram \ref{sfig:Jansenblocks}. They are considered to be constant and proportional to the average number of synapses $J$ between populations. We denote by $\alpha_i$ the related coefficient  (i.e. $J_i=\alpha_i J$). These coefficients can be seen as average probabilities of establishing connections between two populations. On the basis of several neuroanatomical studies (see for instance \cite{braitenberg-schuz:98}) these quantities have been estimated by counting synapses, and their numerical values is given in table \ref{tab:JansenNumerics}.

Note that we consider here constant synaptic weights. The variability in the connectivity weights can be taken into account (see e.g. \cite{faugeras-touboul-etal:09}), and the resulting equation is way more complex than the equation we deal with in the present article. 

\paragraph{Equations of the model}
Following Jansen and Rit's initial work, we consider the three variables $y_0$, $y_1$ and $y_2$ of figure \ref{sfig:Jansenblocks}. To write the system into a set of first-order ordinary differential equation we introduce the derivatives of these variables, $\dot{y}_0,\,\dot{y}_1,\,\dot{y}_2$, noted $y_3$, $y_4$ and $y_5$, respectively. We therefore obtain a system of 6 first-order differential equations that describes Jansen's neural mass model: 

\begin{equation}\label{eq:JROriginal}
\left\{
\begin{array}{ll}
\dot{y_0}(t)=y_3(t) &
\dot{y_3}(t)=Aa\mbox{Sigm}[y_1(t)-y_2(t)]-2ay_3(t)-a^2y_0(t)\\
\dot{y_1}(t)=y_4(t) &
\dot{y_4}(t)=Aa\{p(t)+J_2\mbox{Sigm}[J_1y_0(t)]\}-2ay_4(t)-a^2y_1(t)\\
\dot{y_2}(t)=y_5(t) &
\dot{y_5}(t)=BbJ_4\mbox{Sigm}[J_3y_0(t)]-2by_5(t)-b^2y_2(t).
\end{array}
\right.
\end{equation}

Our study will focus on the variable $y=y_1-y_2$ which models the membrane potential of the pyramidal cells since their electrical activity corresponds to the EEG signal: pyramidal neurons throw their apical dendrites to the superficial layers of the cortex where the post-synaptic potentials are summed, and therefore account for the essential part of the EEG activity \cite{kandel-schwartz-etal:00}. Table \ref{tab:JansenNumerics} summarizes the different numerical values of the original Jansen and Rit's model.

\begin{table}
 \begin{center}
  \begin{tabular}{lll}
  \hline
   Parameter & Interpretation & Value \\
   \hline
   A & Average excitatory synaptic gain & $3.25 mV$ \\
   B & Average inhibitory synaptic gain & $22 mV$ \\
   $1/a$ & Time constant of excitatory PSP & $10 ms$ \\
   $1/b$ & Time constant of inhibitory PSP & $20 ms$ \\
   $\alpha_1,\, \alpha_2$ &Average probability of synaptic contacts in the  & $\alpha_1=1$, $\alpha_2=0.8$\\
   & feedback excitatory loop  \\
   $\alpha_3,\, \alpha_4$ &Average probability of synaptic contacts in the & $\alpha_3=\alpha_4 = 0.25$ \\
   & slow feedback inhibitory loop \\
   $J$ & Average number of synapses between populations & 135\\
   $v_0, \, \nu_{max},\, r$ & Parameters of the Sigmoid $\Si$ & $v_0 = 6mV,\, \nu_{max}=5
s^{-1}$ \\
& & $r=0.56 mV^{-1}$\\
   \hline
  \end{tabular}
  \caption[Numerical values of JR model]{Numerical values used in Jansen's original model}
  \label{tab:JansenNumerics}
 \end{center}
\end{table}

\paragraph{Dimensionless reduction of Jansen's model}
We now change variables in order to reduce the number of parameters in the system. First of all, we chose the time scale of the excitatory population as our new unit time and define $\tau=at$. Therefore the new characteristic scale of the inhibitory variable will be given by the ratio $d=\frac{b}{a}$, $G$ the ratio of postsynaptic amplitudes $B/A$. The dimensionless sigmoidal transform we define is 
\[S(x) = \frac{1}{1+k_0e^{-x}}\]
where $k_0=e^{rv_0}$. We eventually introduce the dimensionless scaled input parameter $P=(r A) \frac{p}{a}$ and 
and the dimensionless scaled connectivity strength $j=(r A )\frac{\nu_{max}}{a}\,J$. The new state variables $(Y_0,X,Y_1,Y_2,Y_3,Y_4,Y_5)$ defined by
\begin{align*}
Y_0(\tau)  &= Jr y_0(\tau/a)\\
Y_i(\tau)  &= r y_i(\tau/a) \quad i=1,\,2 \\
X & = Y_1-Y_2.
\end{align*}
satisfy the equations
\begin{equation}\label{JR_Reduced}
\begin{cases}
\dot{Y_0} &= Y_3 \\
\dot{X}   &= Y_4-Y_5 \\
\dot{Y_2} &= Y_5 \\
\dot{Y_3} &= j\,S(X)-2\,Y_3-Y_0\\
\dot{Y_4} &=P+\alpha_2\,j\,S(\alpha_1\,Y_0)-2\,Y_4-(Y_2+X)\\
\dot{Y_5} &=d\,\alpha_4\,G\,j\,S(\alpha_3\,Y_0)-2\,d\,Y_5-d^2\,Y_2
\end{cases}
\end{equation}
They depend upon the nine dimensionless parameters
$\alpha_i,\,i=1,\cdots,4$, $k_0$, $d$, $G$, $j$, and $P$
as opposed to the 11 parameters of the original model (see table \ref{tab:JansenNumerics}).
The numerical values these parameters corresponding to table \ref{tab:JansenNumerics} are given by:

\begin{equation}\label{eq:Jreduced}
\left\{
\begin{array}{lcllcl}
G & = & \frac{B}{A} = 6.7692 & &  \\
d & = & \frac{b}{a} =0.5 & & \\
\alpha_1 & = & 1 & \alpha_2 & = &0.8  \\
\alpha_3 & = & 0.25 & \alpha_4 & = & 0.25  \\
\log(k_0) & = & r v_0 = 3.36 & &\\
j & = & (r A) \, \frac{\nu_{max}}{a}\,J = 12.285 &  &
\end{array}
\right.
\end{equation}

\subsection{Wendling and Chauvel's model of hyppocampus activity}\label{ssect:WC}
\begin{figure}
 \begin{center}
  \subfigure[Populations and connectivities]{\includegraphics[width=.5\textwidth]{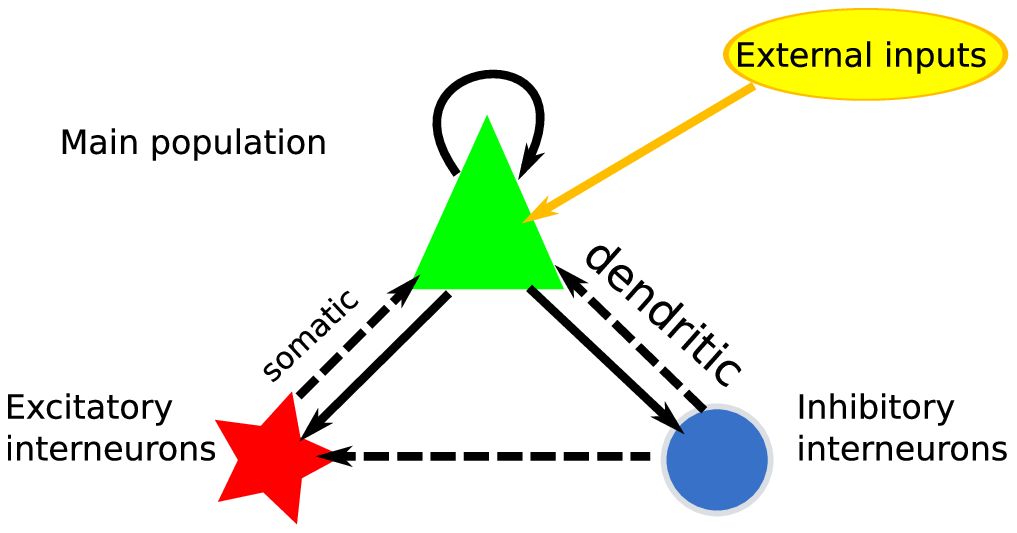}\label{sfig:WCPops}}
  \subfigure[Block diagram]{\includegraphics[width=.7\textwidth]{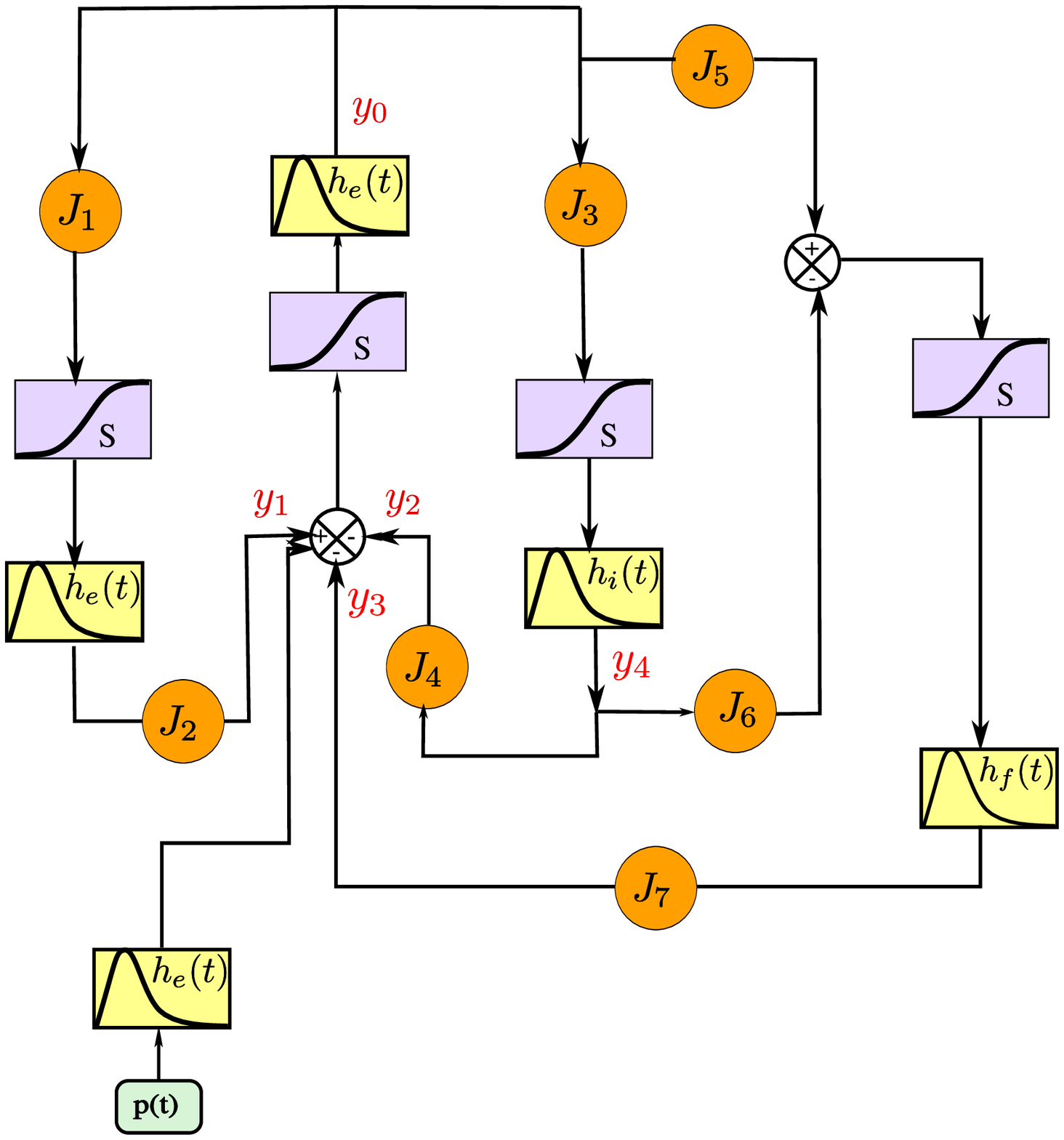}\label{sfig:WCBlock}}
  \caption[Wendling and Chauvel's model]{Neuronal population model based on the cellular organization of the hippocampus. \ref{sfig:WCPops}: Schematic representation of the model. The pyramidal cells population projects to and receives feedback excitatory and inhibitory interneurons, and has a recurrent excitation. \ref{sfig:WCBlock} : Related block diagram.}
  \label{fig:WendlingBlock}
 \end{center}
\end{figure}

One of the main drawbacks of Jansen's model is that it is unable to generate certain types of cortical activity, for instance seen in epileptic activity. As an example, it cannot reproduce the type of fast activity observed at the onset of seizures in limbic structures. In order to propose a better cortical mass model, Wendling, Chauvel and colleagues \cite{wendling-hernandez-etal:05,wendling-chauvel:08} revisited the organization of subsets of neurons and interneurons, focusing on the hippocampus activity. Based on these considerations, they proposed a new neural mass model whose parameters were estimated using real EEG signals. Their model is shown as a block diagram in figure Fig. \ref{fig:WendlingBlock}. The main difference with Jansen's initial model is the addition of somatic-projecting inhibitory neurons ($\rm{GABA}_{\rm{A,fast}}$ receptors).

In this block diagram representation, the PSP functions are given by:
\begin{equation*}
 \begin{cases}
  h_e(t) &= Aa t e^{-at} \\
  h_{i}(t) &= Bb t e^{-bt} \\
  h_{f}(t) &= Cc t e^{-ct} \\
 \end{cases}
\end{equation*}

The equations of the extended model hence read:

\begin{equation}\label{eq:WendlingOriginal}
 \begin{cases}
  \dot{y_0} &= y_5 \\
  \dot{y_5} &= A \, a \Si(y_1-y_2-y_3) - 2 a y_5 - a^2 y_0 \\
  \dot{y_1} &= y_6 \\
  \dot{y_6} &= A \, a \left\{ p(t) + J_2 \Si(J_1 y_0) \right\} - 2ay_6 - a^2 y_1 \\
  \dot{y_2} &= y_7 \\
  \dot{y_7} &= B\, b J_4 \Si(J_3 y_0) - 2 b y_7 - b^2 y_2 \\
  \dot{y_3} &= y_8 \\
  \dot{y_8} &= C \, c J_7 S(J_5 y_0 - J_6 y_4) - 2 c y_8 - c^2 y_3 \\
  \dot{y_4} &= y_9 \\
  \dot{y_9} &= B \, b S(J_3y_0) - 2 b y_9 - b^2 y_4
 \end{cases}
\end{equation}

\paragraph{Numerical values of the parameters}
This model has been fitted using SEEG data, and the authors obtained the values given in the table \ref{tab:paramsWendling}. 

\begin{table}
 \begin{tabular}{lll}
  \hline\\
  Parameter  & Interpretation & Value \\
  \hline\\
  A & Average excitatory synaptic gain & $3.25 mV$ \\
  B & Average inhibitory synaptic gain, slow dendritic inhibition loop  & $22 mV$ \\
  C & Average inhibitory synaptic gain, fast somatic inhibition loop  & $20 mV$ \\
  $\frac 1 a$ & Time constant of average excitatory postsynaptic potentials & $10ms$ \\
  $\frac 1 b$ & Time constant of average inhibitory postsynaptic potentials & $35ms$ \\
  $\frac 1 c$ & Time constant of the filter time delay & $5 ms$ \\
  $\alpha_5$, $\alpha_6$ & Average probability of synaptic contacts in the fast feedback inhibitory loop & $0.1$ \\
  $\alpha_7$ & Average probability of synaptic contacts between slow and fast inhibitory neuron & $0.8$ \\
  \hline
 \end{tabular}
 \caption[Parameters in WC's model]{Parameters interpretations and values of the extended model proposed by Wendling and Chauvel (see \cite{wendling-chauvel:08}). The parameters $\alpha_1, \ldots, \alpha_4, \, J,\, v_0,\, r$ and $\nu_{max}$ have the same interpretations and values as in Jansen's original model, see table \ref{tab:JansenNumerics}} 
\label{tab:paramsWendling}
\end{table}

\paragraph{Reduced Wendling-Chauvel's model}

First of all, one of the most straightforward reduction of the model consists in removing the variables $y_4$ and $y_9$ since they are deduced of $y_2$ and $y_7$ by the simple formulas: $y_2=J_4 y_4$ and $y_7 = J_4 y_9$. To reduce further the model we proceed in the same way as in Jansen and Rit's case. We make the system dimensionless by introducing the new time $\tau=at$ and proceeding to the change of variables: 

\begin{equation*}
 \begin{cases}
  Y_{i}(\tau) = Jr y_i (\tau/a) & i \in \{0,\, 4\}  \\
  Y_{i}(\tau) = r y_i (\tau/a) & i \in \{1,\, 2,\, 3\} \\
  Y_{i}(\tau) = \frac{Jr}{a} y_i (\tau/a) & i \in \{5,\, 9\}  \\
  Y_{i}(\tau) = \frac{r}{a} y_i (\tau/a) & i \in \{6,\,7,\,8 \}  \\
 \end{cases}
\end{equation*}

We denote by $X$ the interesting signal related to the EEG signal: $X = Y_1-Y_2-Y_3$ feeding the pyramidal population and by $Z = \alpha_5 Y_0 - \alpha_6 Y_4 = \alpha_5 Y_0 - \alpha_6/\alpha_4  Y_2 $ the action of the dendritic inhibitory interneuron population on the somatic inhibitory interneuron population. We consider now the variables $(Y_0, X, Z,Y_3, Y_5,Y_6,Y_7,Y_8)$.

These new variables satisfies the following set of differential equations, where we denote for the sake of compactness of notations by a dot the derivative with respect to $\tau$: 

\begin{equation}\label{eq:WendlingReduced}
  \begin{cases}
    \dot{Y_0} &= Y_5 \\ 
    \dot{X} &= Y_6 - Y_7 - Y_8 \\
    \dot{Y_2} &= Y_7 \\
    \dot{Y_3} &= Y_8 \\
    \dot{Z} &= \alpha_5 Y_5 - \alpha_6 Y_9\\
    \dot{Y_5} &= j S(X) - 2 Y_5 - Y_0 \\
    \dot{Y_6} &=  j \alpha_2 S(\alpha_1 Y_0) - 2Y_6 -  (X + Y_2 + Y_3) + P(t)\\
    \dot{Y_7} &= j d_1 G_1 \alpha_4 S(\alpha_3 Y_0) - 2 d_1 Y_7 - d_1^2 Y_2 \\
    \dot{Y_8} &= j d_2 G_2 \alpha_7 S(Z) - 2 d_2 Y_8 - d_2^2 Y_3 \\
    \dot{Y_9} &= j d_1 G_1 S(\alpha_3 Y_0) - 2 d_1 Y_9 - d_1^2 \frac{\alpha_5 Y_0 - Z}{\alpha_6}
  \end{cases}
\end{equation}
%

where $j$, $P$ and $S(x)$ are the same as the one introduced for Jansen and Rit's model. We used the notations:

\begin{equation}
 \left \{
 \begin{array}{ll}
  G_1 = \frac{B}{A} & G_2 = \frac{C}{A} \\
  d_1 = \frac{b}{a} & d_2 = \frac{c}{a} \\
 \end{array}
 \right .
\end{equation}

Using the numerical values in table \ref{tab:paramsWendling} we obtain

\begin{equation*}
 \left \{
 \begin{array}{ll}
  G_1 = 6.76923 & G_2 = 6.15385\\
  d_1 = 0.2857 & d_2 = 2 
 \end{array}
 \right .
\end{equation*}

In this article we are interested in the influence of the other parameters of the model together with the input on the fixed points. More precisely, we are interested in the codimension two bifurcations of this model with respect to $3$ pairs of parameters:
\begin{enumerate}
 \item The input $P$ and the coupling strength $j$.
 \item The input $P$ and the delays ratio $d$.
 \item The input $P$ and the PSP amplitude ratio $G$.
\end{enumerate}

\section{Influence of the total connectivity parameter in Jansen and Rit's model}
\label{sect:bifurcations}
We first study the dynamical properties of Jansen and Rit's model. We recall in the first subsection the main features described by Grimbert and Faugeras in \cite{grimbert-faugeras:06}, and extend their study to codimension two and three bifurcations.

\subsection{Fixed points and stability}
An interesting property of the system \eqref{JR_Reduced} is that the equilibria can be parametrized as a function of the state variable $X=Y_1-Y_2$:

\begin{align*}
Y_0&=j S(X)\\
Y_2&=\alpha_4\,\frac{G}{d}\,j S(\alpha_3 j S(X))\\
P&= X+j S(X)-\alpha_2 j S(\alpha_1 j S(X))
\end{align*}

The Jacobian matrix at the fixed point is also parametrized by $X$ and reads:
\[
J(X) = \left ( \begin{array}{cccccc}
                0 & 0 & 0 & 1 & 0 & 0  \\
                0 & 0 & 0 & 0 & 1 & -1 \\
                0 & 0 & 0 & 0 & 0 & 1  \\
                -1 & j\,S'(X) & 0 & -2 & 0 & 0 \\
                \alpha_1\alpha_2 j S'(j \alpha_1 S(X))  & -1 & -1 & 0 & -2 & 0 \\
                \alpha_3\alpha_4 j d  S'(j \alpha_3 S(X))  & 0 & -d^2 & 0 & 0 & -2d \\
               \end{array}
\right)
\]
Although all the dynamics can be parametrized with the variable $X$, because of the complexity of the sigmoidal function, the analytical bifurcation study is untractable, and one has to make use of a numerical software in order to solve the problem\footnote{Grimbert and Faugeras used the XPPAut software \cite{ermentrout:02}}. 

In the present case, almost all the calculations can be performed analytically in function of the variable $X$. For this reason, for computing accurately the bifurcations of equilibria, we developed our own software implemented using Maple\textregistered  ~in order to identify our codimension two bifurcations. For simple models, the programs we developed give closed-form expressions for the bifurcations points. For the continuation of periodic orbits, the analytical study is no more possible, and we used both the Matlab\textregistered ~ toolbox MatCont \cite{dhooge-govaerts-etal:03} and XPPAut continuation softwares to study the global bifurcations. The algorithms we use and the results we obtain are presented in appendix \ref{appendix:Numerbif}.

\subsection{Codimension 1 bifurcations}
The dependency of the dynamics to the input firing rate has already been studied Grimbert and Faugeras in \cite{grimbert-faugeras:06} when parameters are set as in table \ref{tab:JansenNumerics}. They describe a very rich bifurcation diagram, with the coexistence of two limit cycle, one coming from a Hopf bifurcation, and the other collapsing on the fixed points manifold, as we show in figure \ref{fig:Codim1}.
%

The system features two saddle-node and three Hopf bifurcations. One of the Hopf bifurcations is subcritical, the other two supercritical. The branch of unstable limit cycles originating from this point undergoes a fold bifurcation and connects to a family of stable limit cycles of large amplitude that eventually collide with a saddle-node bifurcation point and disappear via saddle-node homoclinic bifurcation. The periods of these cycles correspond to frequencies in the dimensioned model ranging from $0$ to $5$Hz which is consistent with the frequencies of recorded epileptic spikes. It corresponds to what was interpreted as epileptic oscillatory activity.

The other two Hopf bifurcation share the same family of periodic orbits. The period of these cycles is almost constant, ranging from $9$ to $9.6$ in the dimensionless model, which corresponds in the original model to frequencies in the alpha band. 

\begin{figure}
 \begin{center}
  \includegraphics[width=.5\textwidth]{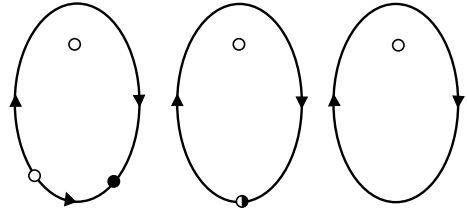}
 \end{center}
 \caption[SNIC bifurcation.]{Projection of a SNIC bifurcation in a subplane of the phase space. Three invariant circles are represented: Left. An heteroclinic orbit. Center. The two fixed points belonging to the invariant circle  merge into a saddle-node fixed point. The resulting homoclinic orbit has an infinite period. Right. The invariant circle turns into a periodic orbit as the fixed points disappear.  }
 \label{fig:SNIC}
\end{figure}

The system undergoes a saddle-node homoclinic bifurcation (SNIC) at the saddle-node bifurcation point, corresponding to a transition between a periodic orbit and an heteroclinic orbit (see figure \ref{fig:SNIC}),
and after this bifurcation, the system presents a family of heteroclinic orbits, which are not linked with oscillations, since the cycle contains a stable fixed point. 

It had been suggested by Grimbert and Faugeras \cite{grimbert-faugeras:06} that this picture was quite sensitive to changes in the parameters. They observed that varying any parameter by more than $5\%$ resulted in drastic changes in the bifurcation diagram and the behaviors. This is why we were interested in understanding better the appearance of these features and its sensitivity in function of different parameters. 

In this section, we will be particularly interested in the influence of the coupling strength $j$.

\subsection{Effect of the coupling strength and the input current}
Let us first study the bifurcations of the system with respect to the pair of parameters $(j,P)$. We first study the bifurcations of equilibria of the system, before studying global bifurcations of cycles and the resulting rhythms generated.

\begin{figure}
 \begin{center}
 \begin{minipage}{.6\textwidth}
  \subfigure[Full Diagram]{\includegraphics[width=\textwidth]{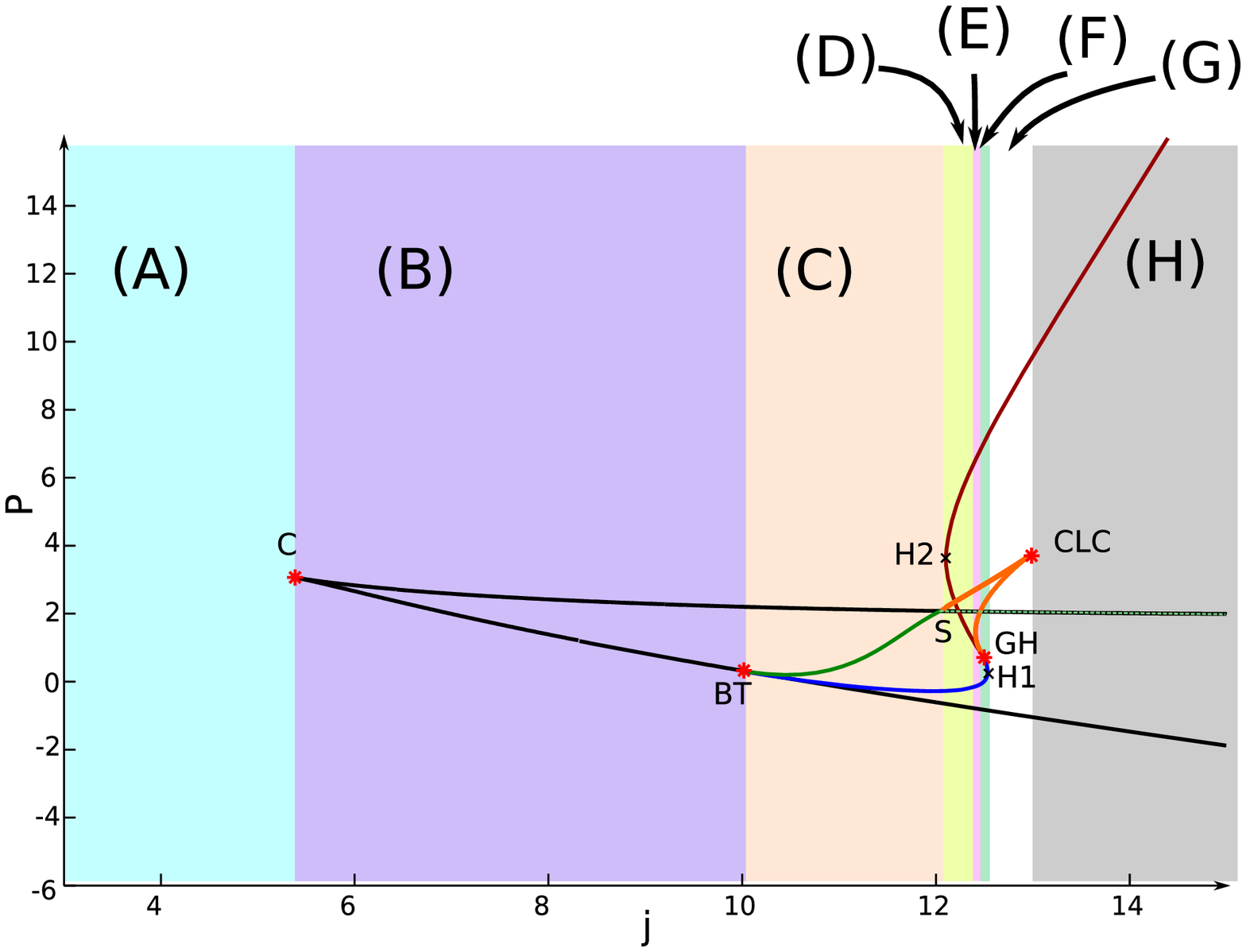}}
  \end{minipage}
  \begin{minipage}{.39\textwidth}
  \subfigure[Parameter zone of interest.]{\includegraphics[width=\textwidth]{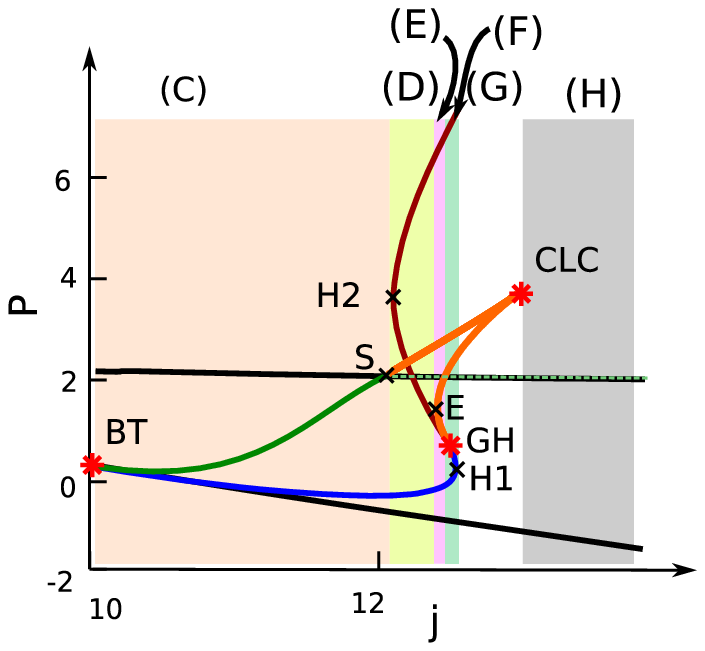}}
  \end{minipage}
 \end{center}
 \caption[Full bifurcation diagram]{Codimension two bifurcations of Jansen's model with respect to the connection strength $j$ and the input value $P$. The behavior of the system for fixed $j$ can be split into  eight zones $(A)\ldots (H)$ described in the text. The black curve corresponds to the saddle-node bifurcations manifold, the point C the cusp bifurcation point, the blue curve corresponds to the  subcritical Hopf bifurcations. It is connected to the saddle-node manifold via subcritical Boganov-Takens bifurcation at the point BT. At this point, the saddle-homoclinic bifurcations curve is plotted in green. It exists while $j$ is inferior to the value related to the Bautin bifurcation. At this point, the saddle-homoclinic bifurcations curve connects to the saddle-node manifold and generates a curve of saddle-node homoclinic bifurcations (dashed green line). The subcritical Hopf bifurcation manifold (b;ie curve) is connected to the supercritical Hopf bifurcation manifold (red curve) through a Bautin bifurcation (point GH). From this bifurcation point there is a manifold of folds of limit cycles represented in orange. The cycles undergo a cusp bifurcation at the CLC point and a saddle-node homoclinic bifurcation at the point S. }
	\label{fig:FullGlobal}
\end{figure}

We numerically observe that the system undergoes the following bifurcations of equilibria (See appendix \ref{appendix:Numerbif} and figure \ref{fig:JansenEquil}):
\begin{enumerate}
 \item A saddle-node bifurcation manifold (a curve represented in black in figures \ref{fig:JansenEquil}, \ref{fig:FICGH} and  \ref{fig:FullGlobal}),
 \item An Andronov-Hopf bifurcation manifold (a curve represented in blue for subcritical and in red for supercritical bifurcations in the same figures),\\
 \item A Cusp bifurcation C,
 \item A Bogdanov-Takens bifurcation BT,
 \item A Bautin bifurcation GH (Generalized Hopf).
\end{enumerate}

The periodic orbits are also explored using a continuation algorithm. We are especially interested in stable cycles which correspond to observable activity in the presence of noise. We chose MatCont continuation package developed by Kuznetsov, Govaerts and colleagues \cite{dhooge-govaerts-etal:03,dhooge-govaerts-etal:03b} and which is very efficient for identifying bifurcations of periodic orbits. The information of the local bifurcations described and the numerical exploration of the bifurcations of limit cycles lead us to identify the following generic bifurcations:

\begin{enumerate}
  \item A saddle homoclinic bifurcations linked with the existence of the Bogdanov-Takens bifurcation point (plain green curve of figure \ref{fig:FullGlobal}). This curve can be locally computed using the normal form of the system at this point, and continued using a continuation algorithm. Along this curve, a branch of limit cycles collapses with the saddle fixed points curve, and the period of the related cycles tends to infinity when approaching this curve. When numerically continuing this curve we observe that that at the value of parameter $j$ related to the Bautin bifurcation, this saddle-homoclinic bifurcation curve collapses with the saddle-node bifurcation manifold giving birth to:
  \item A saddle-node homoclinic bifurcation curve ( dashed green curve in Fig. \ref{fig:FullGlobal}) , and connects to a saddle-node limit cycle at a point we note $S$ in figure \ref{fig:FullGlobal}. 
  \item A Fold of Limit Cycles exists around the Bautin bifurcation point: a family of stable limit cycles and a family of unstable limit cycles collapse and disappear along a nondegenerate fold bifurcation of cycles (orange curve of \ref{fig:FICGH}).
\begin{figure}
 \begin{center}
  \includegraphics[width=.5\textwidth]{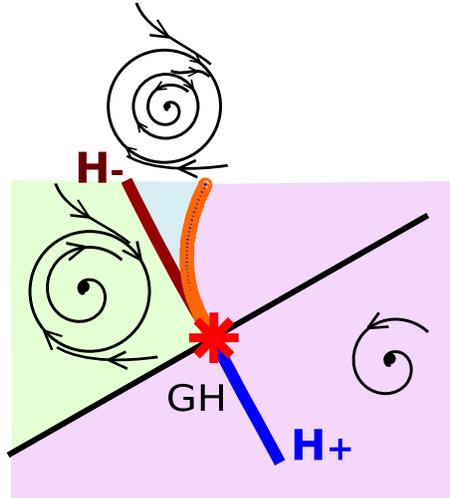}
 \end{center}
 \caption[Bautin bifurcation]{Bautin bifurcation: in region 1 the system has a single stable equilibrium an no cycle. Crossing the Hopf bifurcation boundary $H_-$ to region 2 implies the appearance of a stable limit cycle which survives when we enter region 3. Crossing the Hopf boundary $H_+$ creates an extra unstable limit cycle inside the first one, while the equilibrium regains its stability. Two cycles of opposite stability exist in region 3 and collapse at the curve T through a fold bifurcation of limit cycles that leaves a single equilibrium. The curves are computed in the case of Jansen's model.}
 \label{fig:FICGH}
\end{figure}
From the Bautin bifurcation, the manifold of FIC is continued and we numerically observe that cycles undergo a cusp bifurcation (See figure \ref{fig:FIC}).
The upper branch of limit cycles corresponds to the branch of folds of limit cycles generated at the Bautin bifurcation the cycles shrink to a single point GH. The lower branch connects to the homoclinic saddle-node manifold. At this point, the cycle corresponding to the FIC bifurcation is a saddle-node homoclinic cycle (point S). Note that the curve of folds of limit cycles originating from the Bautin point can be seen as a function of $j(P)$ which is first decreasing then increasing. The point where it changes monotony is noted $E$ in figure \ref{fig:FullGlobal}.b.
\end{enumerate}

The full bifurcation diagram is provided in figure \ref{fig:FullGlobal}. 
The analysis of the bifurcation diagram leads to classify the system into $8$ classes depending on the coupling strength. In each class the system has the same dynamical features and the same qualitative behaviors.
\begin{table}
\begin{center}
  \begin{tabular}{|l|l|l|}
  \hline
   Point & j & P \\
   \hline 
   C & 5.38 & -0.29 \\
   \hline
   BT & 10.05 & -3.07\\
   \hline
   H2 &  12.10 & 0.10\\
   \hline 
   E & 12.38 & 1.21\\
   \hline 
   GH & 12.48 & -2.58 \\
   \hline
   H1 & 12.55 & -3.10 \\
   \hline
   CLC & 12.93 & 3.75\\
   \hline
   \end{tabular}
   \end{center}
   \caption[Numerical values of the bifurcation points]{$(j,P)$ coordinates of the  different bifurcation and special points of the Jansen and Rit's system}
\end{table}

\subsubsection{Neuro-Computational features}
The different classes of parameters respresented in color in figure \ref{fig:FullGlobal} correspond to different responses to varying inputs. Eight classes can be distinguishes. We observe that the original Jansen and Rit's model is in the zone labelled (D) which is of very small extension.

\paragraph{Non-oscillating behaviors}
For $j<j_{H2}$ (zones (A), (B) and (C)), the system does not feature any stable oscillations, and therefore the cortical column will not oscillate. 
\renewcommand{\theenumi}{\Alph{enumi}}
\begin{enumerate}
\item For $j<j_C$ the system has a unique fixed point and no cycle for any input firing rate $P$. Therefore, when the input firing rate  is fixed, for any initial condition, the activity of the column will converge towards this unique equilibrium. In that case the cortical column has a quite trivial behavior: it has no oscillatory activity and converges to rest whatever the input.
\begin{figure}
 \begin{center}
  \subfigure[Equilibria]{\includegraphics[width=.45\textwidth]{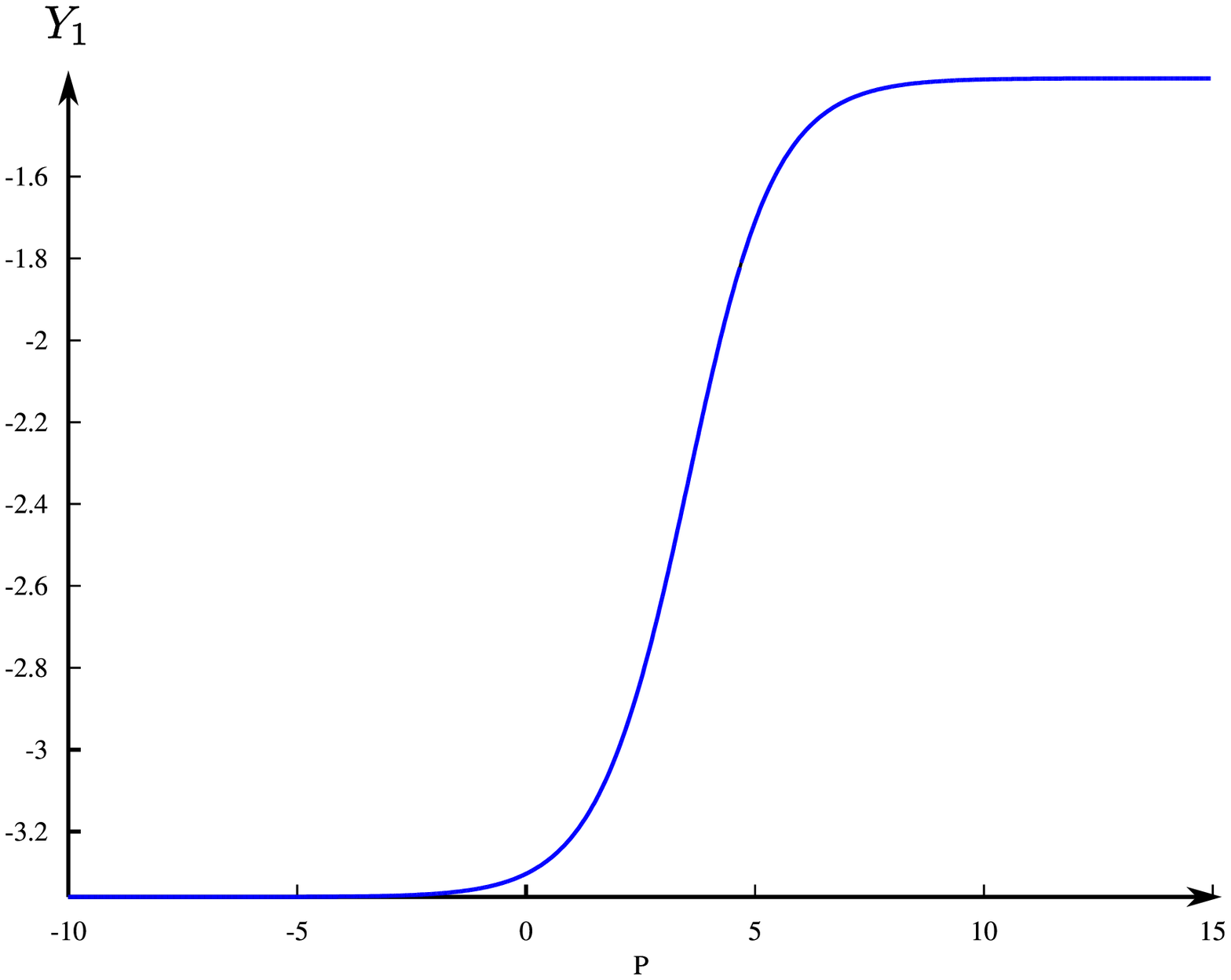}}
  \subfigure[Behaviors]{\includegraphics[width=.45\textwidth]{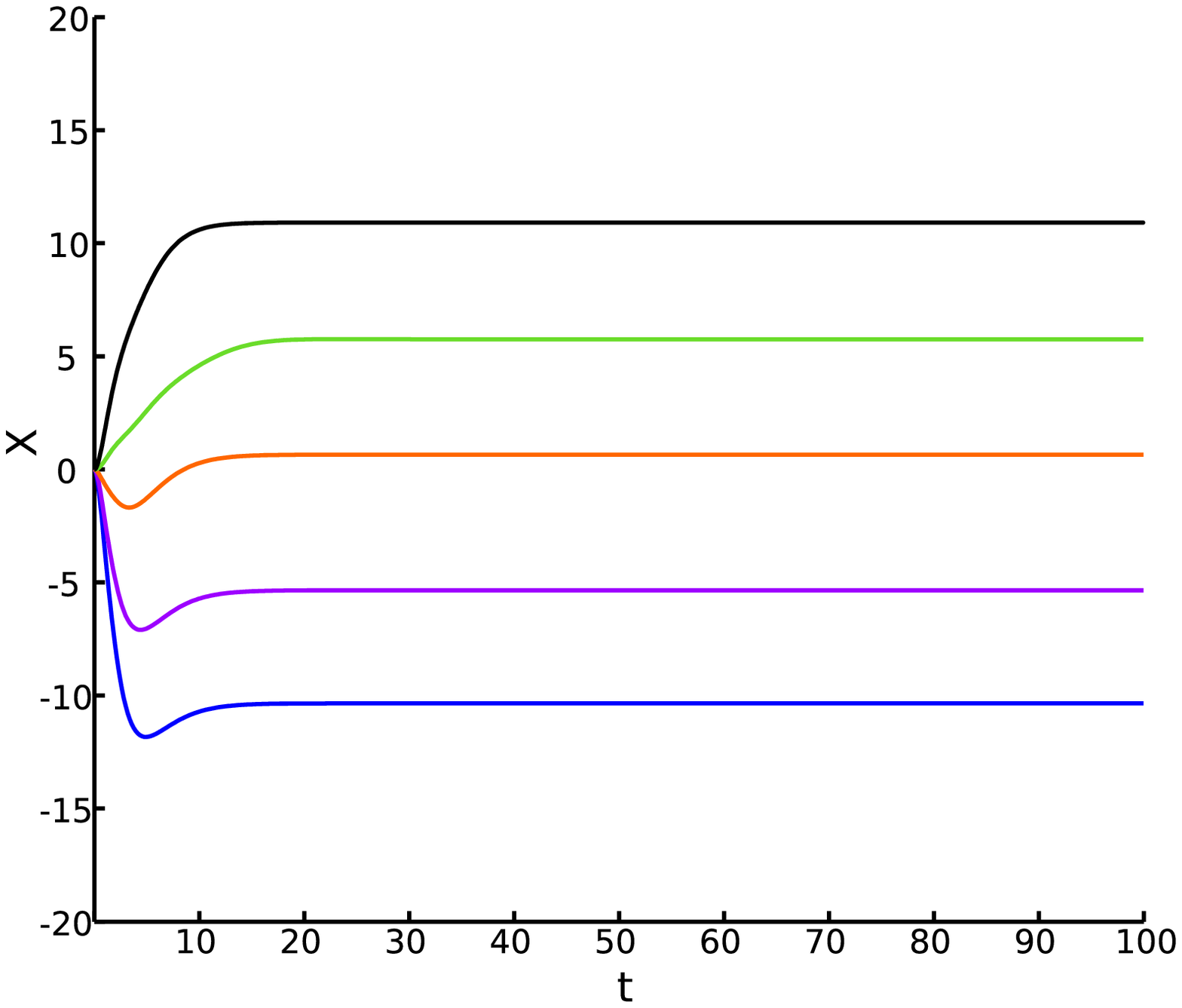}}
 \end{center}
 \caption[Case (A)]{Case (A): $j=4$. (a) equilibria: coordinate $Y_1$ as a function of the input current $P$. For each value of $P$ there exists a unique stable equilibrium. (b) behavior of the system. Variation of the  $X$ coordinate for different values of $P$: the system always converges to the unique equilibrium.}
\end{figure}
 \item For $j_C<j<j_{BT}$, the system undergoes two saddle-node bifurcations when varying the parameter $P$. Depending on the value of $P$, the system has one, two or three fixed points and no cycle (see figure \ref{fig:Bcase}). For $P \not \in [P_1,P_2]$ (i.e. not between the values of the the two saddle-node bifurcations), there is a unique fixed point which is stable and the system converges to this fixed point. When $P\in (P_1,P_2)$, there are three fixed point, one is unstable and the other two stable. The system is therefore bistable: depending on the initial condition, the activity will converge towards one or the other stable fixed point, corresponding to an up-state and down-state activity. The system also presents hysteresis when continuously varying the input in this zone of inputs. Eventually, it can switch between the two stable fixed point if perturbed. 
 \begin{figure}
 \begin{center}
  \subfigure[Equilibria]{\includegraphics[width=.45\textwidth]{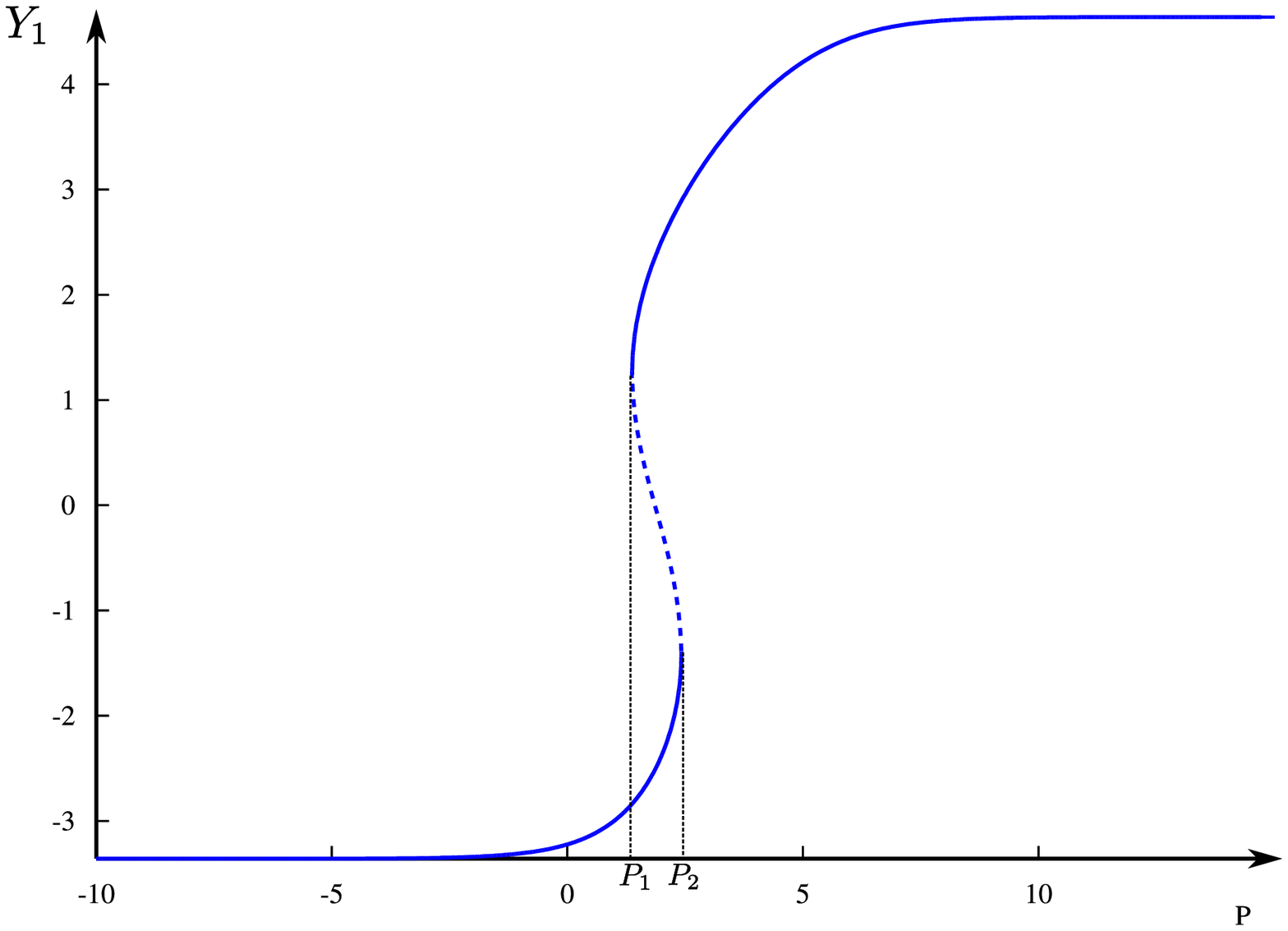}}
  \subfigure[Behaviors]{\includegraphics[width=.45\textwidth]{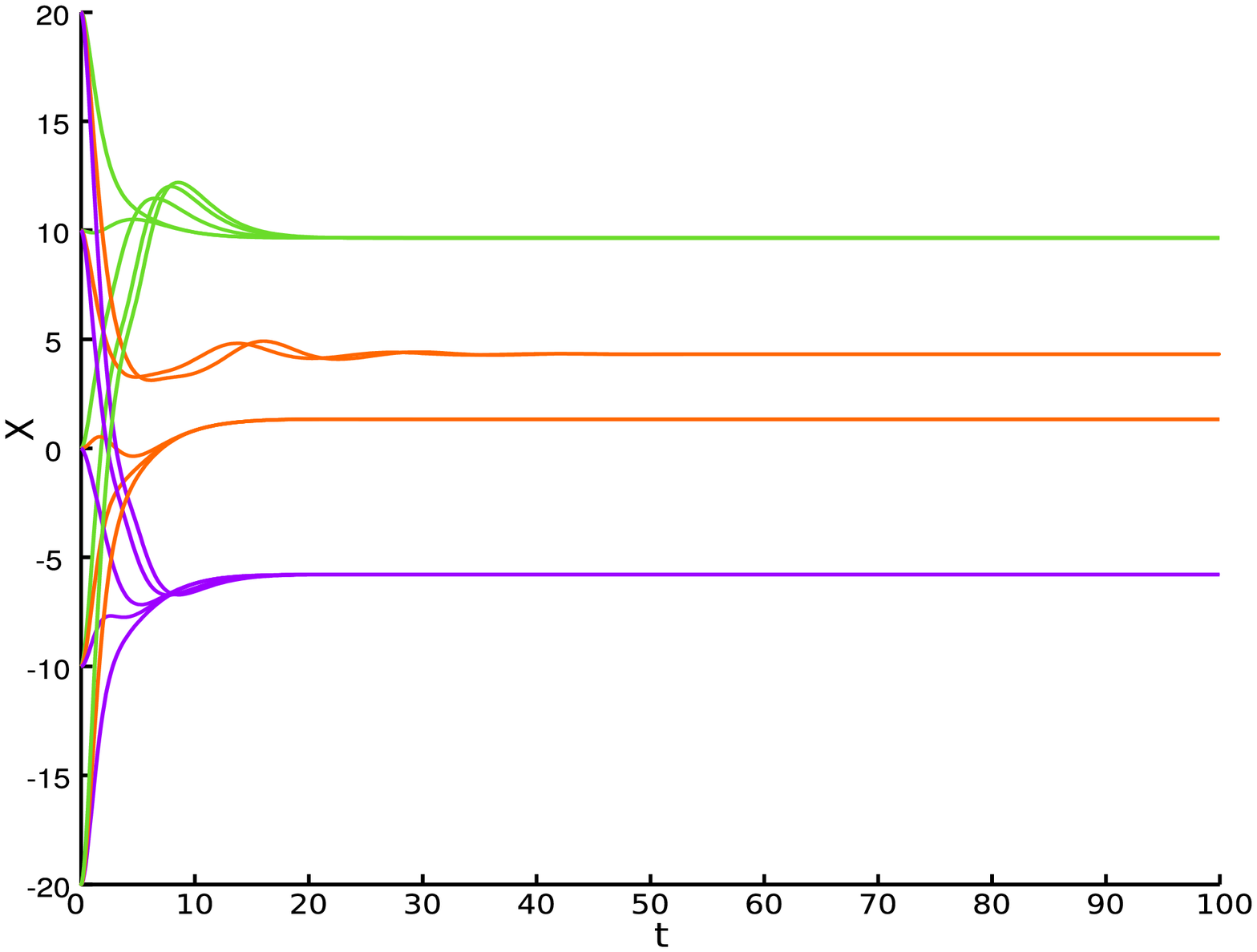}}
 \end{center}
 \caption[Case (B)]{Case (B): $j=8$. (a) equilibria (coordinate $Y_1$ as a function of the input current $P$), and stability. Continuous line: stable equilibrium, dashed line: unstable. For $P\not \in [P_1,P_2]$ there exists a unique stable equilibrium, and for $P\in (P_1,P_2)$ there exist three equilibria, two stable and one unstable. (b): Behavior of the system for different input and initial conditions. Purple: $P=-10$ and green: $P=10$ : the system always converges to the unique equilibrium. Orange: $P=2$: bistability. We observe that the activity resembles a real evoked potential.}
 \label{fig:Bcase}
\end{figure}
 \item For $j_{BT} <j < j_{H1}$: The system has two saddle-node bifurcations and a subcritical Hopf bifurcation (see figure \ref{fig:Ccase}). Therefore, the system has a unique stable fixed point when the input $P$ is not between the two saddle-node bifurcation points (i.e. $P \not \in [P_1,P_2]$), and the system converges towards this fixed point. For $P$ between the first saddle-node bifurcation value and the Hopf bifurcation (i.e. $P\in (P_1,P_H)$), the system has 2 unstable fixed points and a stable fixed point, and generically converges towards the stable fixed point. In the zone between the Hopf bifurcation and saddle homoclinic point ($P\in (P_H,P_{Sh})$) , the system presents two stable fixed points, an unstable fixed point and an unstable limit cycle. Depending on the initial condition, the system will either converge to one or the other fixed point. When $P$ is greater than the saddle-homoclinic bifurcation value and below the greatest saddle-node value ($P\in (P_{Sh},P_2)$), the system has two stable fixed points and an unstable fixed point, and its behavior is similar to the behavior in the previous case. We can see that the system returns to equilibrium via oscillations. Hence in the presence of noise, the system will present oscillations at a certain frequency superimposed to its noisy behavior.
  \begin{figure}
 \begin{center}
  \subfigure[Equilibria]{\includegraphics[width=.45\textwidth]{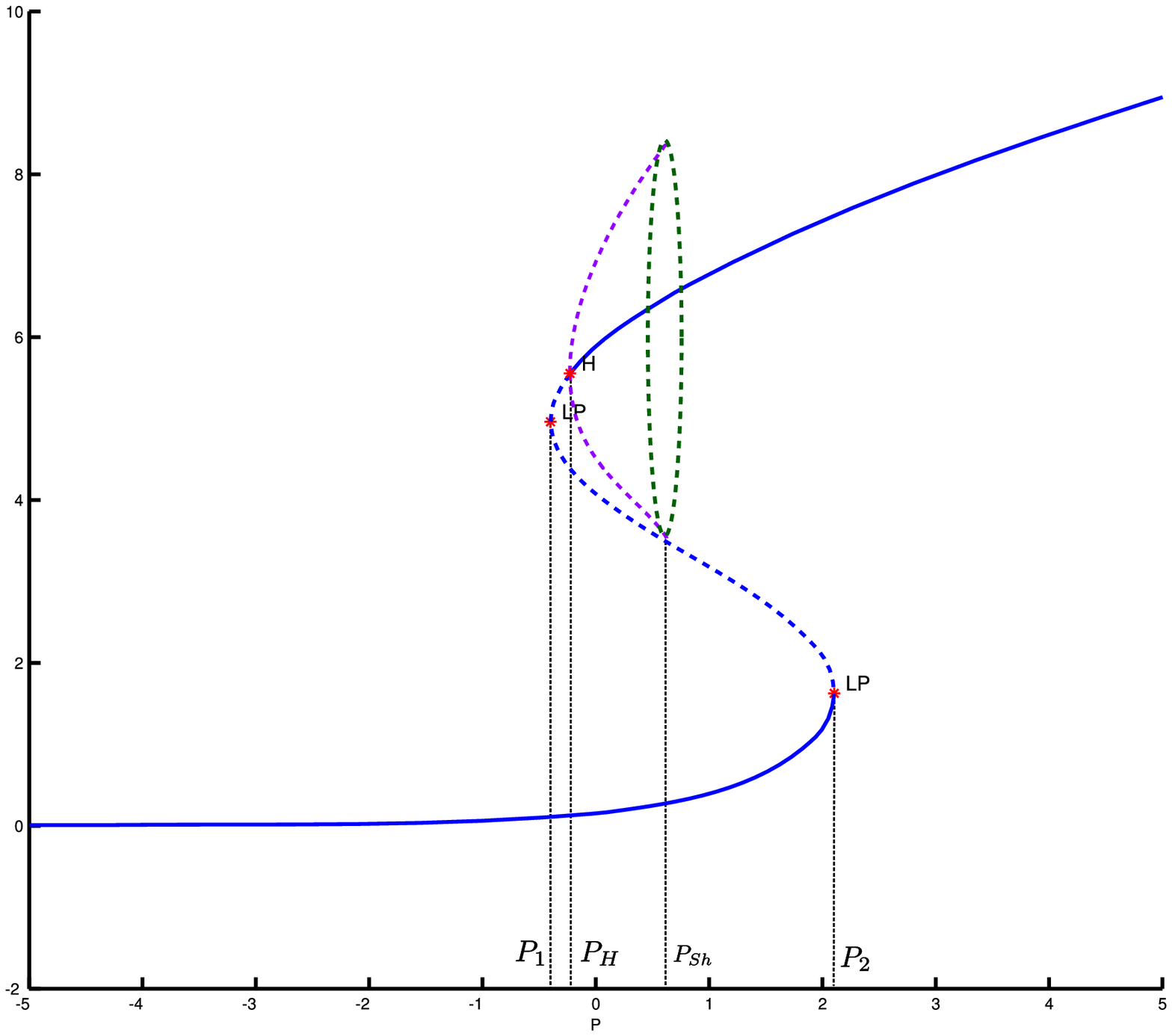}}
  \subfigure[Behaviors]{\includegraphics[width=.45\textwidth]{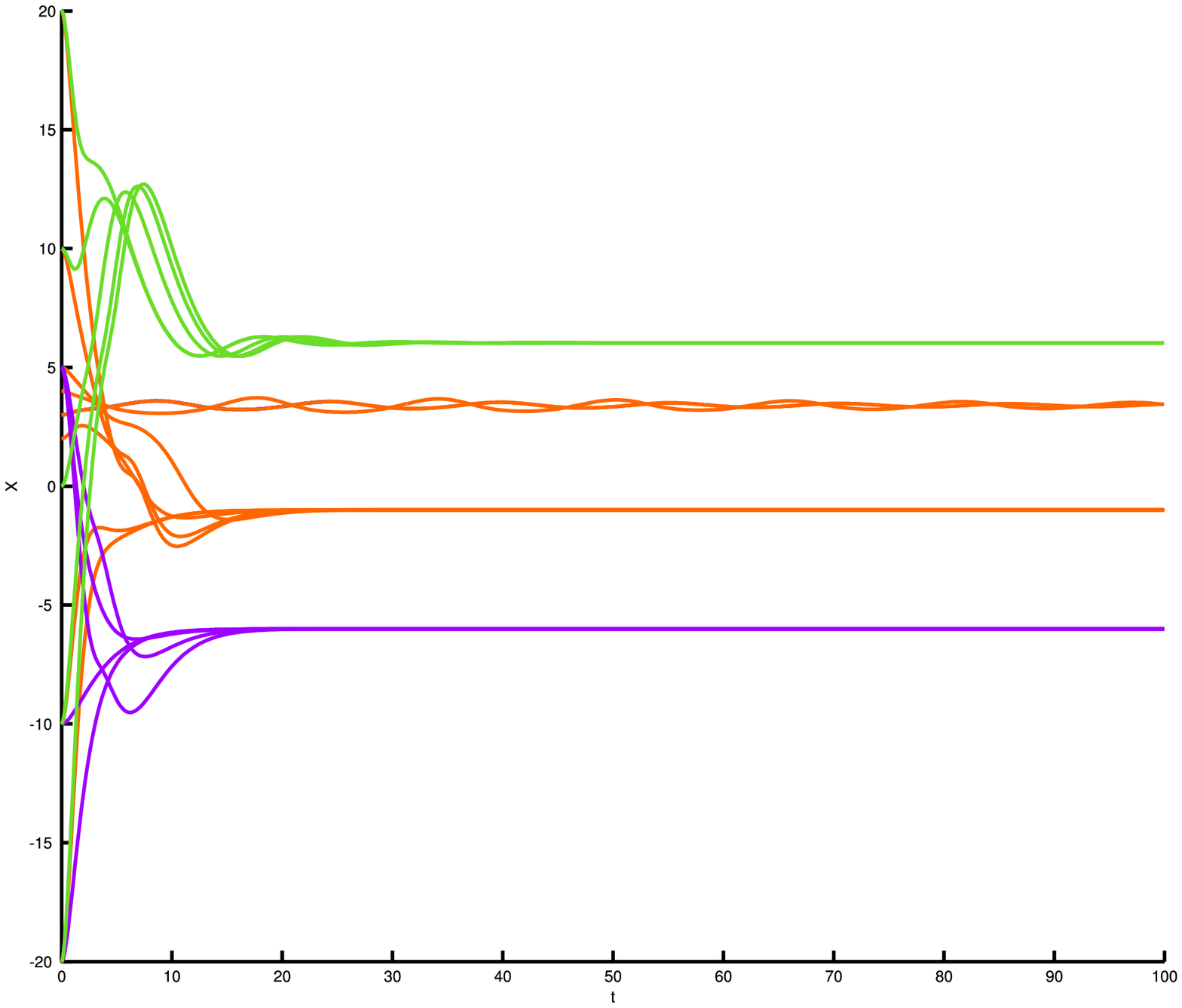}}
 \end{center}
 \caption[Case (C)]{Case (C): $j=11$. (a): Equilibria and bifurcations (coordinate $Y_1$ as a function of the input current $P$), and stability. Blue lines: Continuous : stable equilibrium, dashed : unstable. Dashed purple line: unstable cycles. Green cycle: saddle homoclinic orbit. For $P\not \in [P_1,P_2]$ there exists a unique stable equilibrium, and for $P\in (P_1,P_2)$ there exist three equilibria: for $P \in (P_1,P_H)$ two unstable and one stable, and no limit cycle. For $P \in (P_H,P_2)$, two stable and one unstable. For $P\in [P_H,P_{Sh}]$, there exists an unstable limit cycle. (b): Behavior of the system for different input and initial conditions. Purple: $P=-10$ and green: $P=10$ : the system always converges to the unique equilibrium. Orange: $P=0$: bistability and damped oscillations. We observe that the activity resembles a real evoked potential.}
 \label{fig:Ccase}
\end{figure}
\end{enumerate}

\paragraph{Rhythmic activity}
For values of the connectivity greater than $j_{H_2}$, the system will always present a supercritical Hopf bifurcation and therefore a stable periodic orbit, corresponding to a rhythmic activity of the column.

\begin{enumerate}\setcounter{enumi}{3}
\item For $j \in [j_{H_2}, j_E]$, the system undergoes two saddle-node bifurcations, two supercritical and a subcritical Hopf bifurcations, and one fold of limit cycles. It is the case of the original Jansen and Rit's model  \ref{tab:JansenNumerics}. We label the bifurcation points as in figure \ref{fig:DCase}. 

For $P<P_1$, the system has a unique stable fixed point and for any initial condition, it converges towards this point. 

For $P_1<P<P_{H_1}$, the system has three fixed points, one of which is stable and the other two unstable, and the system still converges towards the unique fixed point. 

For $P_{H_1}<  P<P_{H_2}$, the largest fixed point becomes stable and there is an unstable limit cycle (represented in dashed purple in figure \ref{fig:DCase}). In this region, the system converges towards one of the two stable fixed point depending on the initial condition of the system. 

For $P_{H_2}<P<P_2$, the largest fixed point loses stability via a subcritical Hopf bifurcation and a stable periodic exists with frequency of about 10Hz corresponding to purely alpha activity. In this region, the system either oscillates around the periodic orbit and presents alpha-activity, or converges to the stable fixed point. 

For $P_2<P<P_{FLC}$, the system has no fixed point and two stable limit cycles: the cycle corresponding to the continuation of the subcritical Hopf bifurcation which corresponds to alpha activity, and a large amplitude cycle with a frequency ranging from 0 to 5 Hz which corresponds to an epileptiform activity. The system selects one of these cycles depending on the initial condition, and can switch from one activity to another  when the system is perturbed. Assume that we slowly increase the input $P$. If the system was in the down equilibrium state, it will converge to the epileptic cycle when crossing the bifurcation, and if it was in the up state to the alpha cycle. 

For $P_{FLC}<P<P_{H_3}$, the system has a unique stable trajectory which is a periodic orbit with frequency close to 10Hz, and for any initial condition the system will converge towards this cycle. Eventually, for $P>P_{H_3}$ the system has a unique fixed point and for any initial condition, the state  converges towards this fixed point.
\begin{figure}
 \begin{center}
    \subfigure[Equilibria]{\includegraphics[width=.45\textwidth]{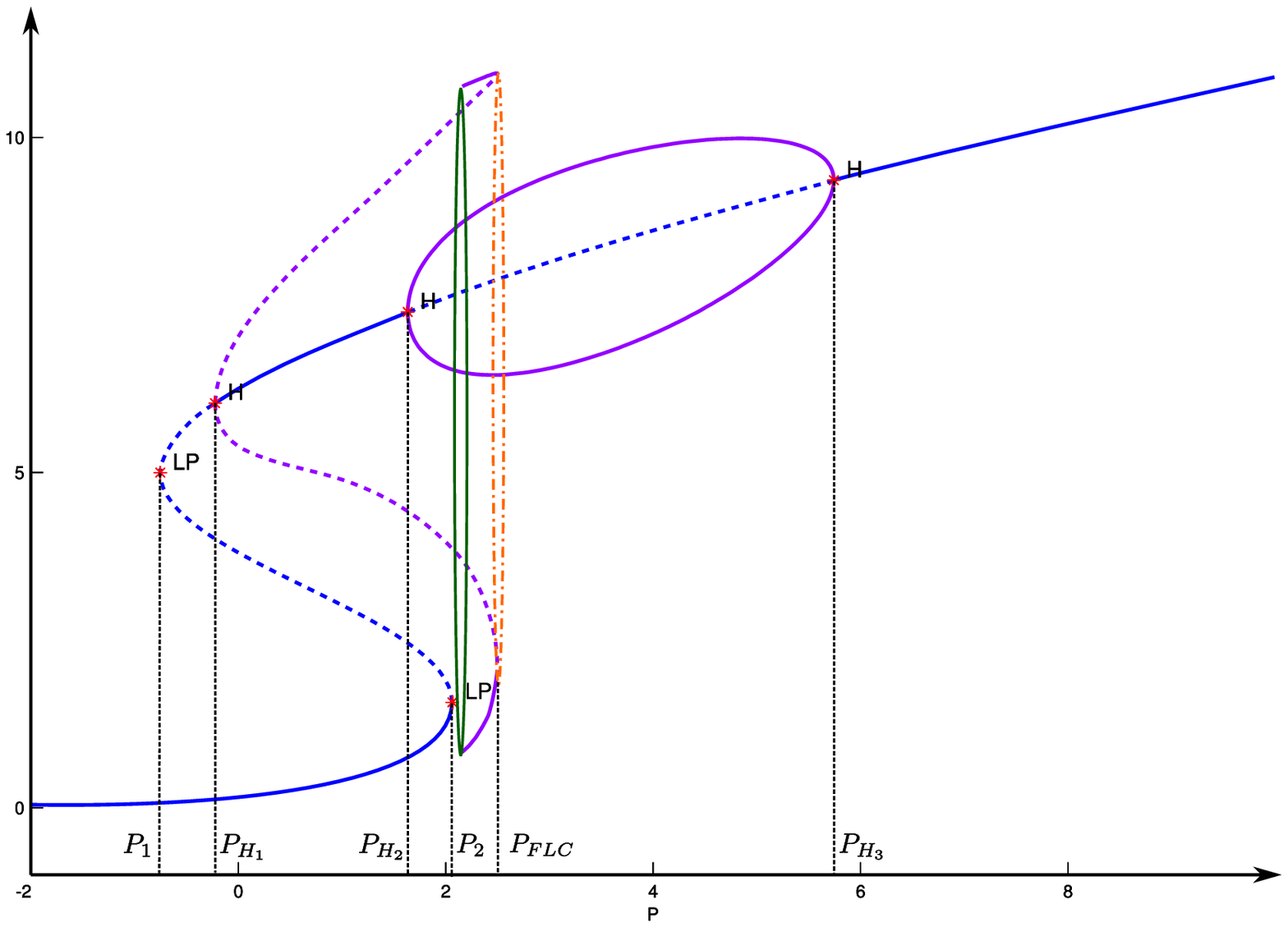}}\quad
    \subfigure[Behaviors]{\includegraphics[width=.45\textwidth]{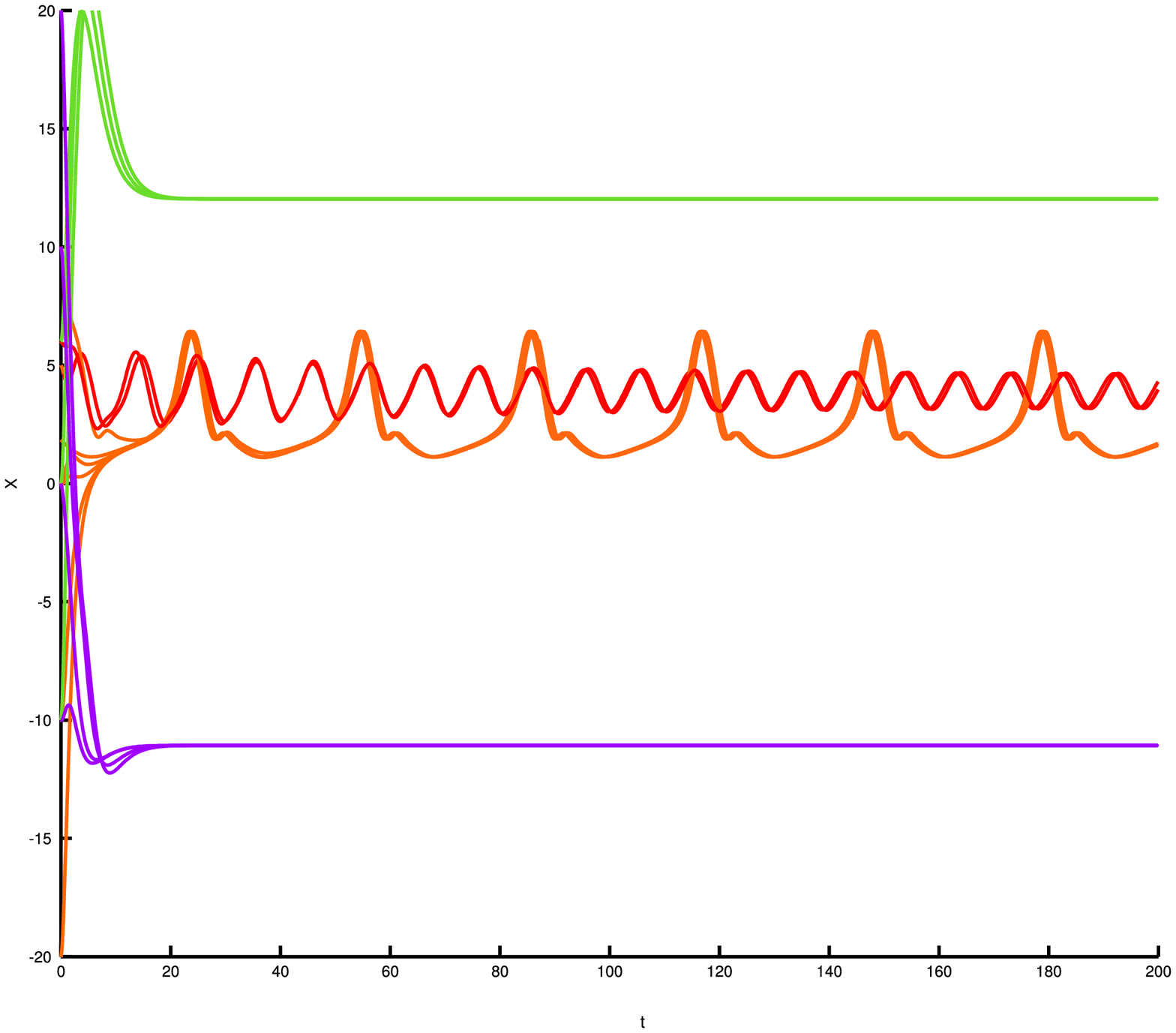}}\\
    \subfigure[Periods]{\includegraphics[width=.45\textwidth]{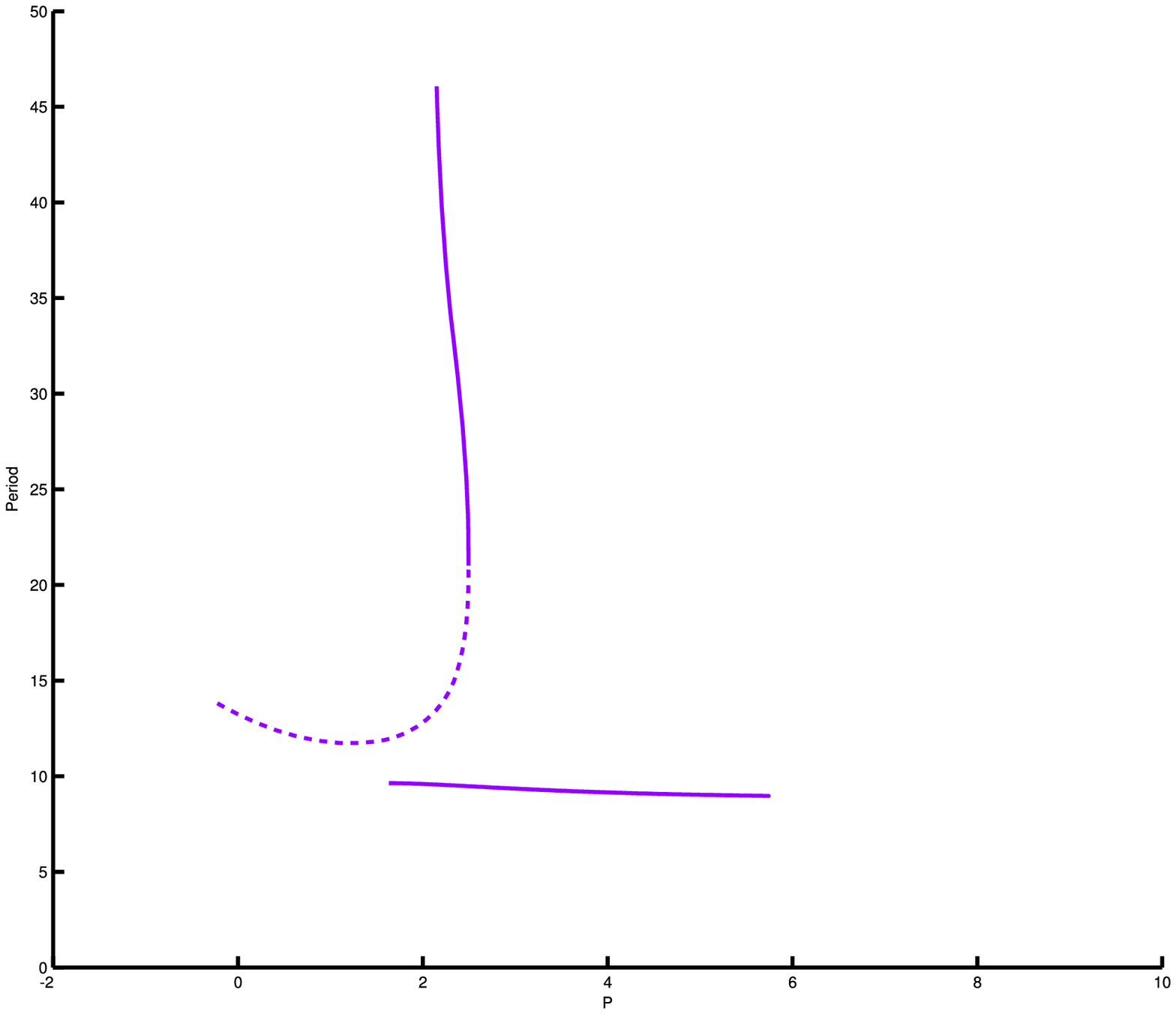}}\\
 \end{center}
 \caption[Case (D)]{Case (D): $j=12.285$. (a): Equilibria and bifurcations (coordinate $Y_1$ as a function of the input current $P$), and stability. Blue lines: Continuous: stable equilibrium, dashed: unstable. Purple lines: dashed: unstable cycles, continuous: stable. Green cycle: saddle homoclinic orbit, orange cycle: fold of limit cycles (see text) (b): Behavior of the system for different input and initial conditions. Purple: $P=-10$ and green: $P=10$: the system always converges to the unique equilibrium. Red and orange: $P=2.3$: bistability of cycles: orange: epileptic spikes, red: alpha activity. (c): period of the limit cycles.}
 \label{fig:DCase}
\end{figure}

 \item For $J \in [j_E, j_{GH}]$, the number and stability of the fixed points are the same as in the previous case. But in this case, the structure of the cycles is more complex. The family of unstable periodic orbits originating from the subcritical Hopf bifurcation is connected to the family of limit cycles of the supercritical Hopf bifurcation $H_2$ associated with the smaller $P$ value, and the branch of limit cycles of the supercritical Hopf bifurcation $H_3$ associated with the greatest $P$ value undergoes two fold bifurcations of cycles and disappears via saddle-node homoclinic bifurcation (see figure \ref{fig:Ecase} and labels herein). For $P<P_{FLC_1}$ the behavior of the system is exactly the same as in case (D): the system generically converges to the unique stable fixed point for $P<P_{H_2}$, and either converges to the stable fixed point or to the stable cycle depending on the initial condition. For $P_{FLC_1}<P<P_{FLC_2}$, the system return to the down-state equilibrium. For $P_{FLC_2}<P< P_2$, the system has a stable fixed point (down-state equilibrium), a stable limit cycle corresponding to alpha-like activity and an unstable limit cycle. In this zone of input the system will experience bistability between alpha-activity and rest. For $P_2<P<P_{FLC_3}$, the system has no stable fixed points and three cycles, one of which is unstable, another one corresponding to alpha activity and the third one to epileptiform activity. For $P>P_{FLC_3}$ the behavior is the same as in the case (D) for $P>P_{FLC}$: while $P<P_{H_3}$ the system presents purely alpha oscillations, and for $P>P_{H_3}$ the system converges to the unique up-state equilibrium.
 \begin{figure}
 \begin{center}
    \subfigure[Equilibria]{\includegraphics[width=.85\textwidth]{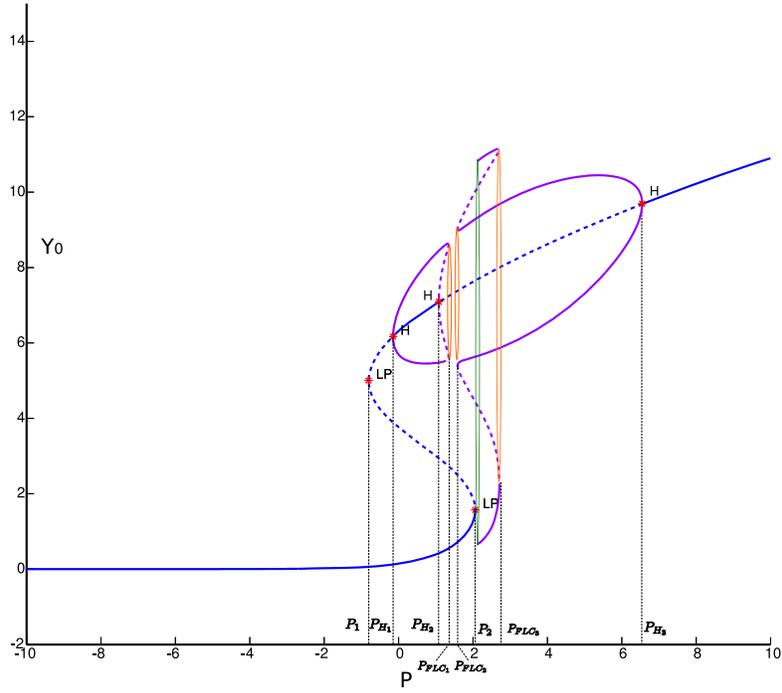}}\\
    \subfigure[Periods]{\includegraphics[width=.45\textwidth]{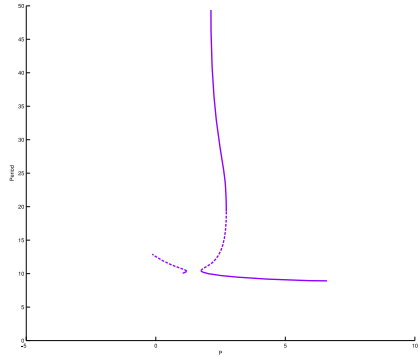}}
 \end{center}
 \caption[Case (E)]{Case (E): $j=12.42$. (a): Equilibria and bifurcations (coordinate $Y_1$ as a function of the input current $P$), and stability. Blue lines: Continuous: stable equilibrium, dashed: unstable. Purple line: dashed: unstable cycles, plain: stable. Green cycle: saddle homoclinic orbit, orange cycles: folds of limit cycles (see text). (b): Behavior of the system for different inputs and initial conditions. Purple: $P=-4$ and green: $P=10$: the system always converges to the unique equilibrium. Red and orange: $P=2.4$: alpha and epileptic activity, light blue: $P=1.1$: rest and alpha oscillations, pink: $P=4$: purely alpha activity. E and time removed because arxiv refuses large files.}
 \label{fig:Ecase}
\end{figure}
 
 \item For $j_{GH}<j<j_{H_1}$, the system has two subcritical Hopf bifurcations whose family of limit cycles are connected, and a supercritical Hopf bifurcation whose limit cycles undergo two saddle-node bifurcations of limit cycles and collapse on the saddle-node manifold via saddle-node homoclinic bifurcation (see figure \ref{fig:Fcase}). For $P<P_{H_1}$ the system converges to the down-state equilibrium. For $P_{H_1}<P<P_{H_2}$, the system is bistable and either converges to the upstate equilibrium or to the downstate equilibrium depending on the initial condition. For $P_{H_2}<P<P_2$ the system converges to the downstate equilibrium. Therefore for $P<P_{H_2}$, the system only presents damped subthreshold oscillation and no real rhythmic activity. For $P_2<P<P_{FLC_1}$ the system is in a pure epileptic activity, for $P_{FLC_1}<P<P_{FLC_2}$ the system presents both epileptic activity and alpha activity, and for $P_{FLC_2}<P<P_{H_3}$ the system only presents alpha activity. In this case, when slowly varying the input $P$, the system will always present epileptic activity: indeed, it is the only stable activity in a certain range of parameters. 
 \begin{figure}
 \begin{center}
    \subfigure[Equilibria]{\includegraphics[width=.85\textwidth]{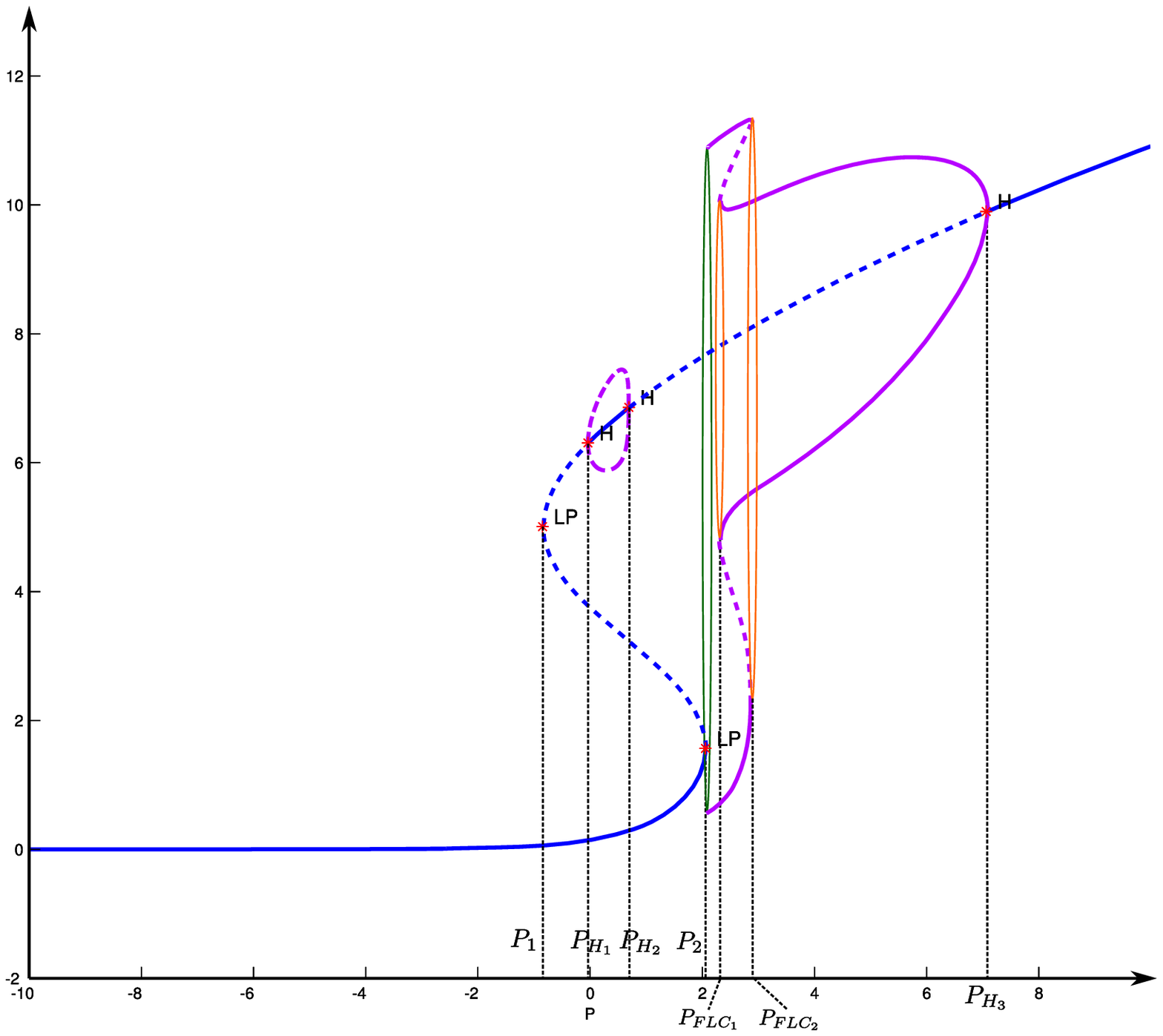}}\\
    \subfigure[Behaviors]{\includegraphics[width=.6\textwidth]{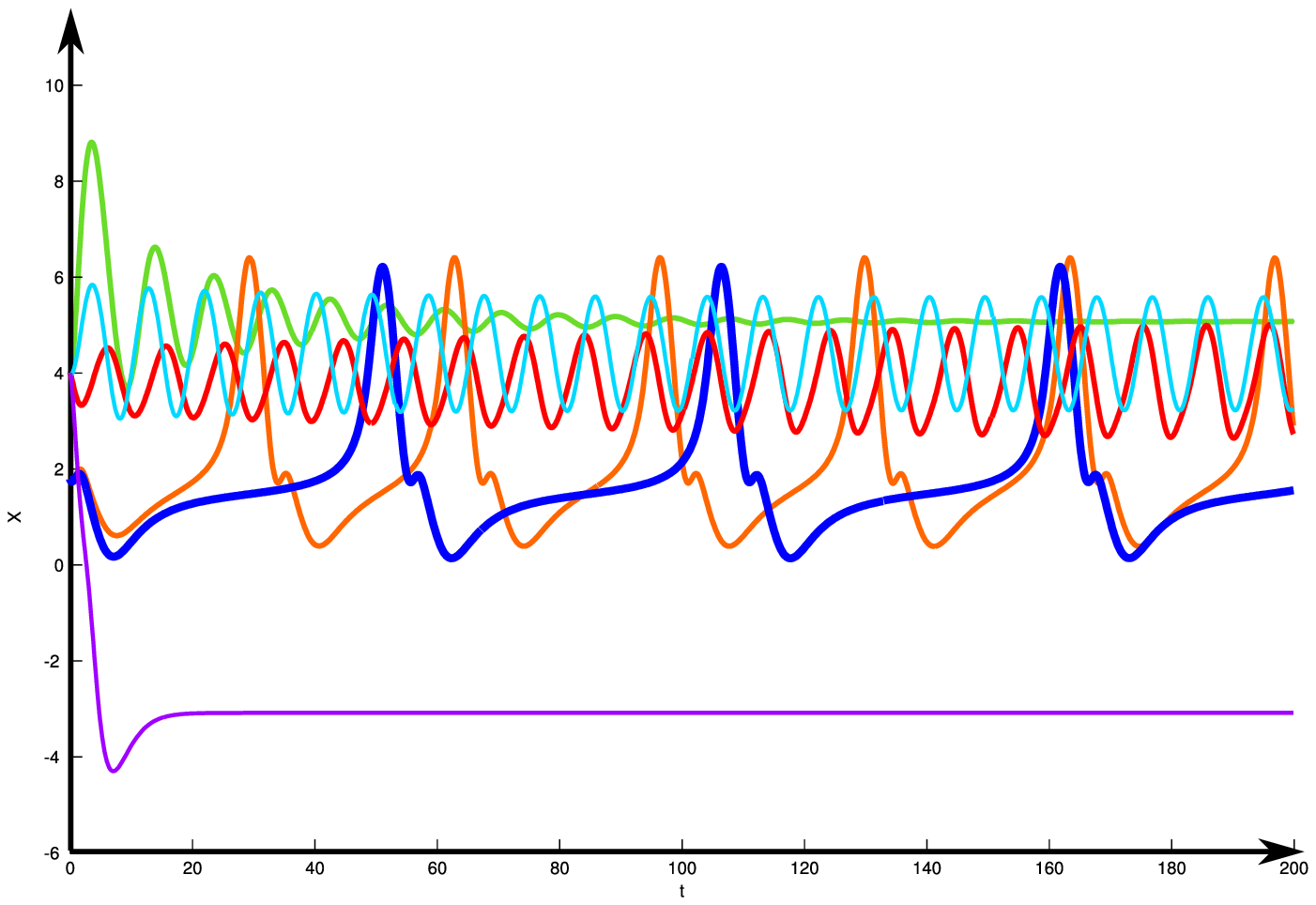}}\quad
    \subfigure[Periods]{\includegraphics[width=.3\textwidth]{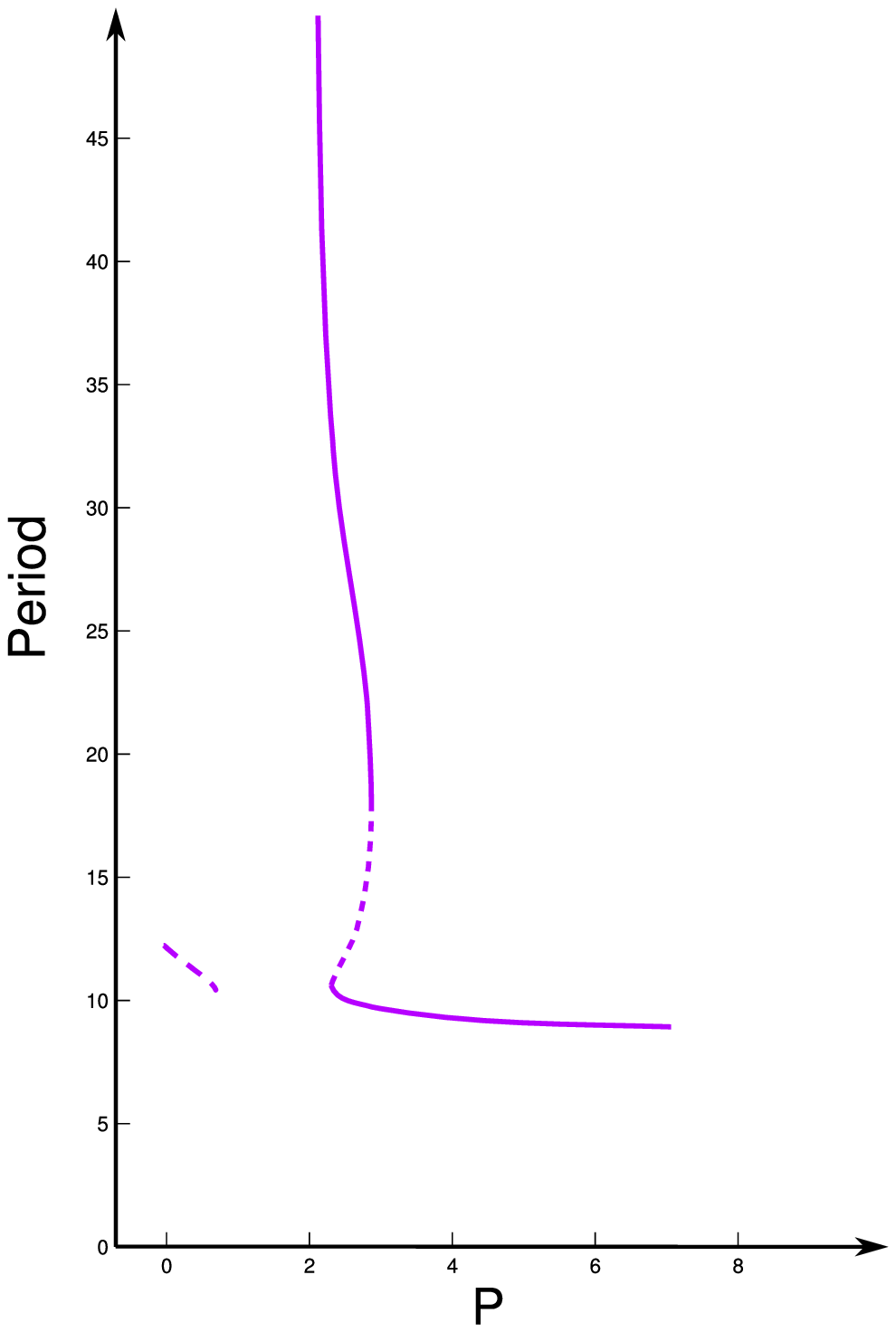}}
 \end{center}
 \caption[Case (F)]{Case (F): $j=12.5$. (a): Equilibria and bifurcations (coordinate $Y_1$ as a function of the input current $P$), and stability. Blue lines: Continuous: stable equilibrium, dashed: unstable. Purple line: dashed: unstable cycles, continuous: stable. Green cycle: saddle homoclinic orbit, orange cycles: folds of limit cycles (see text). (b): Behavior of the system for different inputs and initial conditions. Purple: $P=-2$ and green: $P=9$ : the system always converges to the unique equilibrium. Dark blue: $P=2.05$: purely epileptic activity: slow waves of high amplitude. Red and orange: $P=2.3$: alpha and epileptic activity coexist, light blue: $P=5$: only alpha oscillations.}
 \label{fig:Fcase}
\end{figure}
 \item For $j_{H_1}<j<j_{CLC}$, the system has two saddle-node bifurcations and a supercritical Hopf bifurcation, whose limit cycles undergo two saddle-node bifurcations and disappear via a saddle-node homoclinic bifurcation. For all $P<P_2$, the system always converges to the downstate equilibrium. For $P_2<P<P_{FLC_1}$ the system presents epileptic spikes, and for $P_{FLC_1}<P<P_{FLC_2}$ bistability with epileptic spikes and alpha activity. For $P_{FLC_2}<P<P_{H_1}$ the system presents only a stable alpha activity and for $P>P_{H_1}$ the system always returns to an upstate equilibrium whatever the initial condition (see figure \ref{fig:Gcase})
 \begin{figure}
 \begin{center}
    \subfigure[Equilibria]{\includegraphics[width=.6\textwidth]{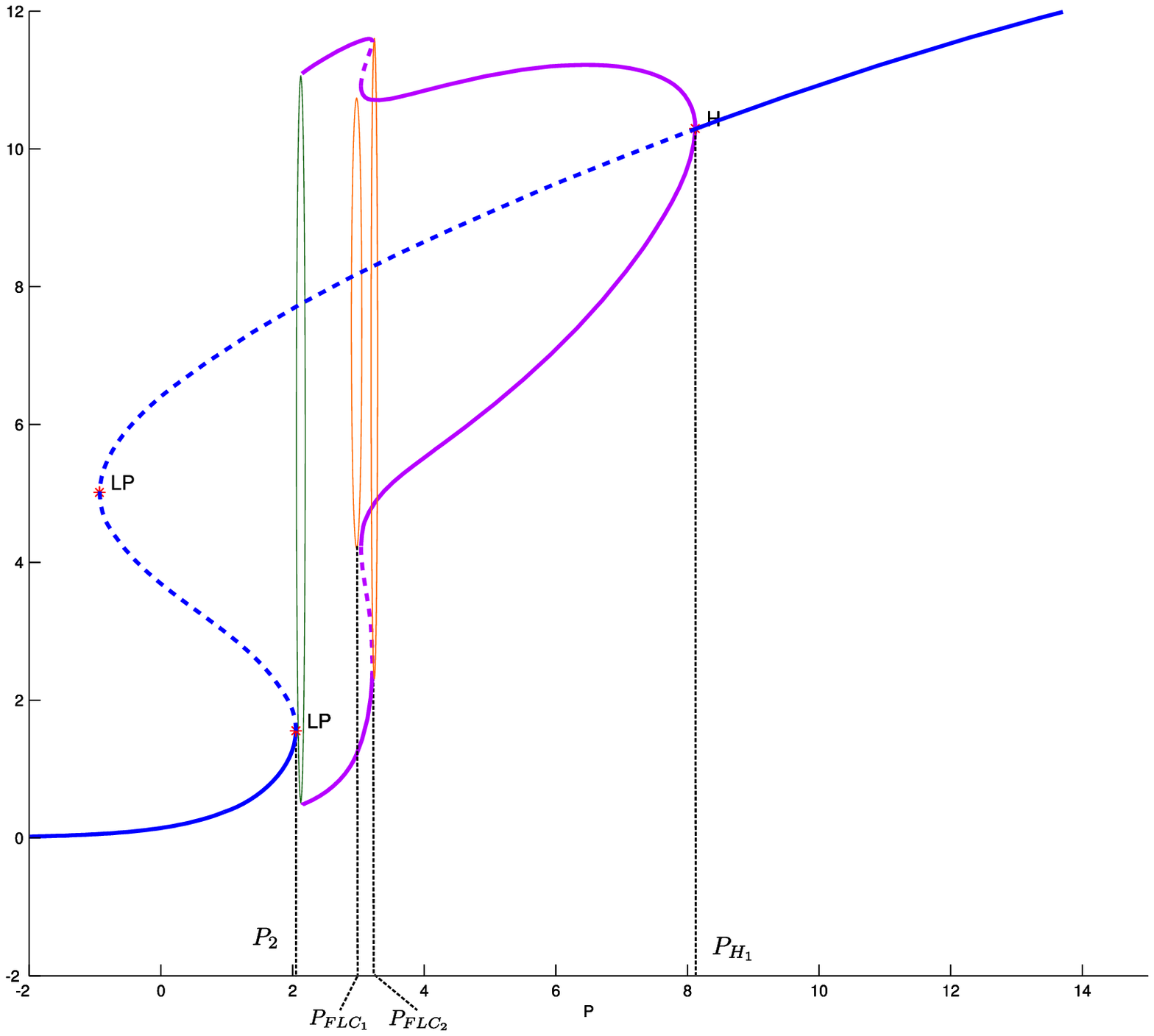}}\\
    \subfigure[Behaviors]{\includegraphics[width=.45\textwidth]{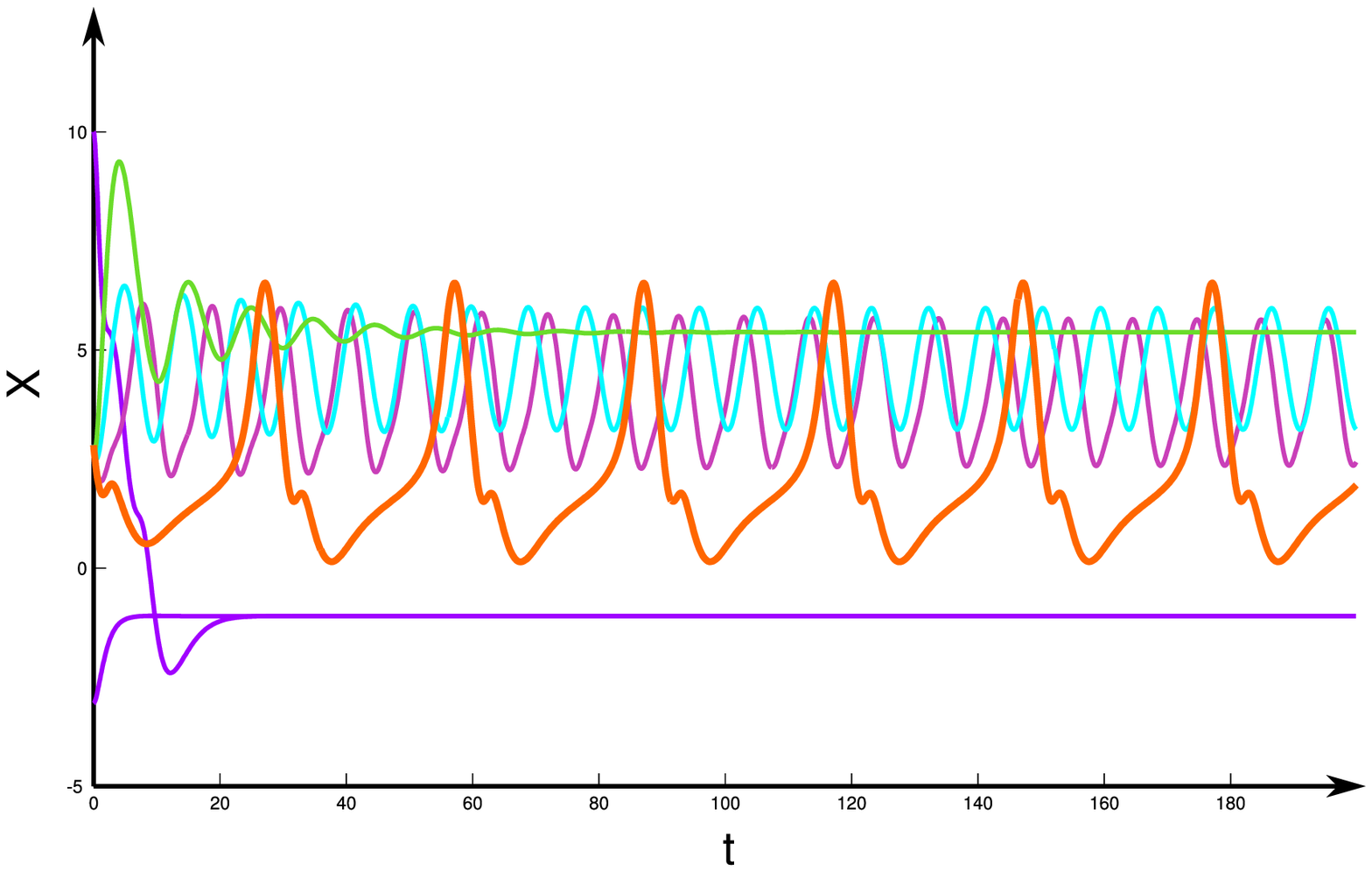}}\quad
    \subfigure[Periods]{\includegraphics[width=.45\textwidth]{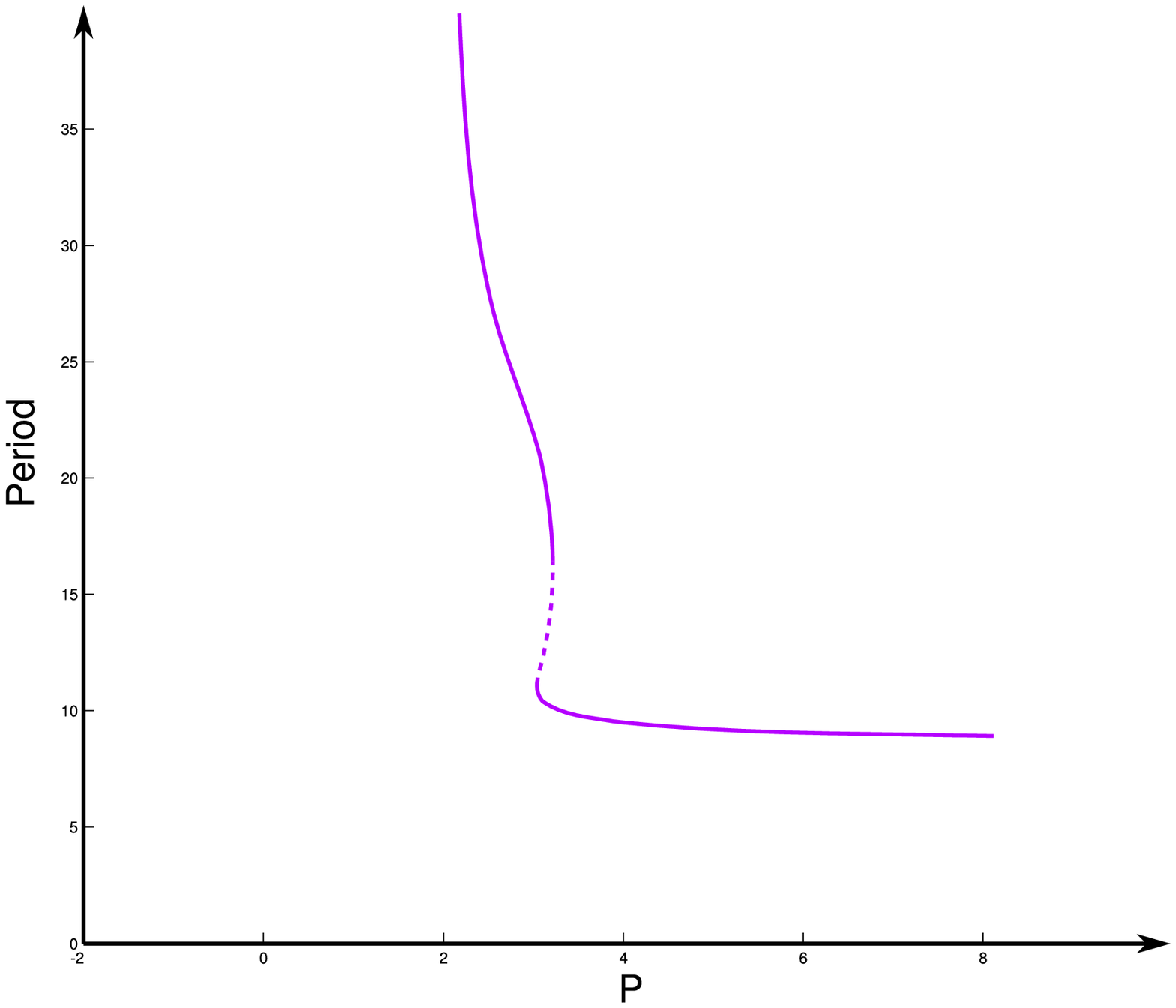}}
 \end{center}
 \caption[Case (G)]{Case (G): $j=12.7$. (a): Equilibria and bifurcations (coordinate $Y_1$ as a function of the input current $P$), and stability. Blue lines: Continuous: stable equilibrium, dashed: unstable. Purple line: dashed: unstable cycles, continuous: stable. Green cycle: saddle homoclinic orbit, orange cycles: folds of limit cycles (see text). (b): Behavior of the system for different inputs and initial conditions. Purple: $P=0$ and green: $P=11$: the system always converges to the unique equilibrium. Orange: $P=2.3$: purely epileptic activity: slow waves of high amplitude. Pink: $P=3.1$: alpha and epileptic activities coexist, light blue: $P=5$: only alpha oscillations.}
 \label{fig:Gcase}
\end{figure}

 \item For $j>j_{CLC}$ the system features two saddle-node bifurcations and a supercritical Hopf bifurcation whose family of limit cycle is regular and disappears via saddle-node homoclinic bifurcation (see figure \ref{fig:Hcase}). For $P<P_2$ or $P>P_{H_1}$ the system converges to the unique fixed point, and for $P_2<P<P_{H_1}$ the system presents oscillations. This case is very interesting from a neuro-computational point of view. Indeed, the family of limit cycles created presents a stable region of oscillations around 10Hz corresponding to alpha activity as in the previous case, for a quite large set of inputs $P_{\alpha}<P<P_{H_1}$. At $P=P_{\alpha}$ it abruptly switches from alpha activity to theta activity where it stays for $P_{\theta}<P<P_{\alpha}$ and eventually switches regularly to delta activity for $P_2<P<P_{\theta}$. In this region of parameter therefore the system presents the rhythm of the normal sleep activity.
 \begin{figure}
 \begin{center}
    \subfigure[Equilibria]{\includegraphics[width=.85\textwidth]{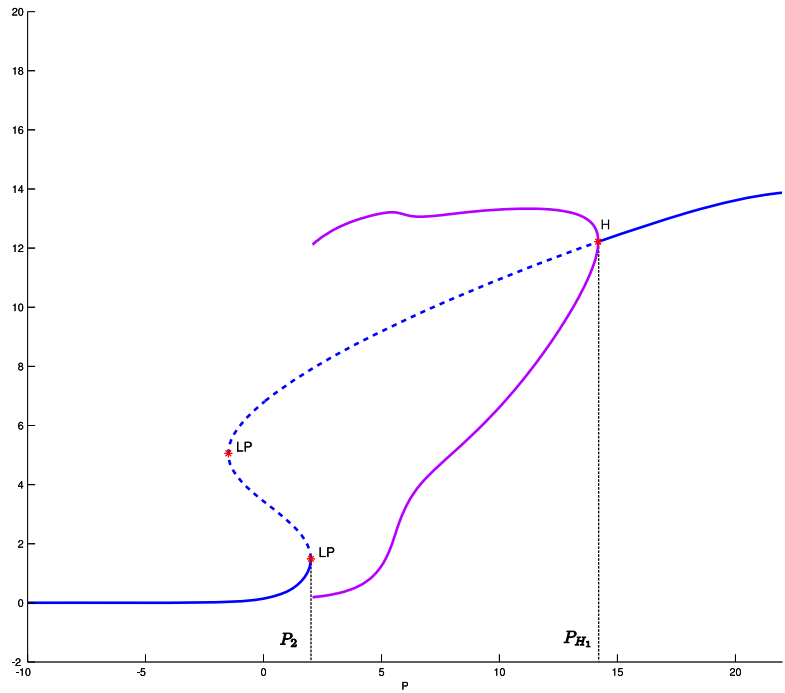}}\\
    \subfigure[Behaviors]{\includegraphics[width=.45\textwidth]{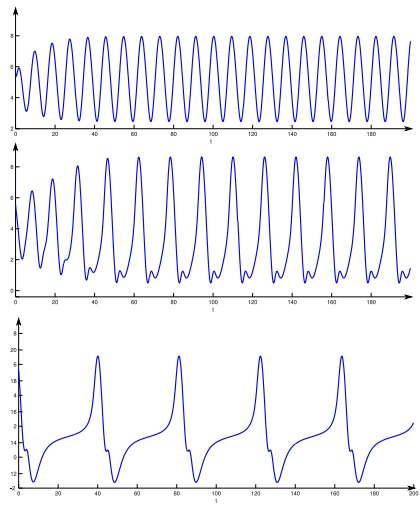}}\quad
    \subfigure[Periods]{\includegraphics[width=.45\textwidth]{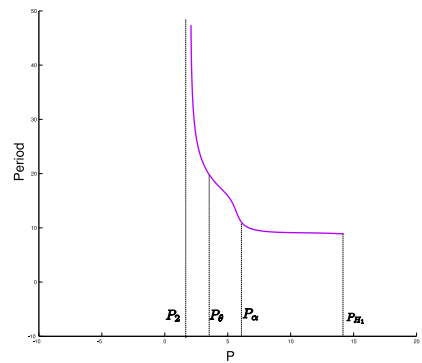}}
 \end{center}
 \caption[Case (H)]{Case (H): $j=14$. (a) Equilibria and bifurcations (coordinate $Y_1$ as a function of the input current $P$), and stability. Blue lines: Continuous: stable equilibrium, dashed: unstable. Purple line: dashed: unstable cycles, plain: stable (see text). (b) Behavior of the system for different inputs and initial conditions. Up: $P=10$: alpha activity, Middle: $P=5$ theta activity. Down $P=2.3$: delta activity.}
 \label{fig:Hcase}
\end{figure}
 
\end{enumerate}
\newpage

\section{Influence of other parameters in Jansen and Rit's model}
\subsection{Effect of the PSP amplitude ratio $G$}

\begin{figure}[!h]
 \begin{center}
  \includegraphics[width=.8\textwidth]{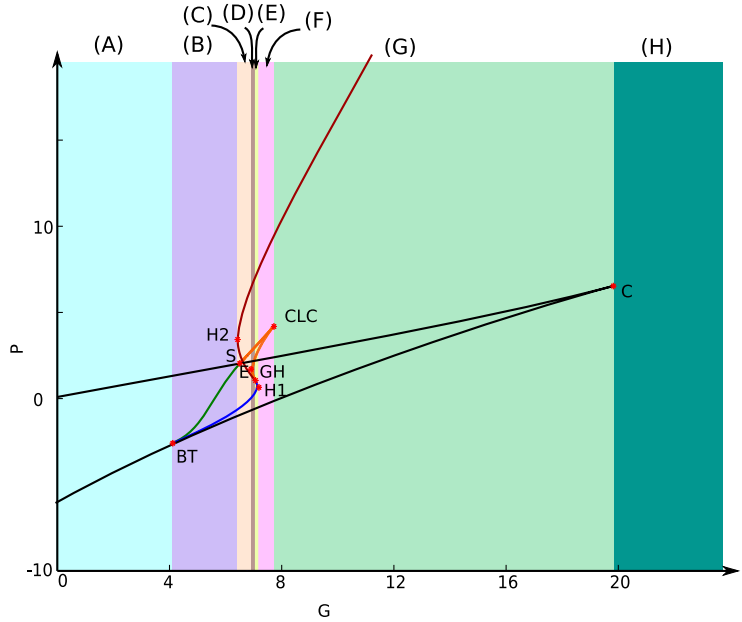}
 \end{center}
 \caption[Jansen and Rit's bifurcations with respect to $G$ and $P$]{Codimention 2 bifurcations in Jansen and Rit's model with respect to the PSP amplitude ratio $G$ and the input $P$. The behavior of the system for fixed values of $G$ can be split into eight zones $(A)\ldots (H)$ described in the text. The black curve corresponds to the saddle-node bifurcations manifold, the point C to the cusp bifurcation point, the blue curve to the  subcritical Hopf bifurcations manifold. It is connected to the saddle-node manifold via a subcritical Boganov-Takens bifurcation at the point BT. At this point, the saddle-homoclinic bifurcations curve is plotted in green and connects to the saddle-node manifold via saddle-homoclinic bifurcation, becoming a saddle-node homoclinic bifurcation curve (dashed green line). The subcritical Hopf bifurcation manifold is connected to the supercritical Hopf bifurcation manifold (red curve) through a Bautin bifurcation (point GH). From this bifurcation point there is a manifold of folds of limit cycles represented in orange. Cycles undergo a cusp bifurcation at the CLC point and a saddle-node homoclinic bifurcation at the point S. (See values of the parameters in table \ref{table:JansenBifsGP})}
 \label{fig:ZonesGP}
\end{figure}

The bifurcation structure as a function of the PSP amplitude ratio $G$ is very similar to the one corresponding to the coupling strength $j$. The system also presents a cusp, a Bogdanov-Takens and a Bautin bifurcation, two branches of saddle-node of limit cycles collapsing at a cusp of limit cycles. The same types of behaviors are observed. We decompose here again the bifurcation diagram into zones depending of $G$, and use the same notations as in the previous case. The bifurcation diagram and the decomposition in zones is given in figure \ref{fig:ZonesGP}.

\begin{table}[!h]
\begin{center}
  \begin{tabular}{|l|l|l|}
  \hline
   Point & G & P \\
   \hline 
   BT & 3.06 & -4.53\\
   \hline
   H2 & 4.27 & 1.05 \\
   \hline
   S & 4.27 & 1.06 \\
   \hline 
   E & 4.85 & -0.43 \\
   \hline 
   GH & 5.07 & -1.34 \\
   \hline
   H1 & 5.23 & -1.78 \\
   \hline
   CLC & 5.69 & 3.12\\
   \hline
   C & 20.51 & 7.29 \\
   \hline
   \end{tabular}
   \end{center}
   \caption[Numerical values of the bifurcation points]{$(G,P)$ coordinates of the  different bifurcation and special points of the Jansen and Rit's system}
   \label{table:JansenBifsGP}
\end{table}
\newpage
The behaviors in the eight zones are very similar. The behavior in the zone (A) is similar to the behavior of the system in the zone (B) in figure \ref{fig:FullGlobal}, zone (B) to zone (C) of figure \ref{fig:FullGlobal}, zone (C) to zone (D) of figure \ref{fig:FullGlobal}, zone (D) to (E), (E) to (F) and (F) to (G) and (G) to (H). The only new kind of behavior corresponds to zone (H) where there exist no saddle node bifurcation but a supercritical Hopf bifurcation which corresponds to the existence of a cycle. 

\subsection{Effect of the delay ratio $d$}
\begin{figure}
 \begin{center}
 \includegraphics[width=.8\textwidth]{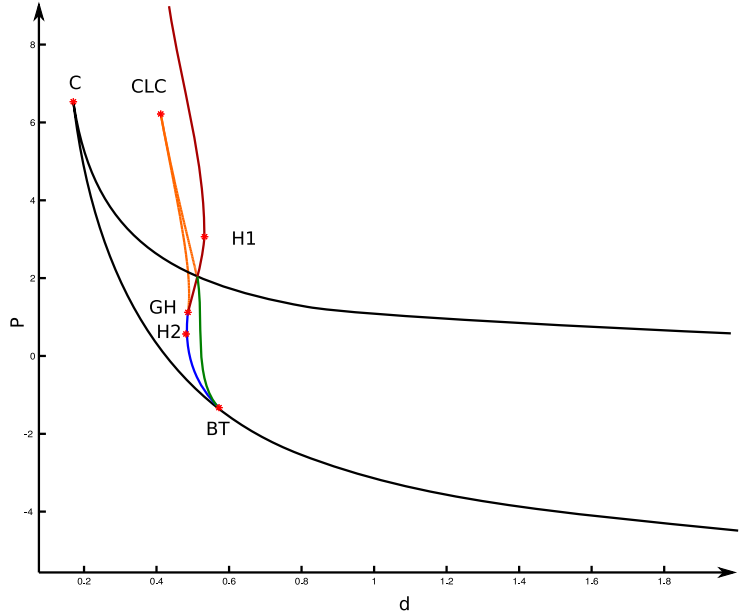}
 \end{center}
 \caption[Jansen's bifurcations in $d$ and $P$]{Codimension 2 bifurcations in Jansen and Rit's model with respect to the delay ratio  $d$ and the input parameter $P$.}
 \label{fig:ZonesdP}
\end{figure}

The delay ratio $d$ features the same bifurcation structure as in the previous two cases with a cusp, a Bogdanov-Takens, a Bautin and a cusp of limit cycles bifurcations, together with saddle-homoclinic and saddle-node homoclinic bifurcations. The bifurcation diagram is shown in figure \ref{fig:ZonesdP}.

\subsection{Sensitivity to the connection probability parameters}

The structure of cycles we studied is unaffected by changes in $\alpha_1$, $\alpha_3$ and $\alpha_4$. For the parameter $\alpha_2$, the picture is slightly different. Indeed, the bifurcation diagram presents a codimension three bifurcation corresponding to the cusp case of the degenerate Bogdanov-Takens bifurcation. With respect to these parameters, the cusp bifurcation point is also a point of a Andronov-Hopf bifurcation. Note that at this very point, the bifurcation diagram is very close to the one of Wendling and Chauvel's model (see section \ref{section:WC}) and therefore will locally generate the same behaviors as those of this more complex model.

This codimension three bifurcation corresponds to the values $DBT: \{P=3.236,\;\alpha_2 = 0.365\}$.

\begin{figure}
\begin{center}
 \subfigure[Codim. 2 bifurcations with respect to $\alpha_3$ and $P$]{\includegraphics[width = .3\textwidth]{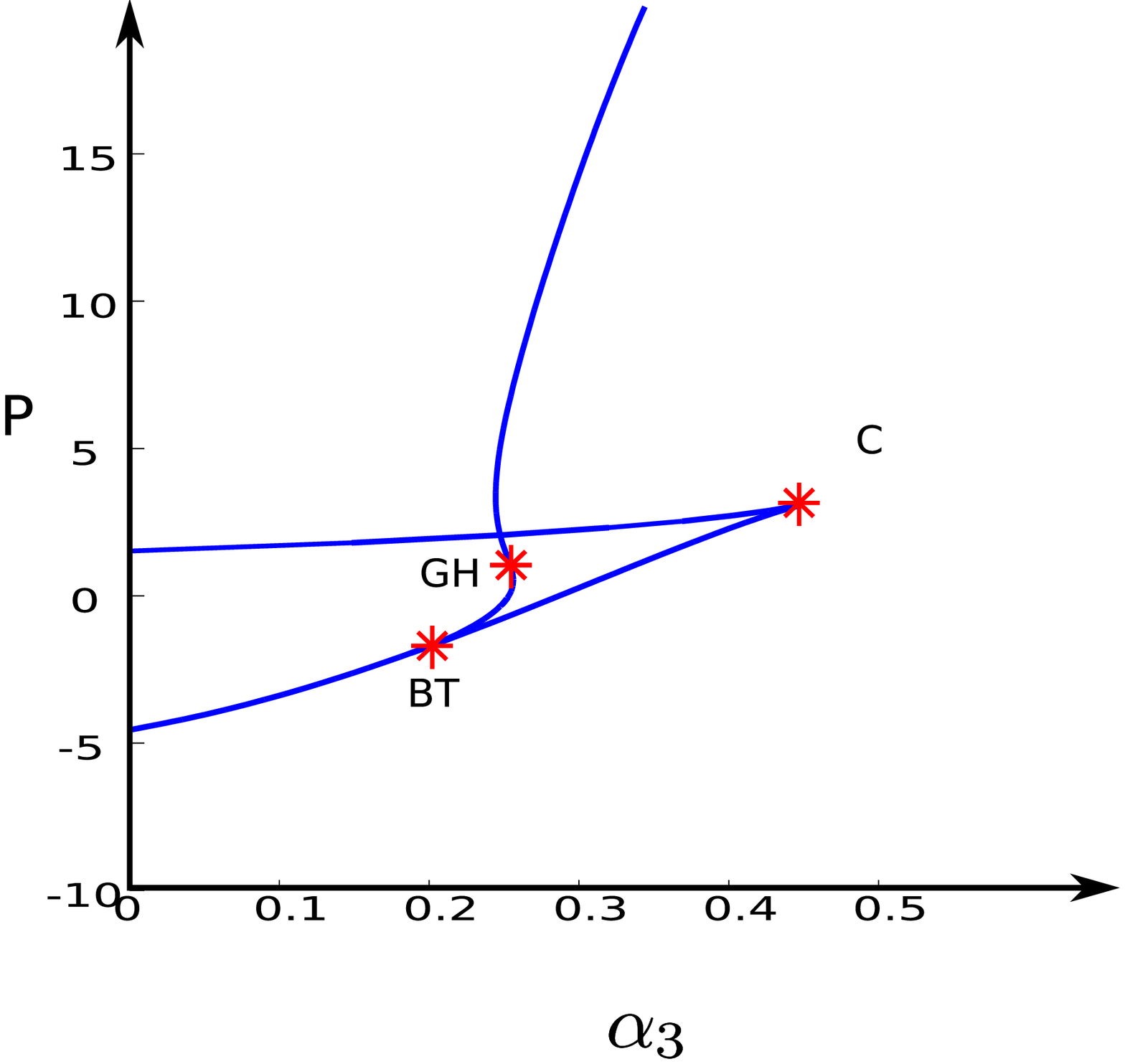}}
  \subfigure[Codim. 2 bifurcations with respect to $\alpha_4$ and $P$]{\includegraphics[width = .3\textwidth]{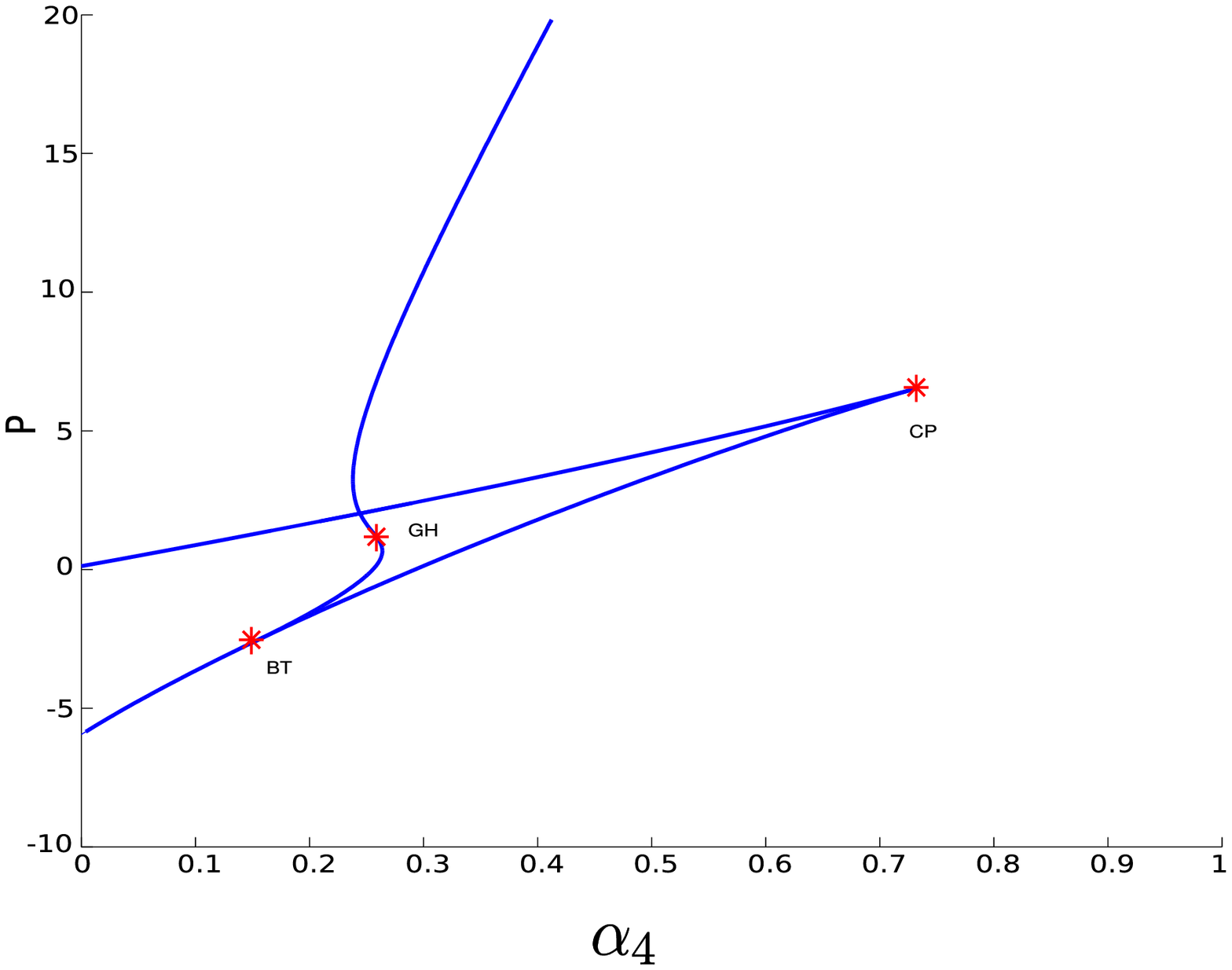}}
 \subfigure[Codim. 2 bifurcations with respect to $\alpha_2$ and $P$]{\includegraphics[width = .7\textwidth]{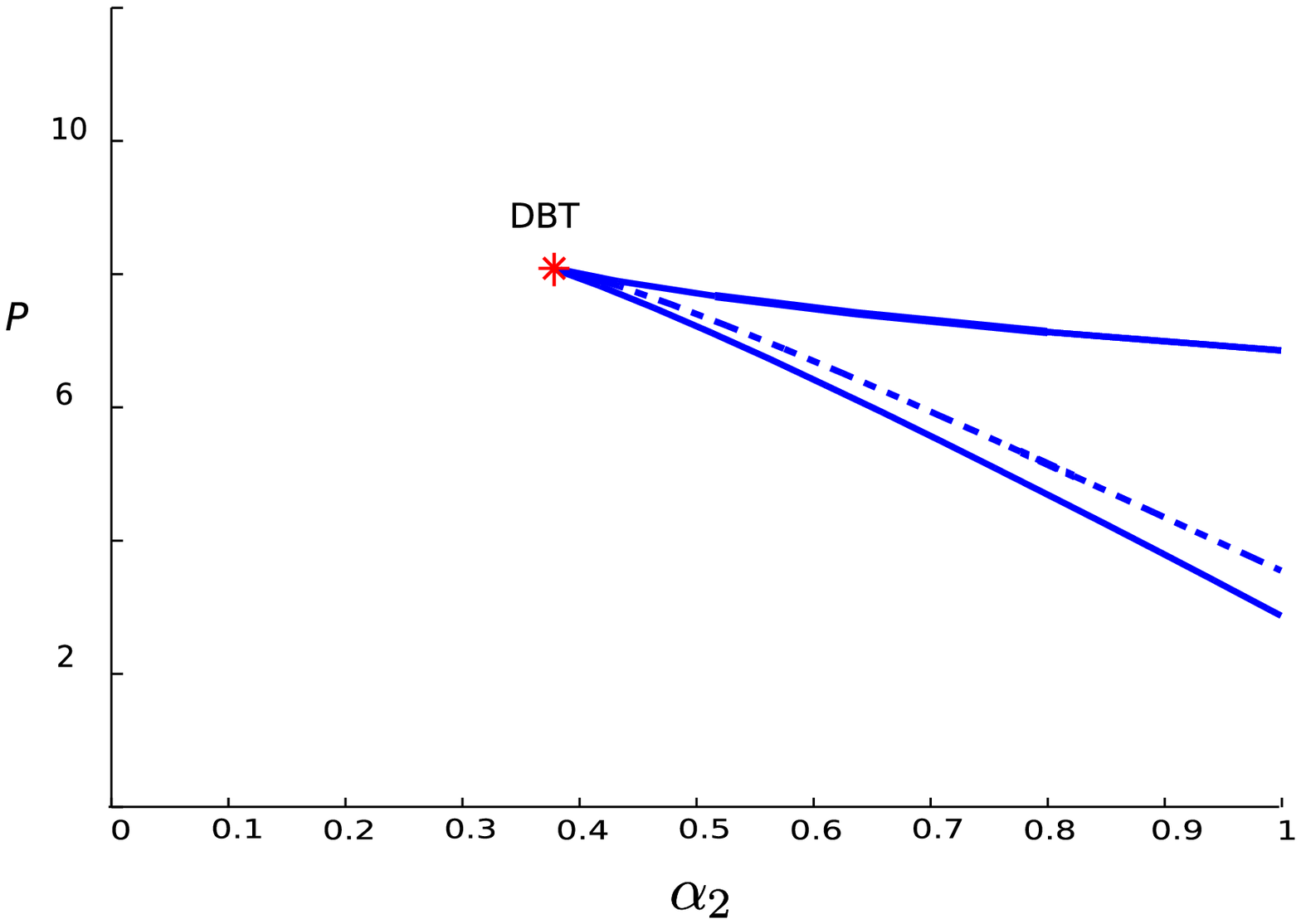}}
\end{center}
\caption[Jansen's bifurcation in $\alpha$s and $P$]{Codimension 2 bifurcations of Jansen and Rit's model with respect to the connection probability parameters $\alpha_1$ (removed because Arxiv refuses large files).,  $i=1,\cdots,4$. In the first three cases no new bifurcations appear: there are only a saddle-node and a Hopf bifurcation manifolds, with possibly Bogdanov-Takens and Bautin bifurcation points. In case (d) we observe a codimension three degenerate Bogdanov Takens (DBT) bifurcation point. }
\end{figure}
%


\section{Analysis of Wendling and Chauvel's model}\label{section:WC}
We now turn to the study of the Wendling and Chauvel's model of hippocampus activity \cite{wendling-chauvel:08}. As stated, this model was initially introduced in order to be able to reproduce the kind of fast activity encountered at the onset of seizures in limbic structures. From many experimental observations, the authors take into account slow and fast currents, which resulted in the model presented in section \ref{ssect:WC}. We show here that codimension one and two bifurcations with the standard parameters for biologically plausible processes do not lead to the emergence of this fast activity. In order to reproduce this kind of activity Wendling and Chauvel used  changes in the amplitude of the PSPs. We present the mathematical emergence of this fast activity via a bifurcation analysis, and discuss the effect of the noise. 

\subsection{Fixed points of the model}
Wendling and Chauvel's dynamical system is even more intricate than Jansen and Rit's, since it has ten dimensions, sigmoidal nonlinearities and component mixing. Nevertheless in this case again, the fixed points of the system can be parametrized as a function of the state variable $X=Y_1-Y_2$. This manifold is given by the following set  equations:

\begin{equation*}
  \begin{cases}
  Y_0 &= jS(X) \\
  Y_2 &= \frac{j G_1 \alpha_4}{d_1} S(\alpha_3 j S(X)) \\
  Z   &= -\frac{j G_1 \alpha_6}{d_1} S(\alpha_3 j S(X)) + \alpha_5 j S(X) \\ 
  Y_3 &= \frac{j G_2 \alpha_7}{d_2} S(-\frac{j G_1 \alpha_6}{d_1} S(\alpha_3 j S(X)) + \alpha_5 j S(X) )\\
  P &=  X + Y_2(X) + Y_3(X)-\alpha_2 j S(\alpha_1 jS(X) ) \\
  Y_i &=0 \qquad \forall i \in \{5 \ldots 9\}
  \end{cases}
\end{equation*}

The input firing rate $P$ is an important parameter. This is why we  always consider bifurcations with respect to this parameter. 

\subsection{Bifurcations of the Wendling and Chauvel model}
\subsubsection{A degenerate codimension three Bogdanov-Takens bifurcation}
We first consider, similarly to Jansen and Rit's case, the codimension two bifurcations with respect to the input $P$ and the total connectivity parameter $j$. 

\begin{figure}
 \begin{center}
  \includegraphics[width=.8\textwidth]{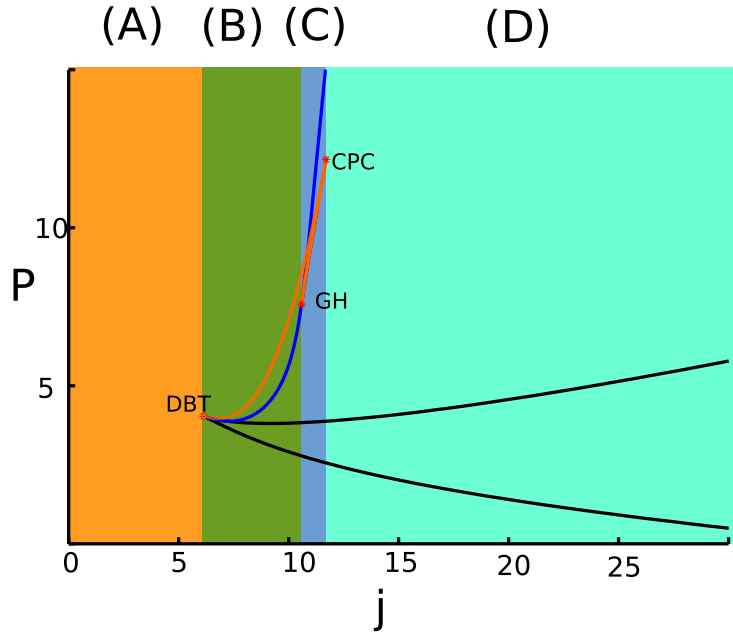}
 \end{center}
\caption[Bifurcations of the Wendling and Chauvel's model]{Bifurcations of the Wendling and Chauvel's model with respect to the parameters $P$ and $j$. The other parameters are given by the values in table \ref{tab:paramsWendling}. We observe a saddle-node bifurcation manifold (plain black curve), a Hopf bifurcation manifold (plain blue curve), a fold bifurcation limit cycles manifold (plain orange curve), a Bautin bifurcation, and a codimension three degenerate Bogdanov-Takens bifurcation (DBT).}
\label{fig:WCBifs}
\end{figure}

We numerically observe that for the original values of Wendling and Chauvel's model, the system features a saddle-node and a Hopf bifurcation manifolds, a cusp, Bautin a codimension three bifurcation: a degenerate Bogdanov-Takens bifurcation. Figure \ref{fig:WCBifs} shows the corresponding bifurcation diagram which is split into four zones labelled from A to D.

\begin{table}
\begin{center}
  \begin{tabular}{|l|l|l|}
  \hline
   Point & j & P \\
   \hline 
   DBT & 6.13 & 4.03 \\
   \hline
   GH & 10.59 & 7.59\\
   \hline
   CLC &  11.71 & 12.15\\
   \hline
   \end{tabular}
   \end{center}
   \caption[Numerical values of the bifurcation points]{$(j,P)$ coordinates of the  different bifurcation and special points of the Wendling and Chauvel's model.}
\end{table}

\begin{enumerate}
 \item In this zone, the system does not present any bifurcation nor any cycle, and the system returns to equilibrium for any initial condition and any input rate.
\begin{figure}
 \centerline{\includegraphics[width=.5\textwidth]{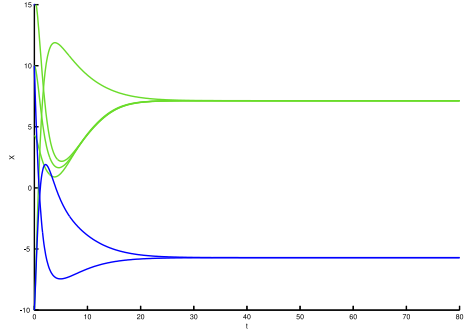}}
\caption{Wendling-Chauvel's model in case (A): unconditional convergence to a unique fixed point, j=5, blue P=-5, green: P=10.}
\label{fig:WCA}
\end{figure}

\item In this zone, the system features two saddle-node bifurcations and a subcritical Hopf bifurcation generating unstable limit cycles. The family of cycles undergoes a fold of limit cycles, and the system displays slow oscillations of large amplitude, corresponding to epileptic-like activity.
\begin{figure}
\begin{center}
 \subfigure[Case (B): Bifurcation diagram] {\includegraphics[width=.4\textwidth]{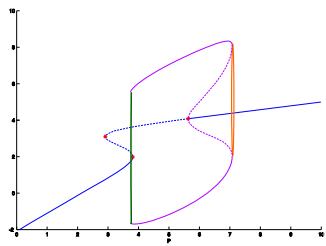}}\\
 \subfigure[Case (B): Periods of the limit cycles] {\includegraphics[width=.4\textwidth]{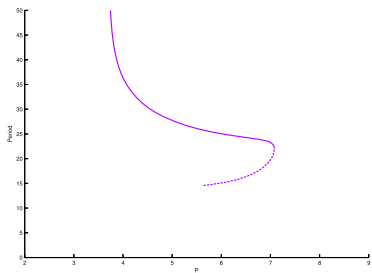}}\quad
 \subfigure[Case (B): Time trajectories] {\includegraphics[width=.4\textwidth]{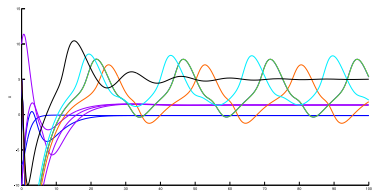}}\\
\end{center}
\caption{Wendling-Chauvel's model in case (B): Epileptic Spikes. Parameters: j=8. (c): blue: P=2, purple: P=3.5, orange P=4, green: P=5.5, cyan P=7, black P=9.}
\label{fig:WCB}
\end{figure}

\item In this case, the system presents alpha activity together with epileptic spikes of low frequency and high amplitude because of the presence of two fold bifurcations of limit cycles. 
\begin{figure}
\begin{center}
 \subfigure[Case (C): Bifurcation diagram] {\includegraphics[width=.49\textwidth]{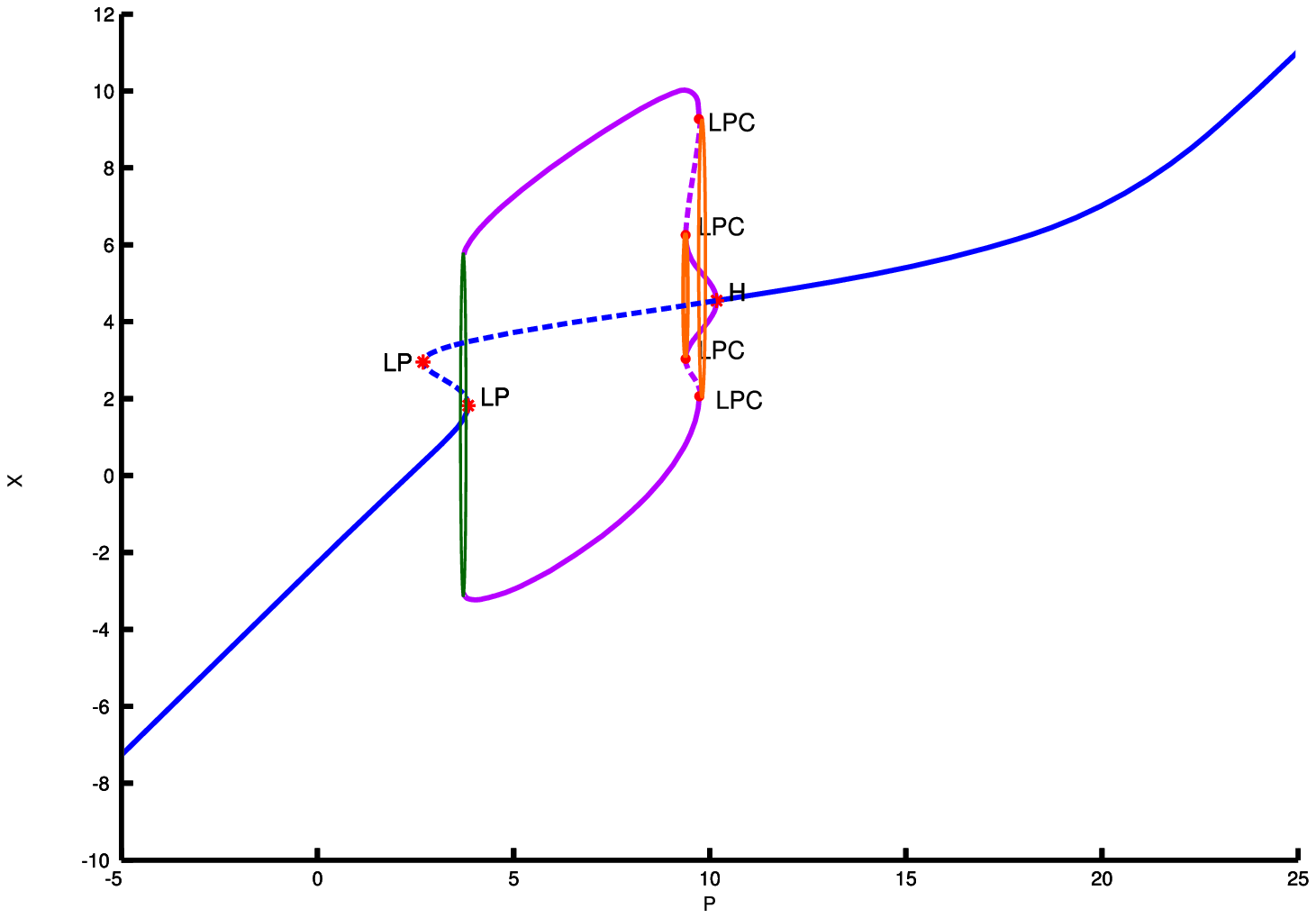}}\quad
 \subfigure[Case (C): Periods of the limit cycles] {\includegraphics[width=.49\textwidth]{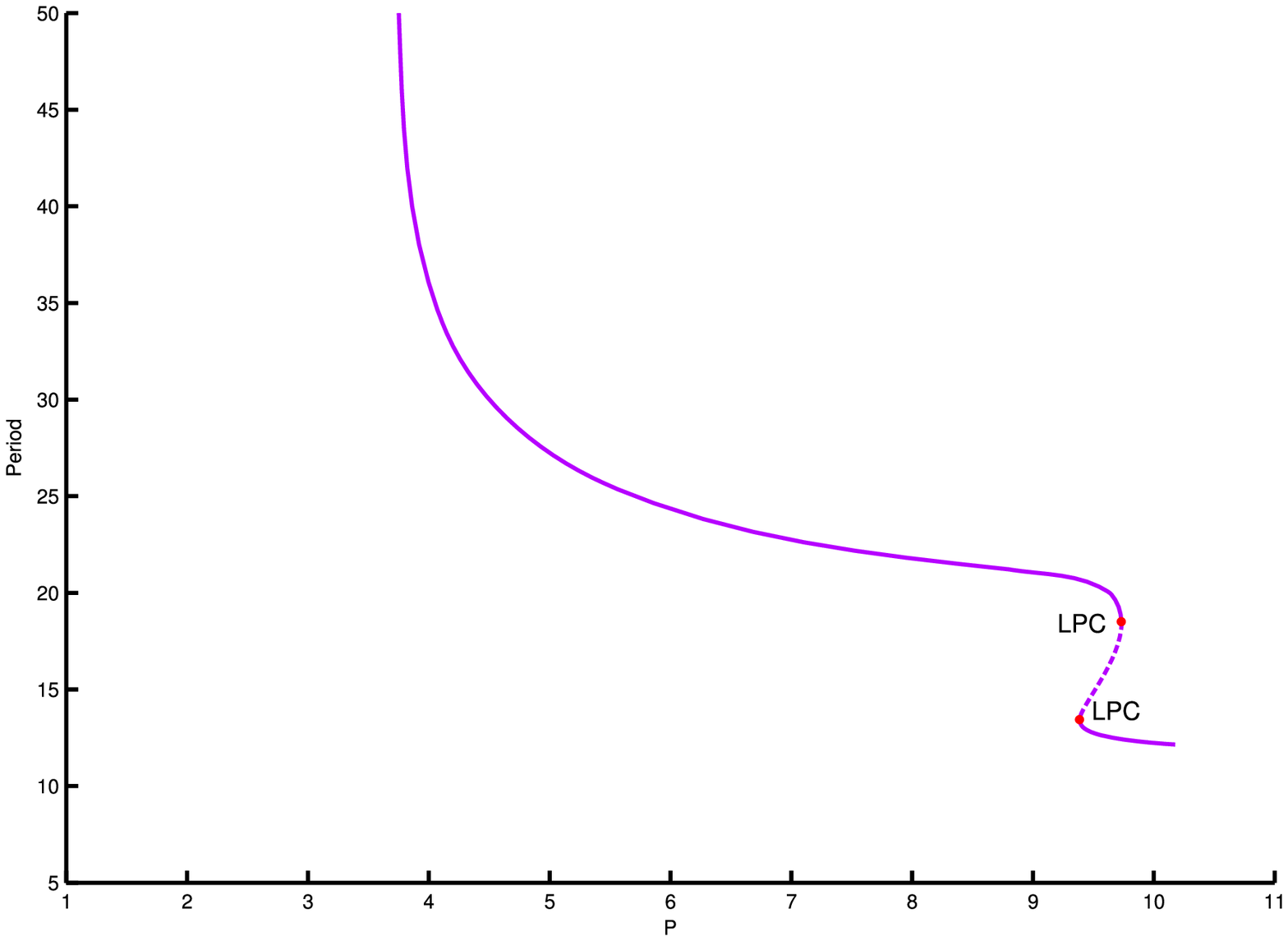}}\\
 \subfigure[Case (C): Time trajectories] {\includegraphics[width=.4\textwidth]{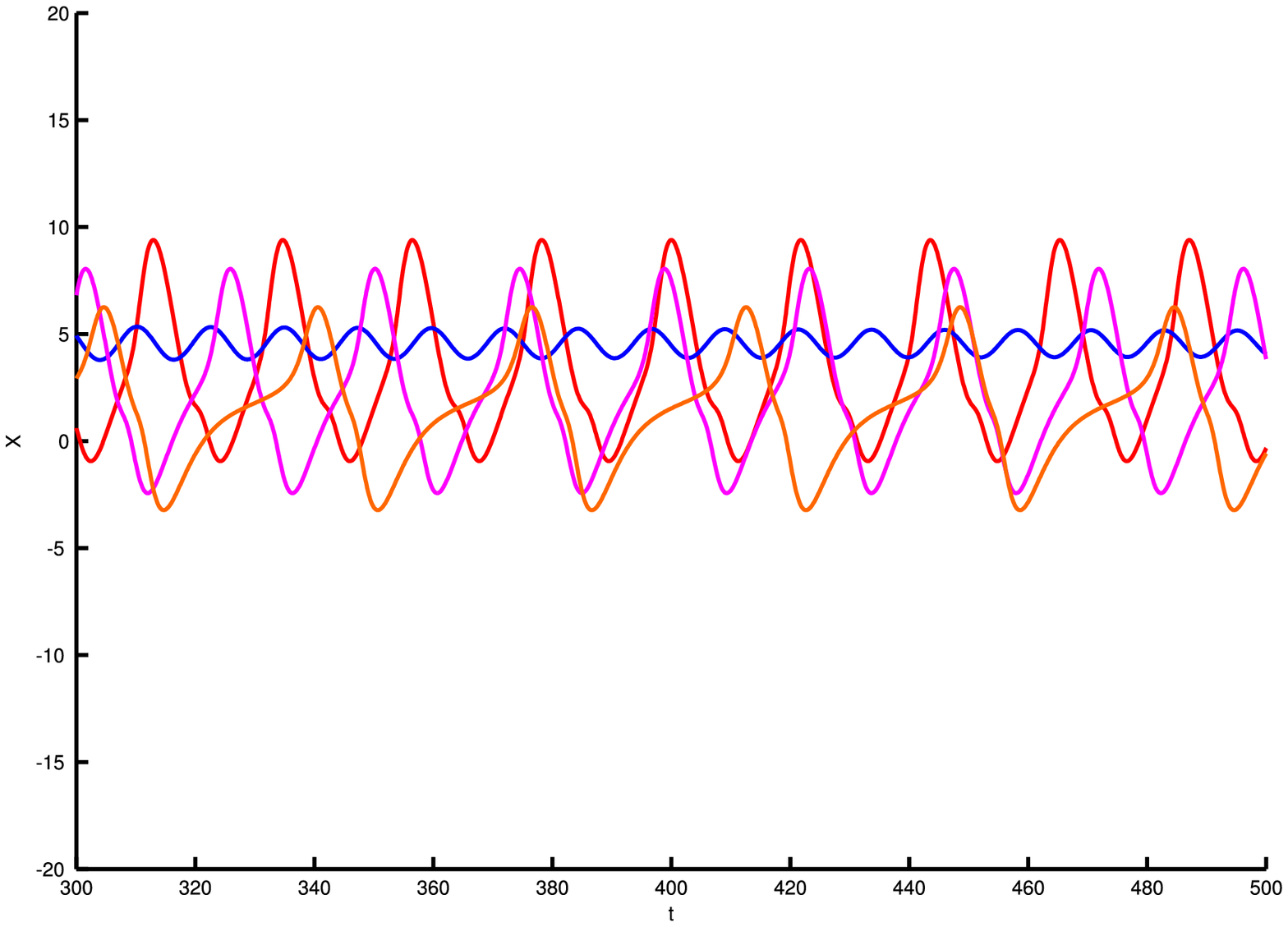}}
\end{center}
\caption{Wendling-Chauvel's model in case (C): alpha and epileptic waves. Blue and red: P=9.7, pink P=7, orange P=5.}
\label{fig:WCC}
\end{figure}

\item Case (D) corresponds to the case where the system goes continuously from alpha activity to slow waves that can be interpreted as sleep waves. 

\begin{figure}
\begin{center}
 \subfigure[Case (D): Bifurcation diagram] {\includegraphics[width=.4\textwidth]{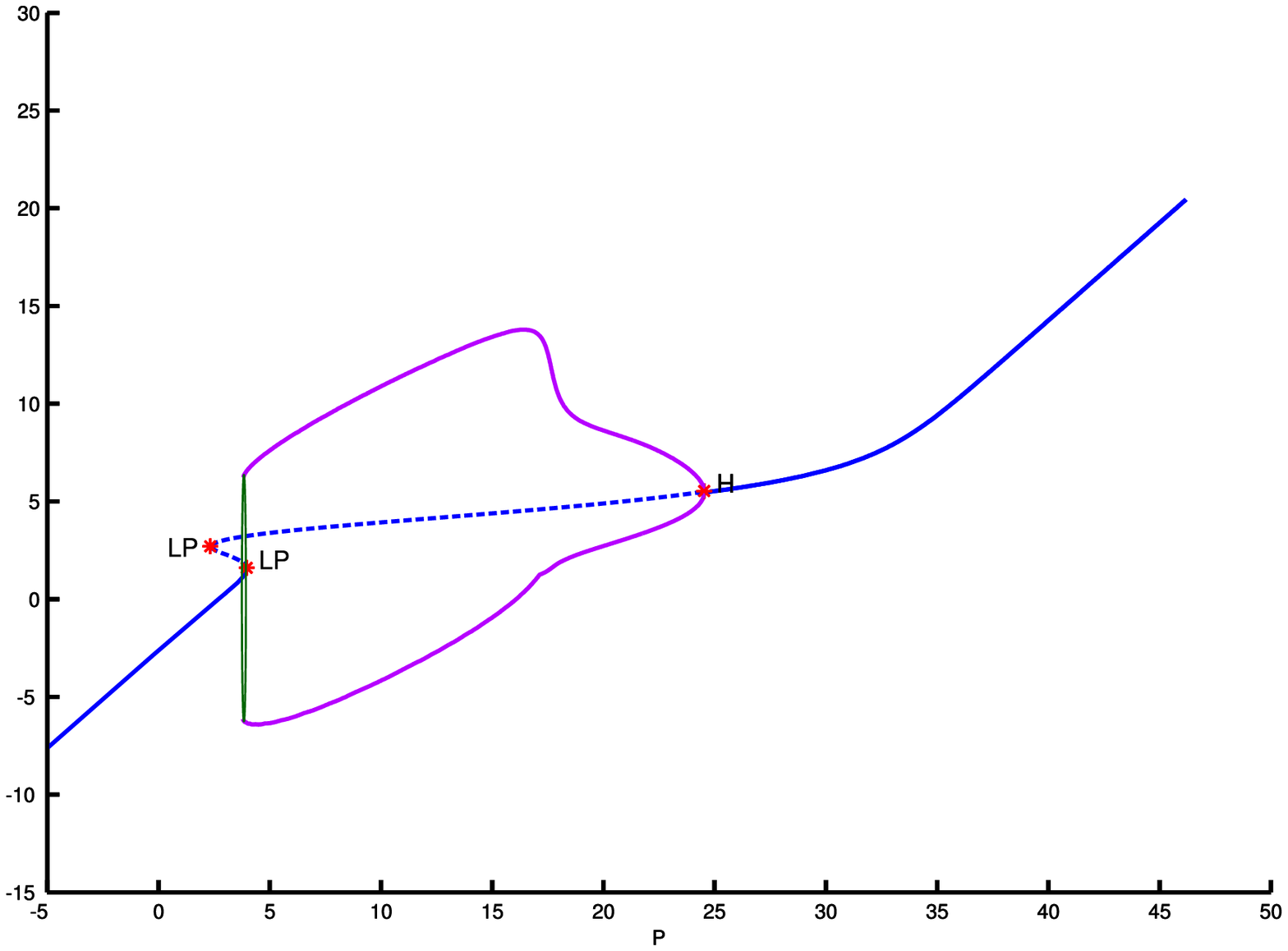}}\quad
 \subfigure[Case (D): Periods of the limit cycles] {\includegraphics[width=.4\textwidth]{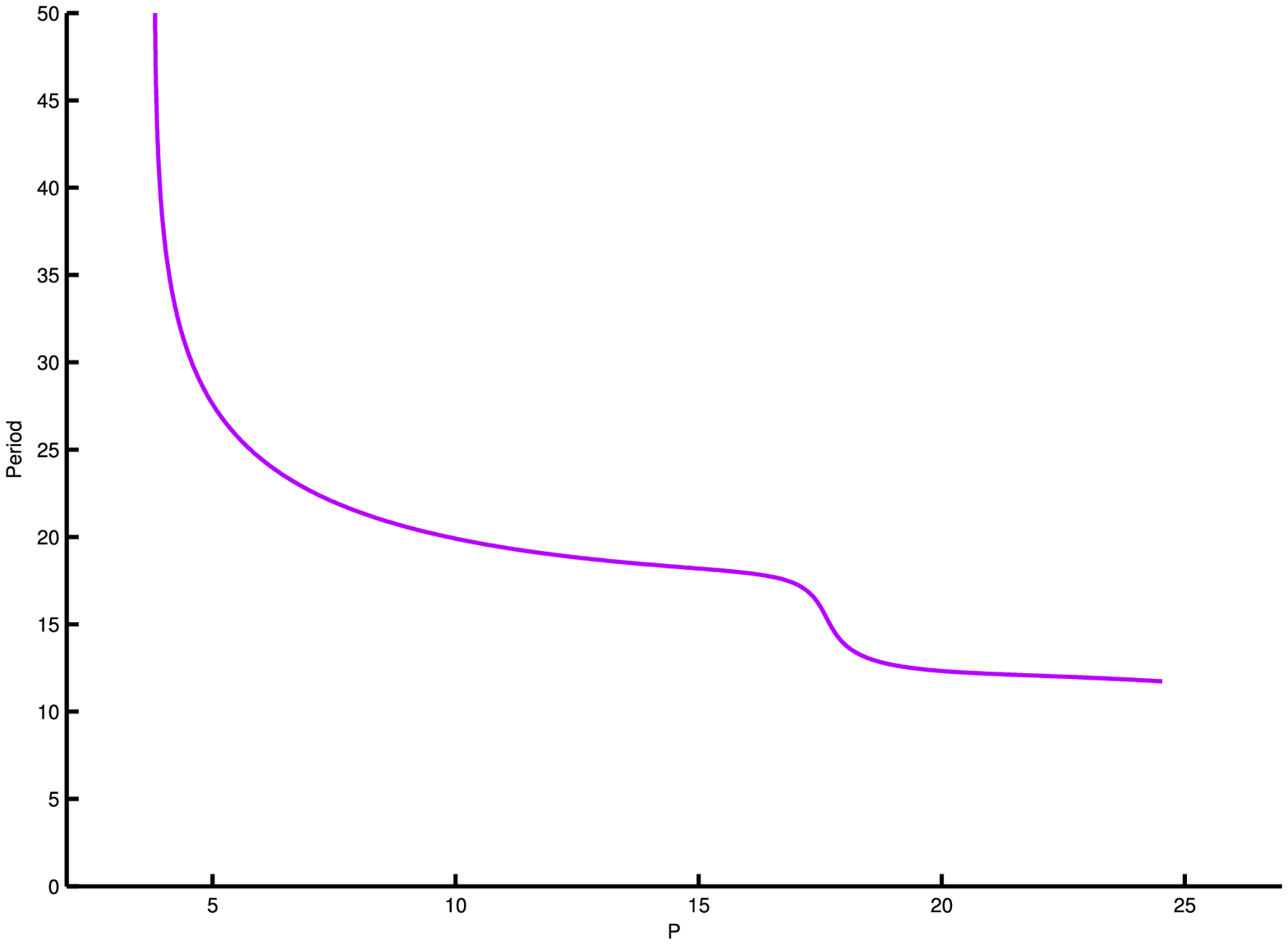}}\\
 \subfigure[Case (D): Time trajectories] {\includegraphics[width=.4\textwidth]{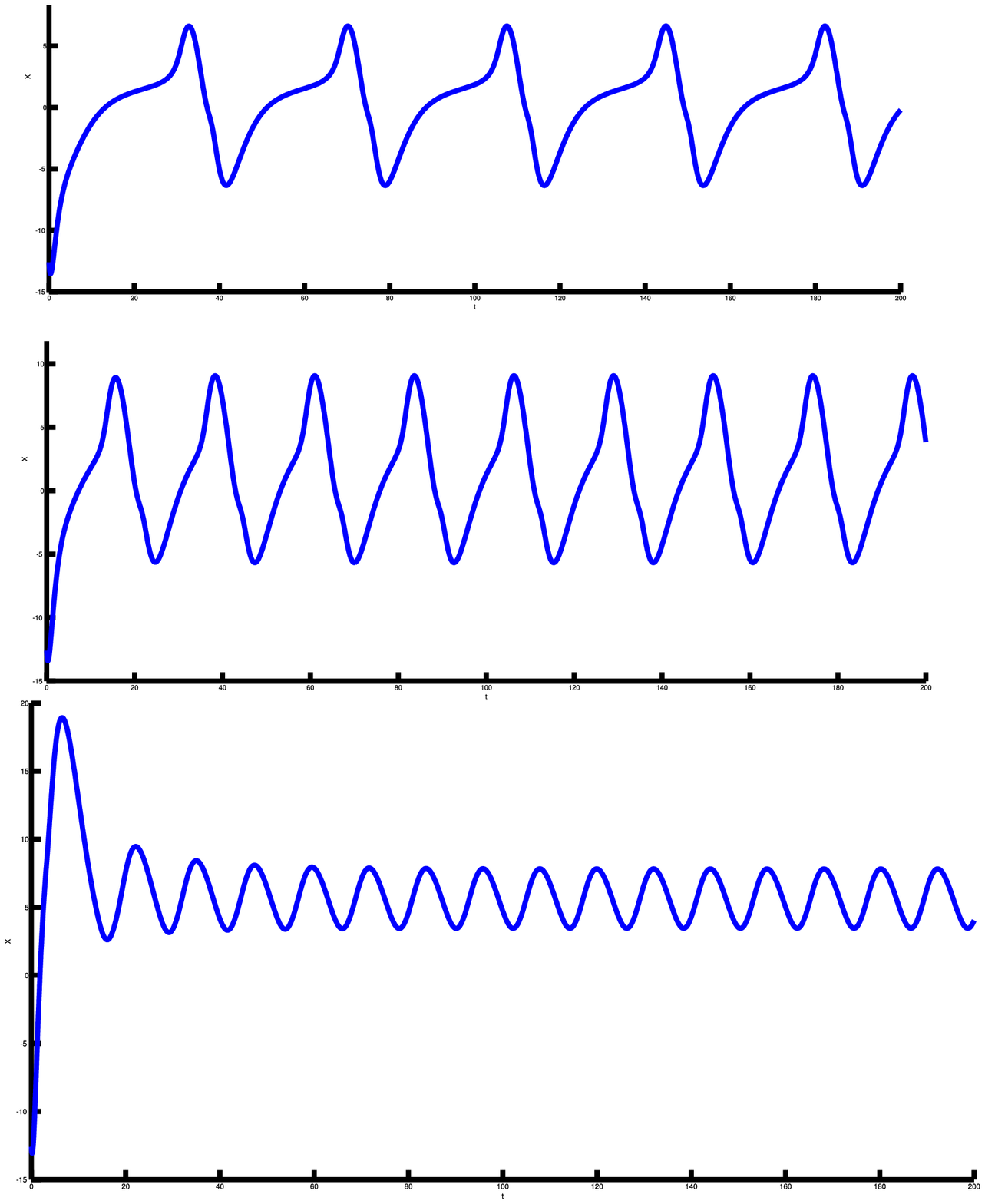}}
\end{center}
\caption{Wendling-Chauvel's model in case (D): Sleep waves, j=13, P=5, 15 and 22. }
\label{fig:WCD}
\end{figure}

\end{enumerate}

We observe that, so far, i.e. in the case where the parameters are set using biophysical values, Wendling and Chauvel's system only features behaviors that are also featured by Jansen and Rit's model. In particular it does not display the kind of fast activity that motivated its development. 

\subsection{Behaviors of Wendling and Chauvel's model and bifurcations}
In \cite{wendling-chauvel:08}, the authors provide a partition of the parameter space $(A,B,C)$ of the original system corresponding to different behaviors of the system, and distinguish between normal activity, narrowband activity, sporadic spikes and fast activity. The partition has been computed simulating different EEG traces for a Brownian input rate. 
The reduction of the model proposed in section \ref{ssect:WC} shows that in the space $(A,B,C)$, bifurcations can be described via a bifurcation analysis in the space $(G_1=B/A, G_2=C/A)$. Nevertheless, one has to be careful that the coefficient $A$ influences the input firing rate since $P=\frac{rA}{Ja}p$. 

\begin{figure}
 \begin{center}
  \subfigure[Small values of $A$: Cycles]{\includegraphics[width=.45\textwidth]{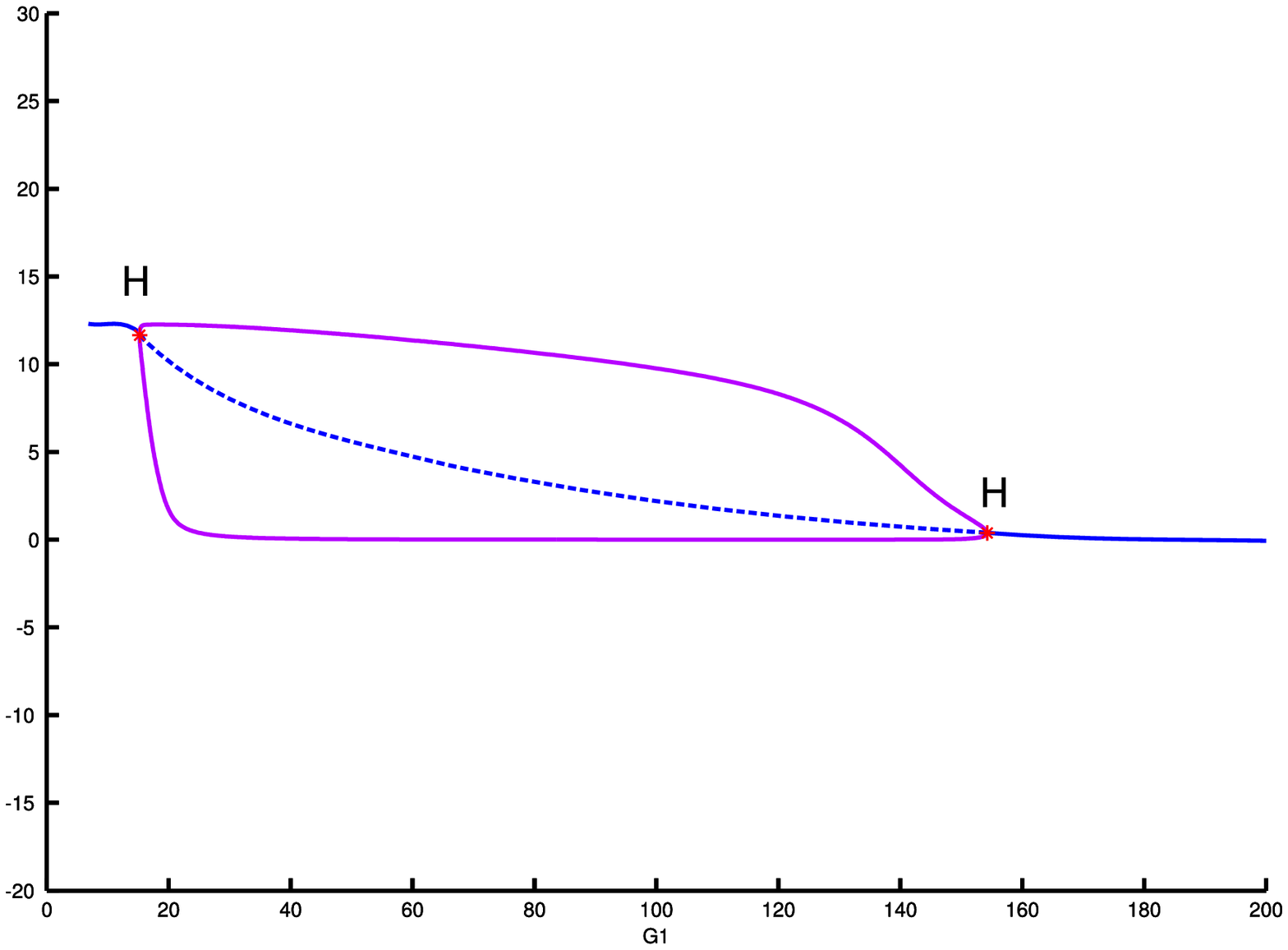}\label{fig:A}}\quad
  \subfigure[Small values of $A$: Periods]{\includegraphics[width=.45\textwidth]{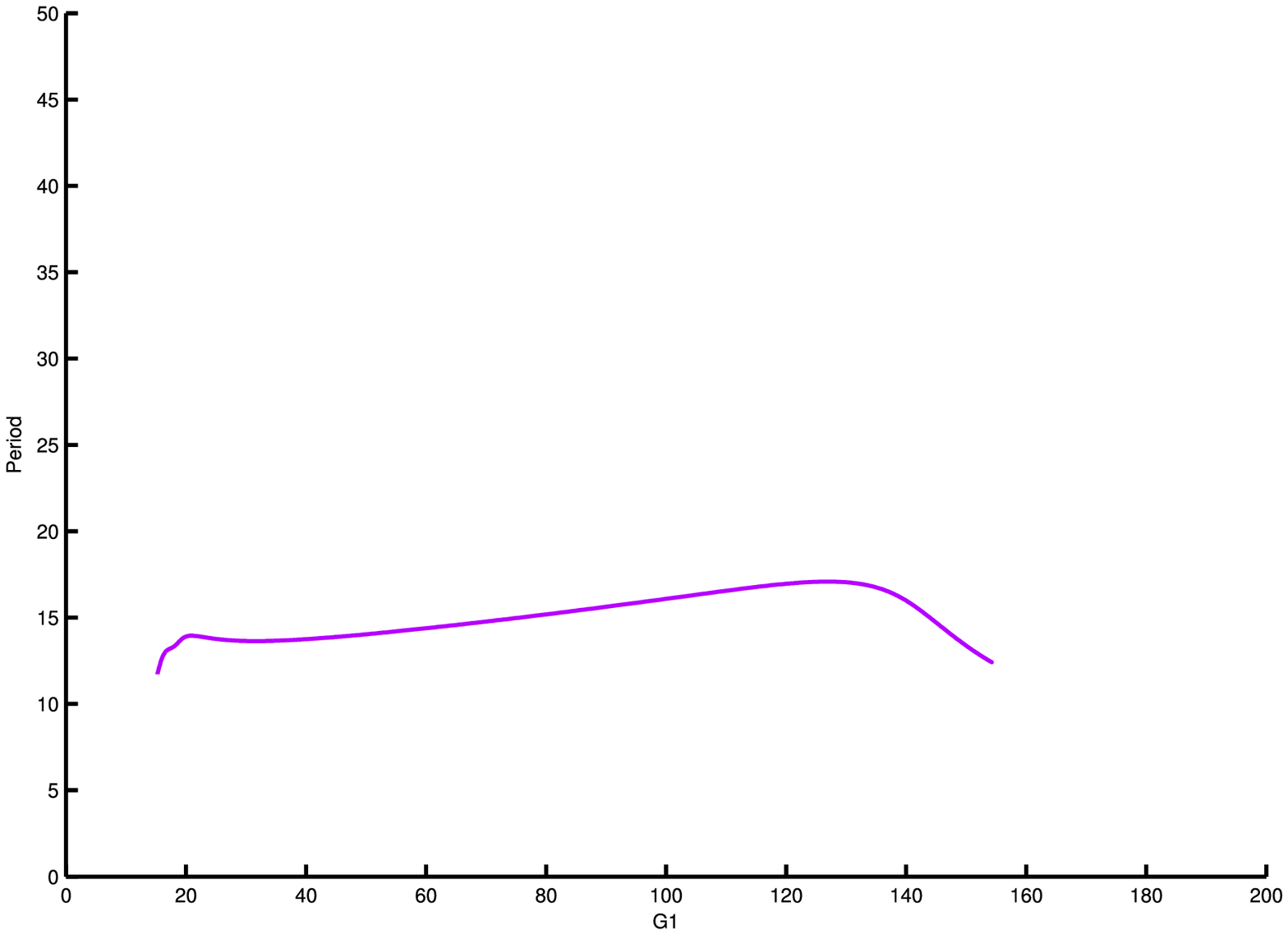}\label{fig:APer}}\\
  \subfigure[Regular values of $A$: Cycles]{\includegraphics[width=.45\textwidth]{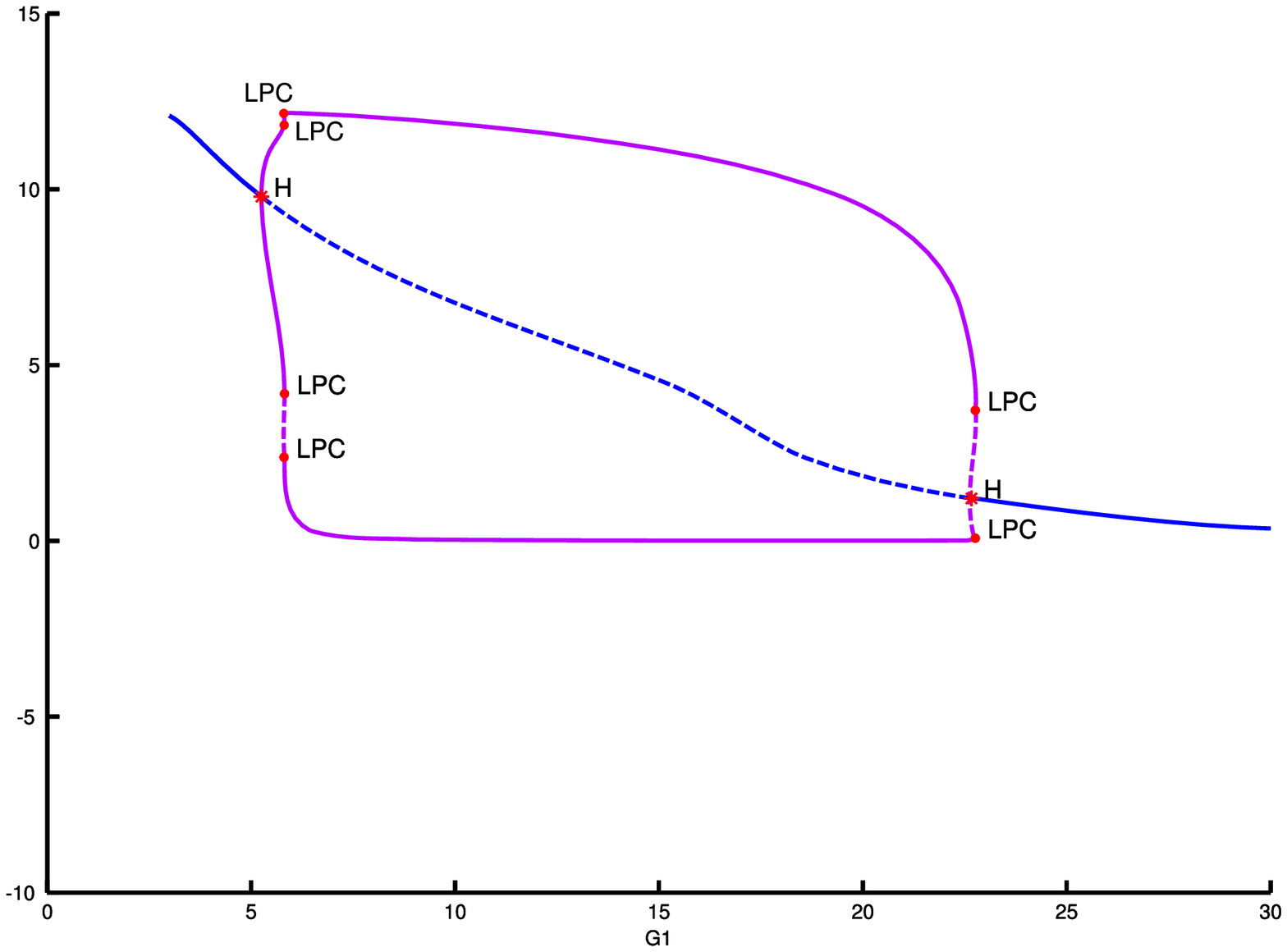}}\quad
  \subfigure[Regular values of $A$: Periods]{\includegraphics[width=.45\textwidth]{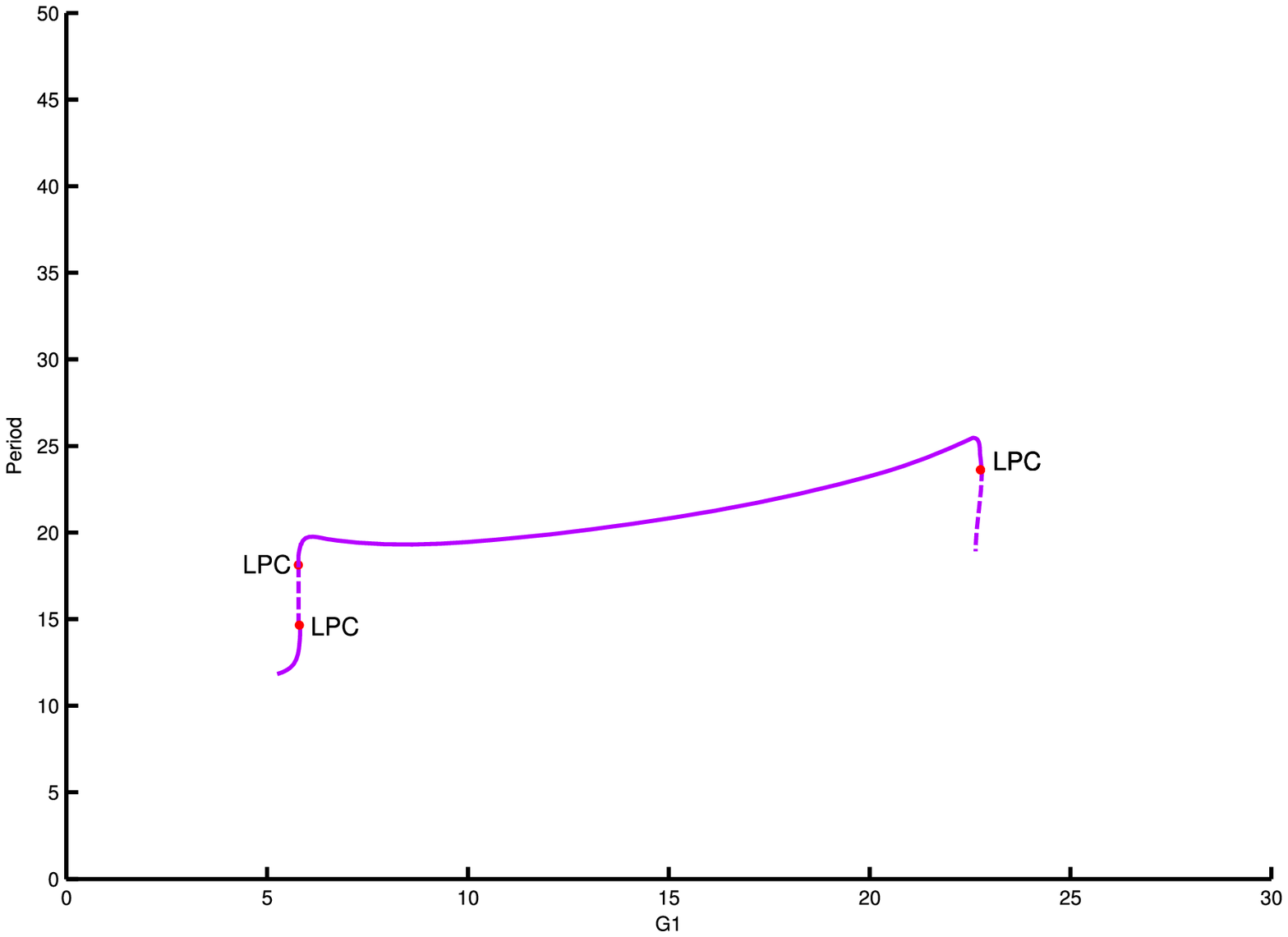}}\\
  \subfigure[Regular values of $A$: Codimension 2 bifurcation diagram with respect to $G_1$ and $G_2$.]{\includegraphics[width=.4\textwidth]{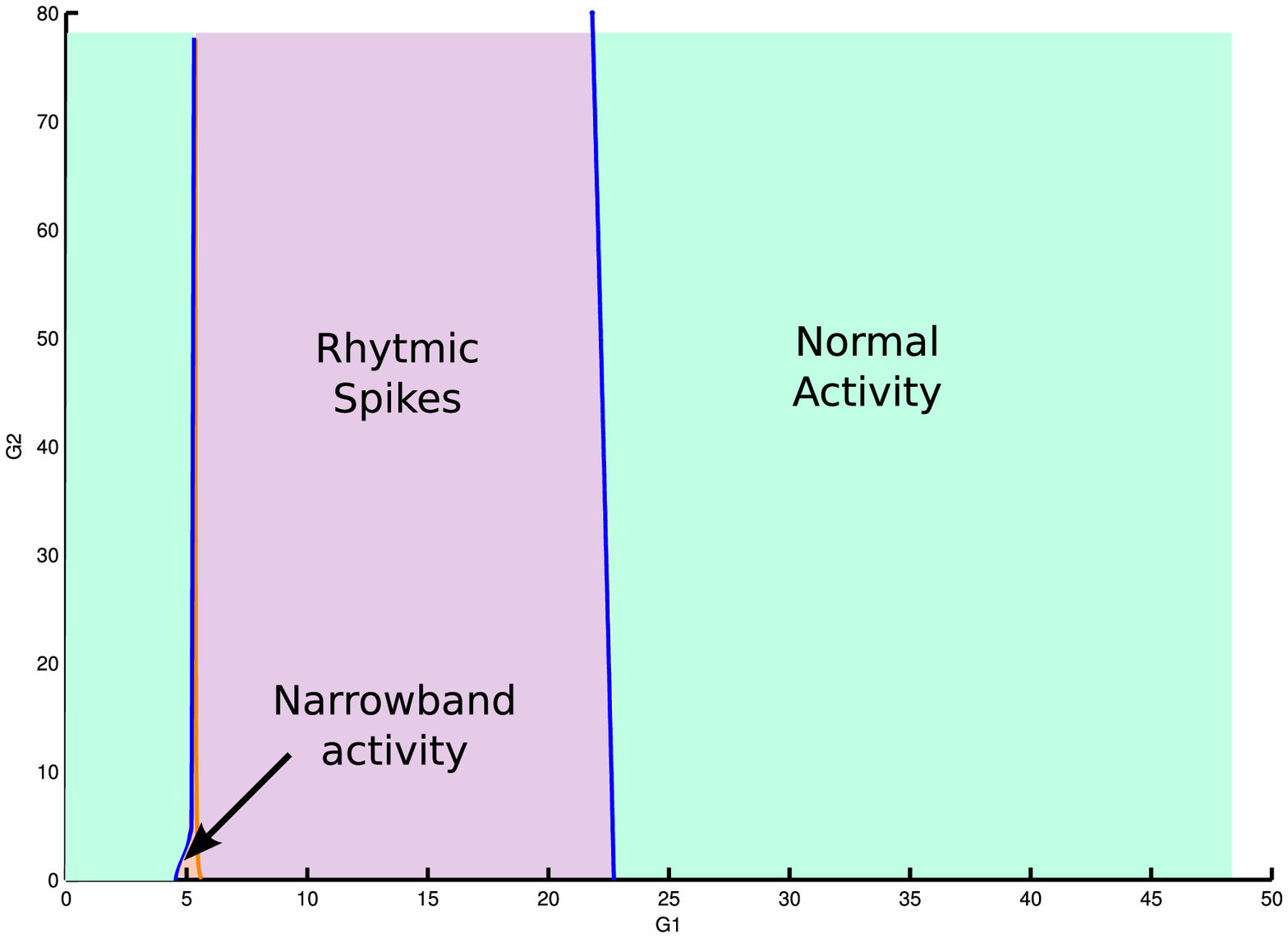}}\quad
  \subfigure[Regular values of $A$: Zoom on Fig (e).]{\includegraphics[width=.4\textwidth]{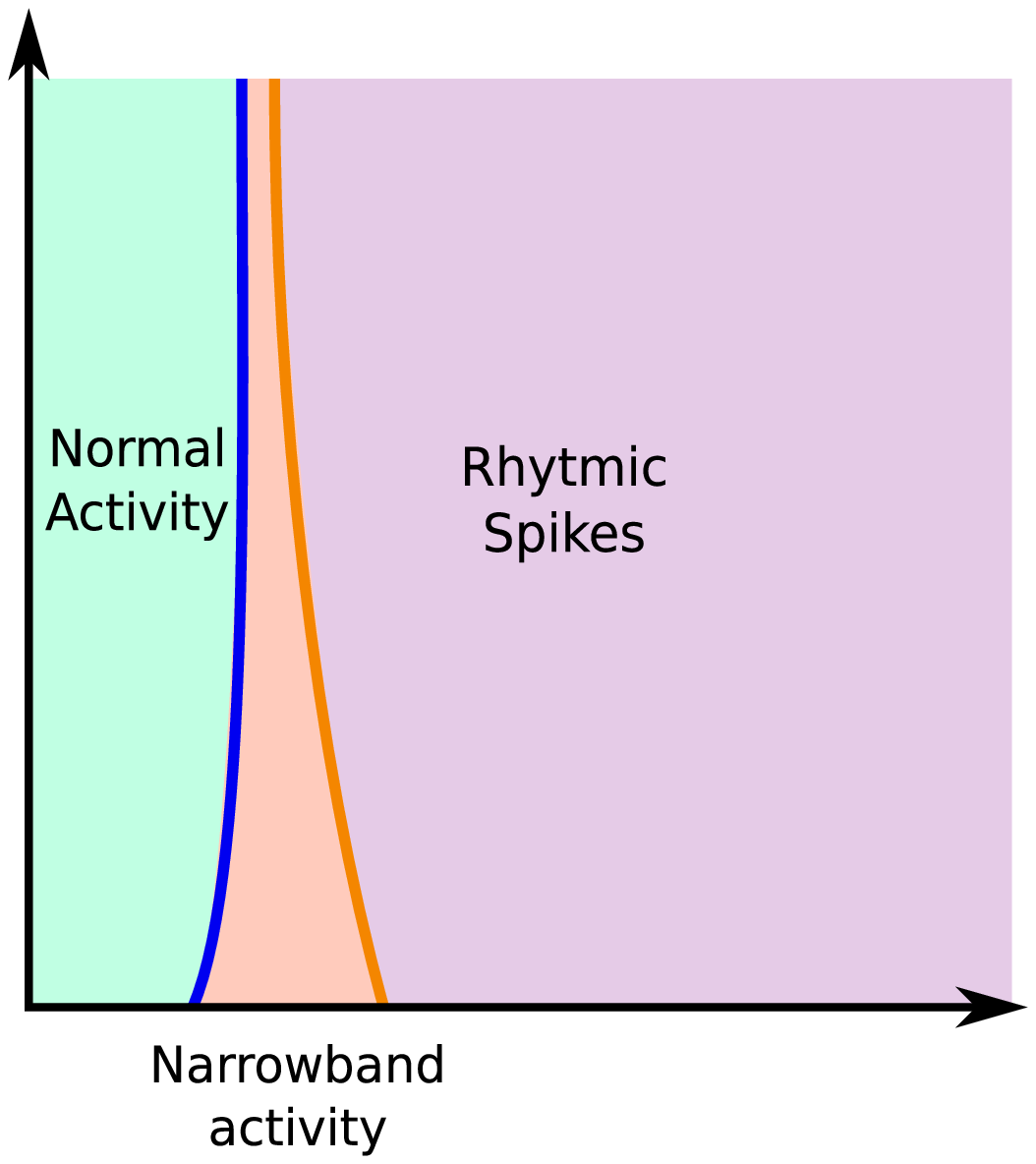}}
 \end{center}
\caption{Behaviours and bifurcations in Wendling and Chauvel's model.}
\label{fig:WendlingParition}
\end{figure}

In the case where $A$ is small, no oscillatory activity is seen in the  region of interest of the parameters: the system is quite far from any bifurcation. For slighter higher values of this parameter, the system presents two Hopf bifurcations which share the same family of periodic orbits leading to narrowband activity $(\alpha,\theta)$, which corresponds to the case $A=3$ in \cite{wendling-chauvel:08}, as represented in figures \ref{fig:A} and \ref{fig:APer}. The fast activity $(\beta,\gamma)$ observed numerically has no deterministic equivalent. It corresponds to a purely stochastic phenomenon: there exists a unique stable fixed point in the region presented by Wendling and Chauvel in \cite{wendling-chauvel:08}, and the Jacobian matrix has complex eigenvalues. The noisy inputs destabilize the system, which moves randomly around the fixed point. The way the system returns to equilibrium corresponds to oscillations at the frequency corresponding to the inverse of the imaginary part of the eigenvalue. This is why the spectrum of the signal generated in the neighboorhood of this fixed point presents a peak in the Fourier spectrum, but this spectrum is quite spread out and the amplitude of the signal is quite low, resembling purely noisy activity.

For larger values of $A$, folds of limit cycle appear. In that case, we observe a small region of narrowband activity $(\alpha,\theta)$, immediately followed by rhythmic epileptic spikes corresponding to a family of limit cycles of large amplitude and low frequency. Sporadic spikes appear in the region where an unstable cycle exists, which is coherent with the interpretation given in the case of Jansen and Rit's model. We also see that it corresponds to a random attraction to a nearby epileptic cycle: sporadic spikes in Wendling and Chauvel's model appear only for regions close to the epileptic activity. Around the Hopf bifurcation, the system features a stable fixed point whose Jacobian matrix has non-real eigenvalues. The imaginary part of the eigenvalue is close to the value corresponding to $\alpha/\theta$ activity near the Hopf bifurcation, and then increases, which corresponds to a transition between narrowband and fast activities. This picture remains valid when increasing further the parameter $A$, as observed in \cite{wendling-chauvel:08}.

\begin{figure}
 \begin{center}
  \subfigure[Codimension two bifurcation diagram]{\includegraphics[width=.45\textwidth]{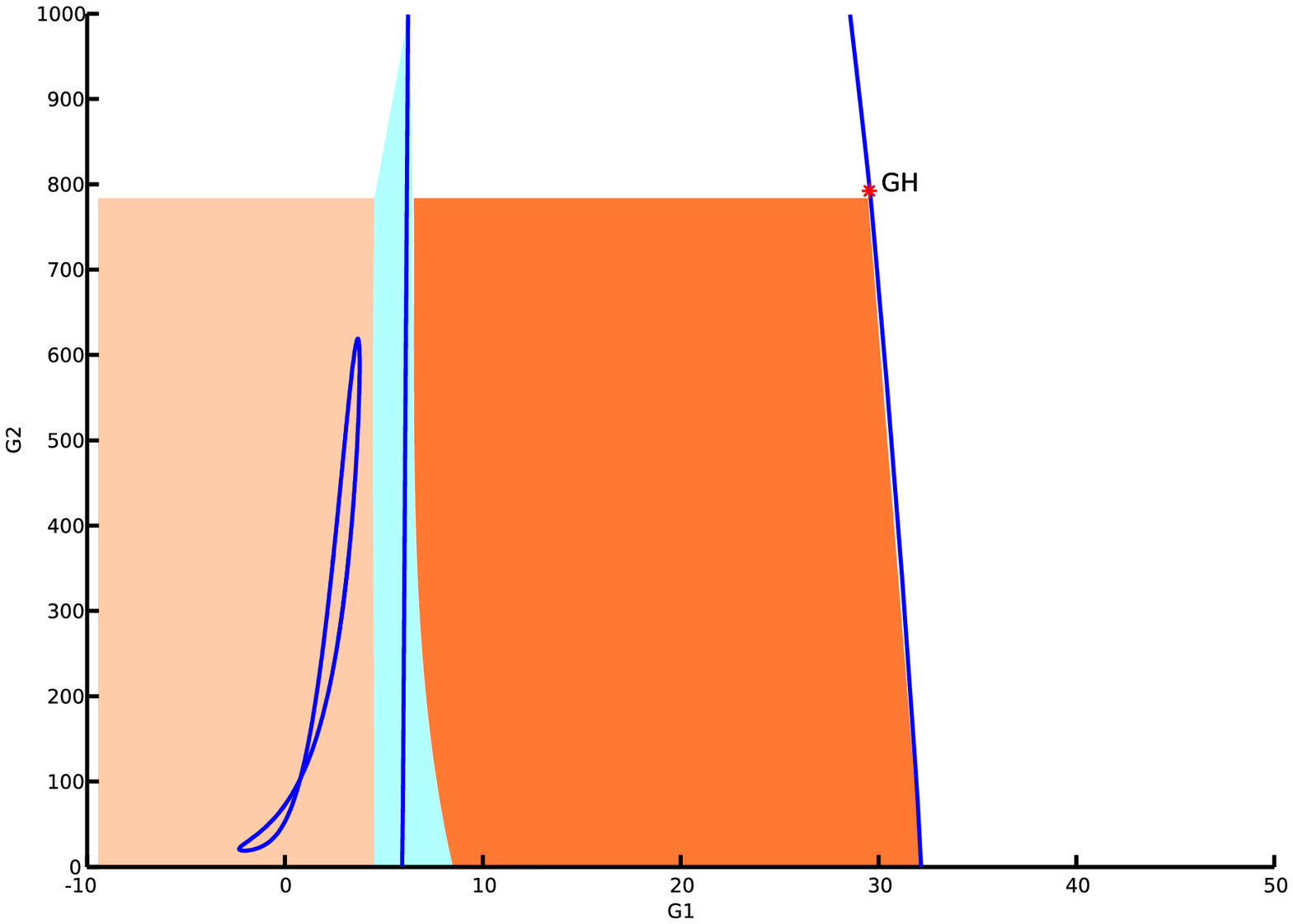}}\quad
  \subfigure[Period of the cycles]{\includegraphics[width=.45\textwidth]{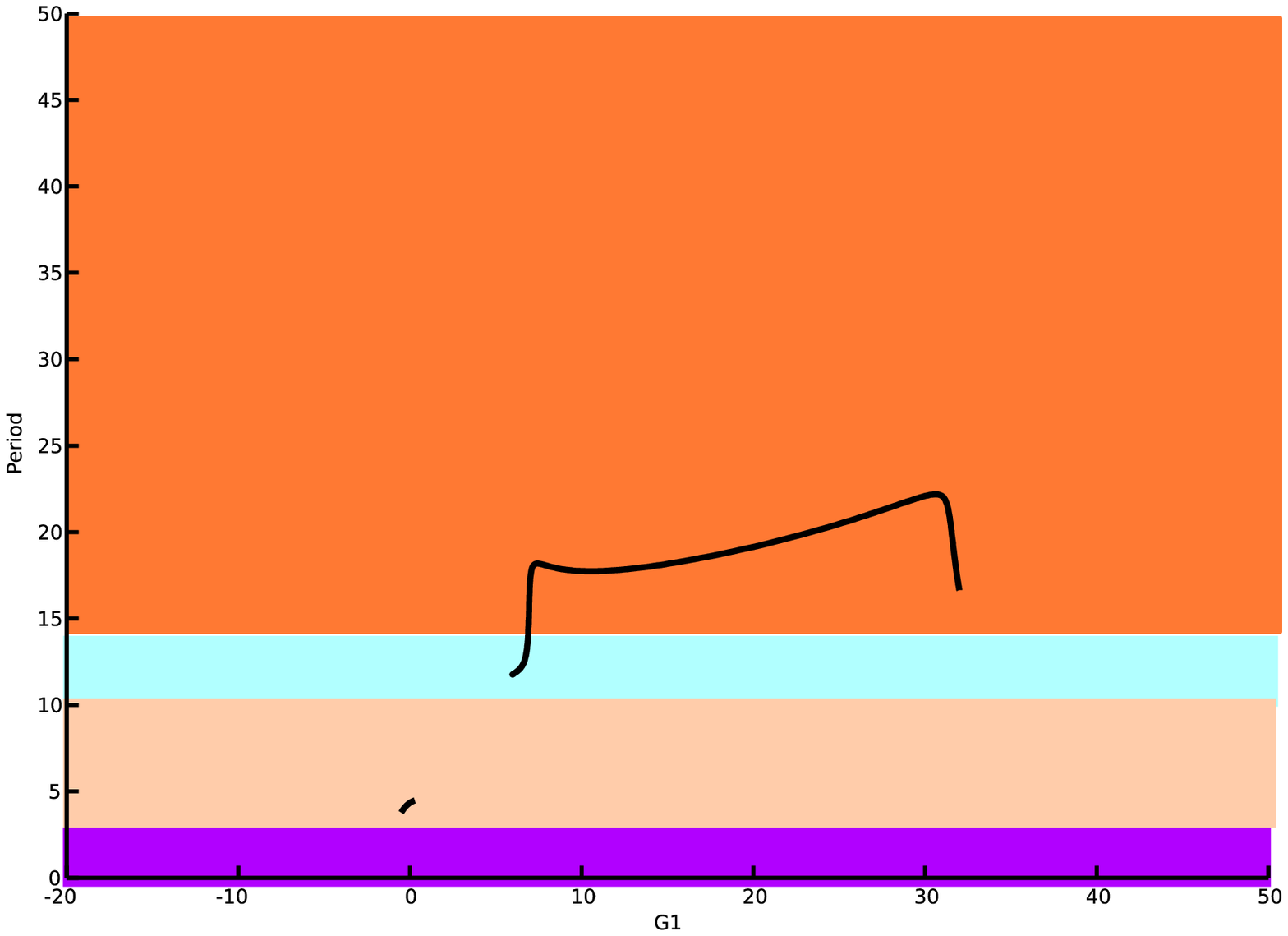}}\\
 \end{center}
\caption{(a): Periods of the limit cycles. Purple band: Very fast waves, Pink: fast activity $(\beta,\gamma)$, blue: narrowband activity $(\alpha, \theta)$, orange: slow activity. (b) Codimension two bifurcation diagram with respect to $G_1$ and $G_2$. We observe numerically the emergence of the different behaviors observed by Wendling and Chauvel. (c) Shape of the cycles in the space $(G_1,Y_0,X)$.}
\label{fig:WendlingStrange}
\end{figure}

For intermediate values of the parameter $A$, we observe deterministic fast activity on a quite strange manifold, as presented in figure \ref{fig:WendlingStrange}, existing for quite a small range of values of the parameter $A$. But the phenomenology stays the same even when deterministic cycles do not exist, as already explained. Note that these cycles do not appear in Wendling and Chauvel's model with the original parameters, since these cycles and these types of rapid period do not appear in the bifurcation diagram in figure \ref{fig:WCBifs}. 


\newpage

\section{Applications: Epileptic seizures dynamics}
Interictal spikes appear at a saddle-node homoclinic bifurcation: for inputs lower than the bifurcation value, the system shows a regular behavior corresponding to the existence of a stable fixed point. Perturbing this equilibrium by taking into account random activity leads to the appearance of non-rhythmic spikes which can be interpreted as interictal spikes. 

\begin{figure}[h]
 \begin{center}
  \subfigure[j=12.285]{\includegraphics[width=.45\textwidth]{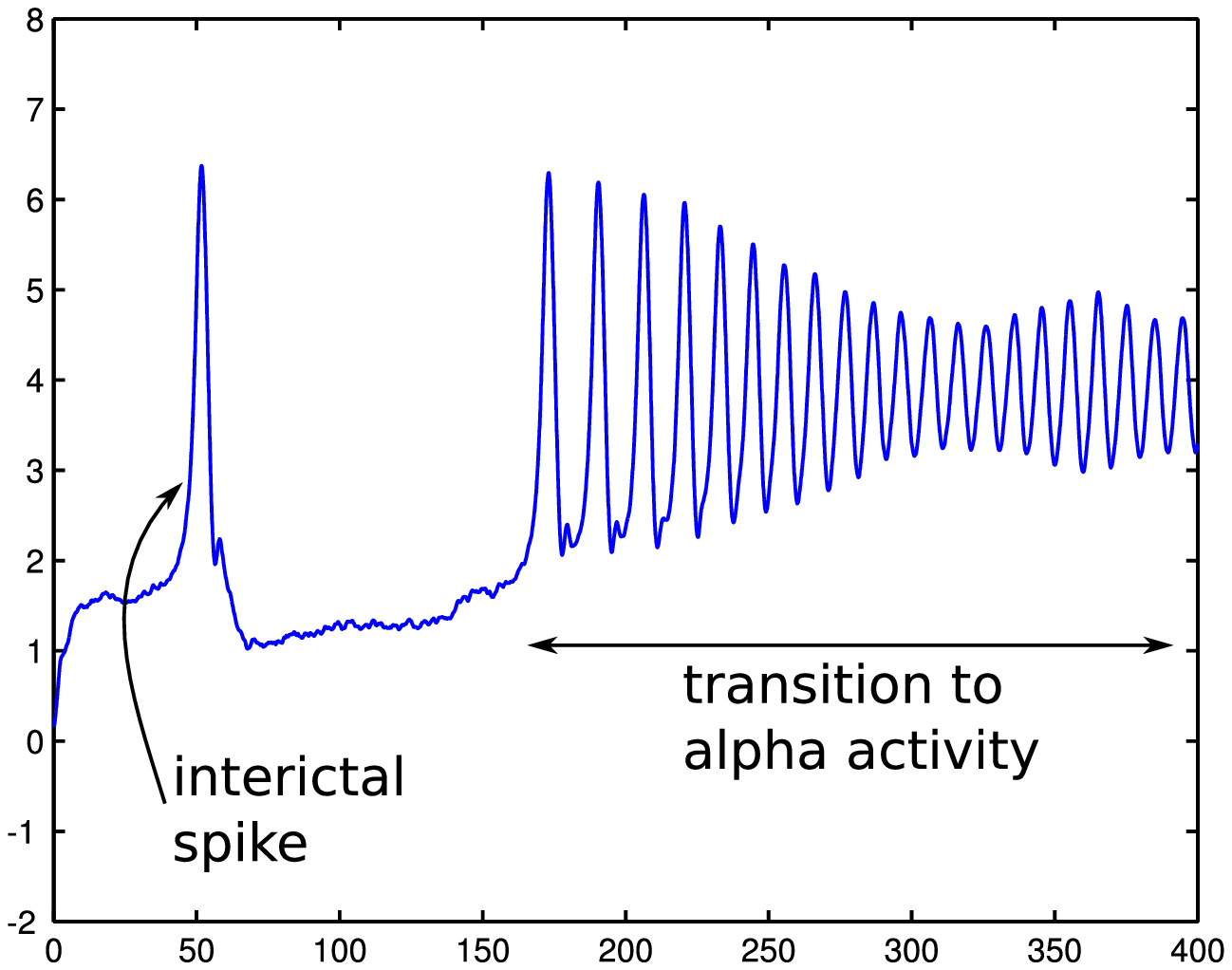}\label{fig:interalpha}}
  \subfigure[j=12.7]{\includegraphics[width=.45\textwidth]{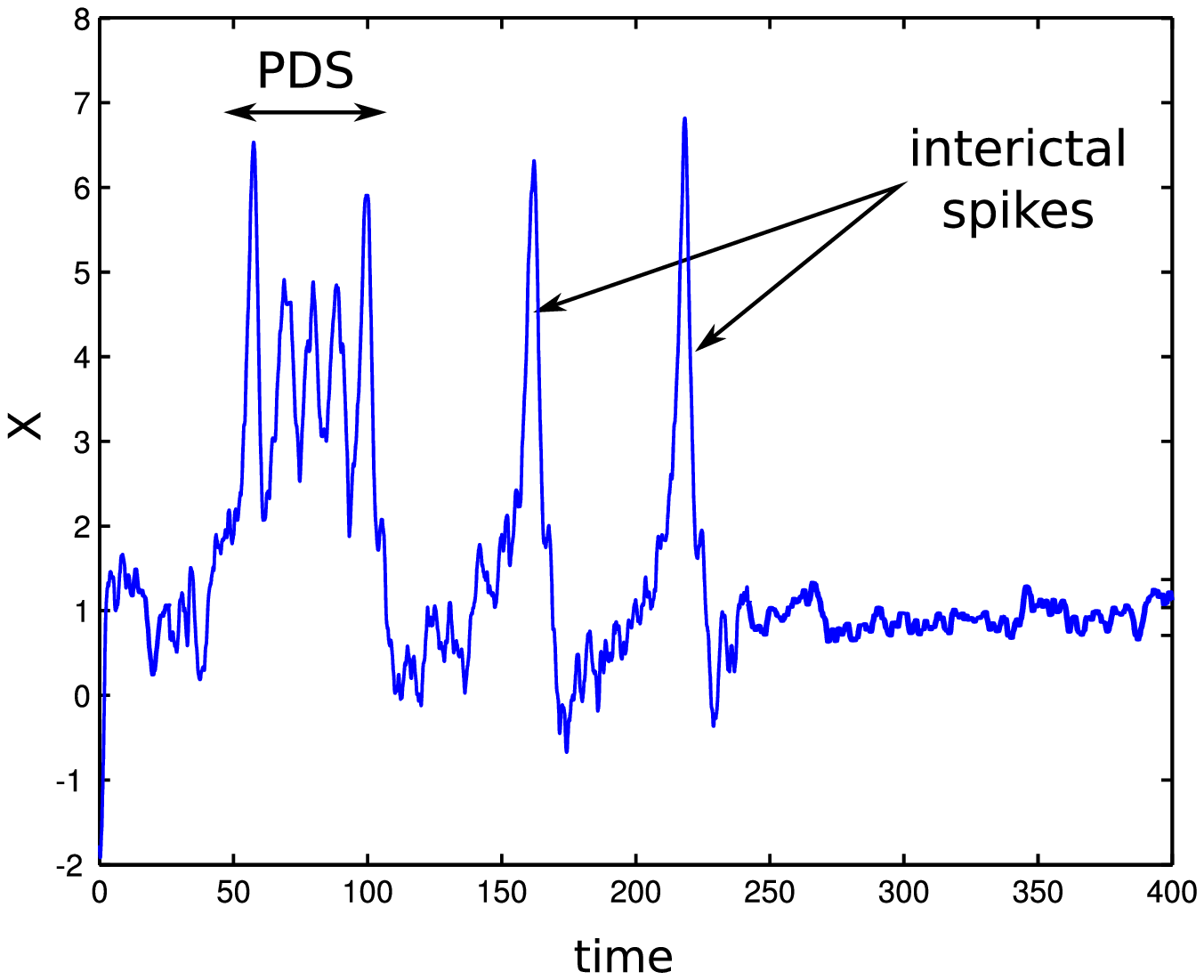}\label{fig:interictalSpikes}}
 \end{center}
\caption{Different behaviours in Jansen's reduced model driven by a Brownian noise $P$ of mean $\mu=1.8$ and standard deviation of $.5$. (a) Interictal spikes leading to oscillatory activity. The parameters are those given in \eqref{eq:Jreduced}.  (b) Interictal spikes and paroxismal depolarizing shift (PDS). The parameters are those given in \eqref{eq:Jreduced} except for $j=12.7$.}
\end{figure}

In the cases (C) to (E) of the diagram in figure \ref{fig:FullGlobal}, interictal spikes appear, but do not always correspond to a pre-ictal activity. Indeed, the destabilization of the stable fixed point may lead the system around the epilepsy limit cycle, but since this cycle is not the only stable behavior, the apparition of random spikes will hardly lead to a seizure. Indeed, it appears that the attraction basin of the epileptic cycles is quite small. When destabilizing the fixed point, we are very likely to present a spike, because at the saddle-node bifurcation point is also a point of the saddle-node homoclinic cycle, part of the branch of limit cycles corresponding to epileptic spikes. But the noise will soon make the solution quit the attraction basin of the cycle and return to the other stable state. In this case, isolated random spikes correspond to some kind of transition from rest to oscillatory activity (see figure \ref{fig:interalpha}). In this case we can observe the emergence of paroxismal depolarizing shifts, as described for instance in \cite[chapter 46]{kandel-schwartz-etal:00}.

In the cases (F) and (G), epilepsy is the only possible behavior in a certain range of input parameters. Therefore, interictal spikes really predict a move towards the epileptic seizure, and therefore are an indication of a transition to a seizure. The type of noise input chosen, the standard deviation of this noise and the equilibrium perturbed impact the structure and the properties of the spike train generated, but the qualitative behavior is always the same. When the inputs keep the system reasonably far from the saddle-node homoclinic bifurcation, spikes have a very small probability to occur and in practice are hardly observed. When getting closer to the bifurcation point, isolated spikes appear, see figure \ref{fig:interictalSpikes}. Their occurence depends on the value of the input firing rate and the noise standard deviation. When increasing further the mean input and getting closer to the bifurcation points, spikes become almost rhythmic and several superimposed action potentials can appear, which are known to be related with paroxismal depolarizing shifts (PDS), see figure \ref{fig:interictalSpikes}. When increasing again the stimulation, we are in an epileptic crisis, and we observe rhythmic high-amplitude epileptic spikes. 

\section{A possibly pathological scheme}
\begin{figure}
 \begin{center}
  \includegraphics[width=\textwidth]{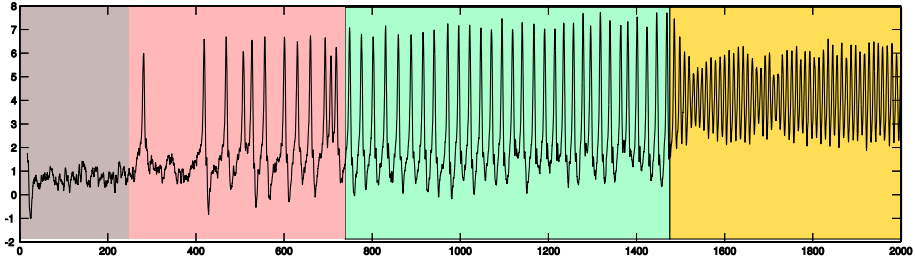}
 \end{center}
 \caption{Seizure in the pathological scheme. We consider an small hyperinnervation corresponding to $j=12.7$. The inputs are considered to be random, with a slowly increasing mean $\mu = 1.5 + 10^{-3}t$ and constant standard deviation of $0.4$. The purple zone corresponds to normal activity, the pink zone to the onset of the seizure: we observe randomly distributed spikes and paroxismal depolarizing shifts. The blue zone corresponds to the seizure itself: we observe rhythmic activity with high amplitude oscillations. The orange zone corresponds to the end of the seizure, with a return to alpha activity. We observe that the signal naturally  waxes and wanes. }
\label{fig:seizure}
\end{figure}

It has been observed by many experimentalists that an inherited epilepsy could be linked with an hyperinnervation (see e.g. \cite{noebels:96,babb-pretorius-etal:89,munoz-mendez-etal:06}). This would correspond in our models to an increase of the value of the total connectivity $j$ (as well as possibly some changes in the proportions coefficients $\alpha_i$). As a first approximation, we observe that small but pathological increases of the total connectivity coefficient $j$ can lead to the ``epileptic parameter zones'' (F) and (G). Let us consider that in the resting state the column receives inputs corresponding to a normal (fixed point) behavior, and that some kind of stimulus results in a slow increase of the mean level of inputs. In that case, the column first behaves normally, then goes through a slow transition to epileptic activity during which interictal spikes and paroxismal depolarizing shifts are likely to occur (see figure \ref{fig:seizure}).

\section{Mathematical interpretation of EEG recordings}
What makes a good neural mass model? This question is very far reaching and the answer clearly depends on the phenomena one wants to account for. From the mathematical point of view, we have been able in the different cases studied in this paper to interpret oscillatory behaviors via a bifurcation analysis of cycles. Lennaert van Veen and David Liley have studied another model of EEG signal in \cite{veen-liley:06}. They observed a bifurcation diagram very similar to the one of Jansen and Rit's model but in their case the saddle connection S corresponds to a Shilnikov Saddle-Node bifurcation \cite{kuznetsov:98}. 

The mathematical approach based on bifurcation theory that we developed in order to understand the behavior of neural masses was already successfully applied to single neuron modeling. In this domain, bifurcations have been related to the excitability properties of the underlying dynamical system since the pioneering paper of Rinzel and Ermentrout \cite{rinzel-ermentrout:89} in the late 80s, and the canonical models studied by Ermentrout and Kopell. This approach was then generalized to the study of the excitability properties of the neurons (see, e.g., \cite{izhikevich:07, izhikevich:00} for reviews). This approach naturally led to introduce such phenomenogical models of neurons  as the nonlinear twodimensional integrate-and-fire neuron \cite{izhikevich:03, brette-gerstner:05, touboul:08, touboul-brette:08b}. In the case of neural masses there are two outstanding questions that we need to answer: first, can we relate EEG behaviors with bifurcations of the underlying dynamical system? and second, can we propose a ``minimal model'' that is able to reproduce phenomenologically EEG signals? We next address these two questions.

\subsection{Bifurcations and oscillatory behaviors}
EEG signals are characterized by different phenomena which can be linked to mathematical objects:
\begin{itemize}
\item The normal activity is related to the existence of a stable fixed point. 
\item Oscillatory activity, whatever its frequency band, is related to the existence of a family of limit cycles.
\item Fast activity can be either related to the existence of a family of limit cycles having a high frequency, or to the destabilization through noise of a stable equilibrium whose Jacobian matrix has complex eigenvalues.
\item Epileptic spikes correspond to low frequency, large amplitude, oscillations. In general these oscillations appear suddenly when varying the parameters which corresponds to the advent of a seizure. Low frequency oscillations can be related to the existence of homoclinic orbits or homoclinic bifurcations which are for instance caused by the presence of a Bogdanov-Takens bifurcation. 
The simplest mathematical object accounting for a sudden appearance of large amplitude oscillations is the fold bifurcation of limit cycles.
\item In some cases it can be interesting for the model to feature a bistability between an epileptic activity and a rhythmic activity. The simplest bifurcation accounting for this coexistence is the cusp of limit cycles. 
\end{itemize}
Note that all these mathematical objects are present in the two models that are studied in this paper.
\subsection{A minimal model?}
A suitable model for reproducing EEG activity could therefore contain the birfurcations just cited. In order to obtain a simple model of this kind, we build on the observation that both Jansen and Rit's and Wendling and Chauvel's models feature a codimension three bifurcation with the universal unfolding of the degenerate Bogdanov-Takens bifurcation. This bifurcation corresponds to the cusp case of the singularity of planar vector fields with nillpotent Jacobian matrix in the cusp case. This bifurcation was studied by Dumortier, Roussarie and Sotomayor in \cite{dumortier-roussarie-etal:87}, and then complemented with the elliptic, saddle and focus cases by Dumortier and collaborators in \cite{dumortier-roussarie-etal:91}. 

The center manifold of this bifurcation is two-dimensional and the normal form of this bifurcation reads on the center manifold:
\begin{equation}\label{eq:DBTNormForm}
\begin{cases}
\dot{x} &= y \\
\dot{y} &= x^2 + \alpha + y (\beta+\gamma \,x \pm x^3)
\end{cases}
\end{equation}

This bifurcation is generic, as proved in \cite{dumortier-roussarie-etal:91}, therefore any system that can be put in this form with higher order terms will present the same bifurcation diagram. Dumortier and collaborators proved also that locally around the Degenerate Bogdanov-Takens point (DBT), the system presents a curve of cusp bifurcations, Bogdanov-Takens bifurcations, two saddle-node bifurcation curves, a Hopf bifurcation curve, a curve of fold bifurcation of limit cycles, a Bautin (generalized Hopf) bifurcation curve and a saddle-node homoclinic curve. 

Hence, when varying the input rate and the connectivity strength, this model features exactly the same bifurcations as Wendling and Chauvel's with their original parameters. 

\section{Conclusion}
We have presented in this article some results arising from the study of bifurcations in neural mass models in order to understand the origin of some brain rhythms and epileptic seizures. The study of the different cycles of the system is of great interest. The effect of the different parameters on the model behaviour can be understood in the light of these bifurcations. The effect of noise on such systems, and the relations between the dynamics of these deterministic models and the related mean field equation is a noble endeavor which will probably provide great insights on the effect of noise in the apparition of seizures and in the rhythms of the brain.

\appendix
\section{A symbolic and numerical bifurcation algorithm}\label{appendix:Numerbif}
In this appendix we present an algorithm to compute formally or numerically local bifurcations for vector fields implemented in Maple\textregistered. This algorithm is based on the closed form formulations for genericity and transversality conditions given in textbooks such as \cite{guckenheimer-holmes:83, kuznetsov:98}. For our numerical experiments, we implemented a precise and efficient solver of equations based on dichotomy. This algorithm controls the precision of the solution we are searching for, and is way faster than the native \texttt{fsolve} Maple application. It is the one we used in this article to study the Jansen and Rit's and Wendling and Chauvel's models bifurcations.



\subsection{Saddle-node bifurcation manifold}
We recall that given a dynamical system of the type $\dot x = f(\mu,x)$ for $x \in \R^n$ and $\mu \in \R$, the systems undergoes a saddle-node bifurcation at the equilibrium $x=x_0, \mu=\mu_0$ if and only if (see e.g. \cite[Theorem 3.4.1.]{guckenheimer-holmes:83}): 

\renewcommand{\theenumi}{(SN\arabic{enumi})}
\begin{enumerate}
 \item\label{case:SN1} $D_x f(\mu_0, x_0)$ has a simple $0$ eigenvalue. Denote by $v$ (resp $w$) the right (resp. left) eigenvector.
 \item\label{case:SN2} $\langle w, \partial f/\partial \mu) (x_0, \mu_0) \rangle  \neq 0$
 \item\label{case:SN3} $\langle w, (D_x^2 f(\mu_0,x_0) )(v,v) \rangle  \neq 0$
\end{enumerate}

The algorithm we use to identify the saddle-node bifurcation manifold is a straightforward application of this theorem, and consists in:
\begin{itemize}
 \item solving the implicit equation $\det(D_x f(\mu_0, x_0)) = 0$ on the fixed-points manifold. We obtain the equilibria where the Jacobian determinant vanishes, i.e. where there is a null eigenvalue. 
 \item We then compute the left and right eigenvectors of the Jacobian matrix at these points:
 \begin{itemize}
  \item check that the dimension of the eigenspace related to the eigenvalue $0$ is 1
  \item check conditions \ref{case:SN2} and \ref{case:SN3}. 
 \end{itemize}
\end{itemize}

At the points where the determinant of the Jacobian matrix vanishes and the conditions are not satisfied, we are led to consider the two classical cases of the cusp and the Bogdanov-Takens bifurcations.

\subsection{Cusp bifurcation}
At a cusp bifurcation point, the Jacobian matrix of the system has a null eigenvalue and the first coefficient of the normal form (given by condition \ref{case:SN3}) vanishes. In that case, under the differential transversality and genericity conditions of \cite[lemma 8.1]{kuznetsov:98}, a smooth change of coordinates will put locally the system in the normal form of the cusp bifurcation. Our algorithm numerically checks the two conditions at these singular points.

\subsection{Bogdanov-Takens bifurcation}
If along the saddle node manifold a second eigenvalue vanishes, under the genericity and transversality conditions of \cite[chapter 7.3]{guckenheimer-holmes:83}, a smooth change of coordinates will locally put the system in the normal form of the Bogdanov-Takens bifurcation. These conditions are also numerically checked by our algorithm.

\subsection{Andronov-Hopf bifurcation manifold}
Changes in the stability of fixed points can also occur via Andronov-Hopf bifurcations. In this case, the real part of an eigenvalue crosses $0$ but not its imaginary part. Theoretically, the dynamical system $\dot x = f(x,\mu)$ undergoes a Hopf bifurcation at the point $x=x_0, \mu=\mu_0$ if and only if (see \cite[Theorem 3.4.2]{guckenheimer-holmes:83}):
 
\renewcommand{\theenumi}{(H\arabic{enumi})}
\begin{enumerate}
 \item\label{case:H1} $D_x f(x_0, \mu_0)$ has a simple pair of pure imaginary eigenvalues and no other eigenvalues with zero real part. Denote by $\lambda(\mu)$ the eigenvalue which is purely imaginary at $\mu_0$.
 \item\label{case:H2} $l_1(x_0, \mu_0)\neq 0$ where $l_1$ is the first Lyapunov exponent.
 \item\label{case:H3} $\frac{d}{d\mu} (Re(\lambda(\mu)))\vert_{\mu=\mu_0} = d \neq 0$
\end{enumerate}

In this case, checking condition \ref{case:H1} is not as easy as condition \ref{case:SN1}. Different methods are available in order to compute Hopf bifurcation points (see \cite{guckenheimer-myers-etal:96,guckenheimer-myers:96}). Our algorithm is based on computing the bialternate product of the Jacobian matrix of the system $2\,J(j,P)\odot Id$ where $Id$ is the identity matrix. It is known that the determinant of this matrix vanishes if and only if the Jacobian matrix has two opposed eigenvalues (its eigenvalues are $\lambda_i + \lambda_j$ where $\lambda_i$ and $\lambda_j$ are eigenvalues of the initial matrix $J(j,P)$). 
Therefore at these points where the bialternate product vanishes we have to check that there exists a purely imaginary eigenvalue to avoid cases where the system has two real opposed eigenvalues. Even with this step, the bialternate product method is much more efficient than other methods based on the characteristic polynomial such as Kubicek's method \cite{kubicek:80}. 

\subsection{Bautin bifurcation}
If along the line of the Hopf bifurcations the first Lyapunov coefficient vanishes and a subcritical Hopf bifurcation becomes supercritical when varying the parameters, under the differential conditions of \cite[theorem 8.2]{kuznetsov:98}, the system undergoes a Bautin bifurcation. At the points where the first Lyapunov coefficient vanishes, the algorithm numerically checks these conditions.

\bibliographystyle{plain}

\end{document}